\documentclass[12pt]{amsart}
\usepackage{amssymb}
\usepackage{amsthm}
\usepackage{euscript}
\newcommand{\ip}[2]{\langle #1,#2\rangle}

\theoremstyle{plain}
\newtheorem{lem}[subsection]{Lemma}

\newtheorem{thm}[subsection]{Theorem}
\newtheorem{prop}[subsection]{Proposition}

\newtheorem{cor}[subsection]{Corollary}

\begin{document}
\title[Spectrum in  topological algebras and  old theorems]{Spectrum in alternative topological algebras and a new look at old theorems}
\author{Bamdad R. Yahaghi}

\address{Department of Mathematics, Faculty of Sciences, Golestan University, Gorgan 19395-5746, Iran}
\email{ bamdad5@hotmail.com,  bbaammddaadd55@gmail.com}


\keywords{  nonassociative algebras,  (classical) division algebras, flexible/power-associative/alternative algebras,   topological/normed/C*-algebras,  quadratic/locally complex  algebras, bilinear forms, spectrum, hyperinvariant linear manifolds,  Frobenius/Hurwitz/Zorn/Gelfand-Mazur theorems}
\subjclass[2010]{
46H70, 17A35, 17A45, 	17A75,  13J30, 46H35, 17D99}

\bibliographystyle{plain}

\begin{abstract}
In this paper, we consider real and complex algebras as well as algebras over general fields. In Section 2,  we revisit and prove several results on (quadratic) algebras over general fields.  As an example, we demonstrate that a quadratic algebra over a field of characteristic not $2$ is flexible if and only if it is proper--a concept introduced in this paper. In Section 3, we show how to develop the spectral theory in the context of complex (resp. real) one-sided alternative topological algebras.   Along the way,  spectral theory is used to prove several topological counterparts of the well-known theorems of Frobenius, Hurwitz, Gelfand-Mazur, and Zorn. This is done, for example, in the context of left (resp. right) alternative topological algebras whose duals separate their elements. In Section 4, we consider real and complex algebras in various topological settings and reconsider and prove several results. For instance, it is shown that given a $ 1 < k \in \mathbb{N}$, on any locally complex algebra,  there exists a unique nonzero vector space norm, say, $\|.\|$,  satisfying the identity $\|a^k\| = \|a\|^k$  on the algebra.   In Section 5, among other things, we revisit and slightly strengthen the celebrated theorems of Frobenius,  Zorn, Gelfand-Mazur, and Hurwitz, and also give slight extensions of their topological counterparts, e.g., theorems of Albert, Kaplansky, and Urbanik-Wright to name a few, in several settings.  

\end{abstract}

\maketitle 

\bigskip

\begin{section}
{\bf Introduction}
\end{section}

\bigskip

A classical theorem of Frobenius \cite{Fro} asserts that {\it up to isomorphisms of real algebras,  the only finite-dimensional associative real division algebras are $\mathbb{R}$, $\mathbb{C}$, and $\mathbb{H}$};  see \cite[page 125]{A}, \cite[Theorem 2.1]{BSh}, and \cite{P}.  In 1930,  Zorn \cite{Zo}  extended Frobenius' theorem to alternative real division algebras as follows. {\it Up to isomorphisms of real algebras, the only finite-dimensional alternative real division algebras are $\mathbb{R}$, $\mathbb{C}$, $\mathbb{H}$, and $\mathbb{O}$}; see \cite[2.5.29]{GP} and \cite[page 262]{EHHKMNPR}. The Gelfand-Mazur theorem states that {\it up to isomorphisms of complex (resp. real) algebras, $\mathbb{C}$ is (resp. $\mathbb{R}$, $\mathbb{C}$, and $\mathbb{H}$ are) the only associative  normed division complex algebra (resp. real  algebras)}; see \cite[Section 8.4]{EHHKMNPR}. Finally, Hurwitz's classical theorem \cite{Hu} asserts that {\it   up to isometric isomorphisms of real algebras, the only finite-dimensional unital absolute-valued real algebras whose absolute values come from inner products are $\mathbb{R}$, $\mathbb{C}$, $\mathbb{H}$, and $\mathbb{O}$.} 
These classical theorems and their extensions have a vast history. For a detailed and comprehensive account of these theorems and the existing literature on them and their extensions in various algebraic and topological settings, we refer the reader to Sections 2.5.4, 2.6.3, and 2.7.4 of \cite{GP}. We also refer the reader to Chapters 8, 9, 10, and 11 of \cite{EHHKMNPR}. In writing parts of this paper, we were  influenced by the reference \cite{GP}.

In this paper, we consider real and complex algebras in both algebraic and topological settings, and also consider algebras over general fields. We demonstrate how to develop spectral theory in the setting of left (resp. right)  alternative topological algebras, particularly normed algebras. Later,  spectral theory is used to prove several  counterparts of the theorems of Gelfand-Mazur, Frobenius, Hurwitz, and Zorn, e.g.,  in the context of  left (resp. right)  alternative topological  algebras whose topological duals separate their elements. We revisit and prove several results on (quadratic) algebras over general fields and on real and complex algebras in various topological contexts. Additionally, among other results, we revisit and slightly strengthen these classical theorems and present slight extensions of their topological analogues in the settings already considered by other authors.  It is worth mentioning that the references of this paper are by no means complete, and perhaps very incomplete. Those references that are not cited in the text, of course, are listed for their general relevance to the subjects touched on in this paper.

\bigskip

Let us set the stage by recalling some standard notations and definitions. Customarily,  we use the symbols $\mathbb{R}$, $\mathbb{C}$,  $\mathbb{H}$, and $\mathbb{O}$ to, respectively, denote the sets, in fact the real division algebras, of real, complex, quaternion, and octonion numbers. We call these sets the {\it numerical division algebras}. In this paper,  we mostly consider  algebras over $\mathbb{F}$, by which, unless otherwise  strictly stated, we mean  $\mathbb{R}$ or $\mathbb{C}$.

A vector space  $\mathbb{A}$ over a field $F$ together with a multiplication coming from an $F$-bilinear function on  $\mathbb{A}$  is said to be an {\it $F$-algebra}. Throughout,   by an algebra we mean an arbitrary nonzero algebra not necessarily associative or commutative or unital. For $ S \subseteq \mathbb{A}$, we use the symbol $ \langle S \rangle_F$ to denote the $F$-linear subspace spanned by the subset $S$.   An element $ a$ of an $F$-algebra  $ \mathbb{A}$ is said to be {\it a flexible element} if  
 $(ab)a  =  a(ba)$ for all $b \in  \mathbb{A}$. An $F$-algebra  $\mathbb{A}$ is said to be {\it flexible} if it consists of  flexible elements. An element $ a$ of an $F$-algebra  $ \mathbb{A}$ is said to be {\it a left (resp. right) alternative element} if  $ a(ab) = (aa)b $ (resp. $ b(aa) = (ba)a$) for all $ b \in  \mathbb{A}$.  The $F$-algebra  $\mathbb{A}$ is said to be {\it left (resp. right) alternative} if  $ a(ab) = (aa)b $ (resp. $ a(bb) = (ab)b$) for all $ a, b \in  \mathbb{A}$, which are respectively called the left and the right alternative identities. The $F$-algebra  $\mathbb{A}$ is said to be {\it alternative} if it is both left and right alternative;   $\mathbb{A}$ is said to be {\it one-sided alternative} if it is  left or right alternative. By a theorem of E. Artin, \cite[Theorem 3.1]{Sch}, {\it an $F$-algebra  $\mathbb{A}$ is alternative if and only if for every $a, b \in  \mathbb{A}$, the $F$-subalgebra generated by the elements $ a$ and $b$ is associative.}   Also, by \cite[Lemma 2]{A2} or by setting $c=a$ in Proposition \ref{1.3}(i), {\it a  left (resp. right) alternative $F$-algebra  $\mathbb{A}$ is  alternative if and only if it is  flexible.}  An $F$-algebra  $\mathbb{A}$ is said to be {\it power-associative} if, for every $a \in  \mathbb{A}$, the $F$-subalgebra generated by the element $ a$ is associative. By a theorem of Albert, \cite[Lemma 1]{A2}, {\it one-sided alternative $F$-algebras are power-associative whenever  the characteristic of $F$ is not $2$.} If the characteristic of the ground field $F$ is not $2$, {\it the Jordan product} of   $ \mathbb{A}$, denoted by $\circ$, is defined by $ a \circ b = \frac{1}{2}(ab + ba)$. Following \cite{GP}, we use the symbol $ \mathbb{A}^{{\rm sym}}$ to denote the algebra $ \mathbb{A}$ together with its vector space operations and its Jordan product.  Clearly,  {\it the algebra $ \mathbb{A}^{{\rm sym}}$ is   commutative, and if $ \mathbb{A}$ is power-associative, then so is $ \mathbb{A}^{{\rm sym}}$. }

As is usual,  for an algebra $\mathbb{A}$ over a field $F $, we use the symbols $\mathcal{L}(\mathbb{A})$  (resp. $ \mathbb{A}' := \mathcal{L}(\mathbb{A}, F) $) to denote the set of all $F$-linear operators (resp. functionals)  from  $\mathbb{A}$  into $\mathbb{A}$ (resp.  $F$); $ \mathbb{A}' $ is called the {\it algebraic dual} of the algebra $\mathbb{A}$.  The algebra  $\mathbb{A}$ is said to be {\it unital or to have an identity} if its multiplication has an identity element, denoted by $1 \not= 0$, where $0$ denotes the identity element of the addition operation of the algebra. An element $ e \in \mathbb{A}$ is said to be a {\it left (resp. right) identity element} if $ ea = a $ (resp. $ae = a$) for all $ a \in \mathbb{A}$. The algebra  $\mathbb{A}$ is called {\it  associative (resp. commutative)}  if its multiplication is associative (resp. commutative).  A nonzero element $ a \in  \mathbb{A}$ is said to be {\it invertible} if there exists a unique element of the algebra, denoted by $a^{-1}$, satisfying the relations $ a a^{-1} = a^{-1} a = 1$. {\it Left (resp. right) invertible elements} of an algebra are naturally defined. The symbols  $\mathbb{A}^{-1}$ (resp. $\mathbb{A}_l^{-1}$, $\mathbb{A}_r^{-1}$) are  used to denote the set of all invertible ( resp. left invertible, right invertible) elements of the algebra  $\mathbb{A}$.   Note that if an algebra  $\mathbb{A}$ is associative, then the uniqueness of the inverse element is a redundant hypothesis in the definition of the invertible elements of  $\mathbb{A}$. It turns out that the same holds for left (resp. right) alternative algebras; see Proposition \ref{1.3}(iv) below.

An algebra $\mathbb{A}$ is said to be {\it algebraic} if the subalgebra generated by any element of it is finite-dimensional. A unital algebraic algebra over a field $F$ is said to be {\it locally a field extension of}  $F$ if the unital subalgebra generated by any element of it is a field extension of $F$, equivalently, the algebra is power-associative and that the subalgebra generated by any element of it has no nontrivial zero-divisors. A unital algebra $\mathbb{A}$ over a field $F$  is said to be {\it quadratic} if the unital subalgebra generated by any element of it is at most two-dimensional.  In this paper, for the most part, we deal with quadratic algebras that are locally field extensions of their ground fields. The reader should beware that  there is no universal consensus on the definition of quadratic algebras. For instance, according to \cite[page 764]{A3}, for a given field $F$, a unital $F$-algebra is said to be quadratic if the $F$-algebra generated by any nonscalar element of it is a quadratic field over $F$, and hence in this sense real quadratic algebras, following \cite{BSS}, are what we will call locally complex algebras. On the other hand,  what we call quadratic algebras according to \cite[page 322]{A2} are algebras of degree two.  It is readily checked that quadratic algebras are power-associative.  By a theorem of A. Albert, \cite[Theorem 1]{A2}, {\it every quadratic left (resp. right) alternative algebra over a field of characteristic not $2$  is alternative}.   {\it If $ \mathbb{A}$ is quadratic, then so is $ \mathbb{A}^{{\rm sym}}$.}  An $F$-algebra $ \mathbb{A}$ is said to be {\it quadratic relative to a left (resp. right)  identity element $ e \in \mathbb{A}$} if   for any  $a \in \mathbb{A}$, there are $ r , s \in F$ such that $ a^2 = r a + s e$. Given an algebra  $\mathbb{A}$ over a field $F$,  by a {\it bilinear form on} $\mathbb{A}$, we mean  an $F$-bilinear function $ \langle . , . \rangle : \mathbb{A} \times \mathbb{A}  \longrightarrow F$, in which case the function $ n : \mathbb{A} \longrightarrow F$, defined by $ n(a) = \langle a , a \rangle$ ($a \in  \mathbb{A}$) is called {\it the quadratic form induced by the bilinear form} $ \langle . , . \rangle$. A bilinear form $ \langle . , . \rangle$ on $\mathbb{A}$ is said to be {\it symmetric} if $ \langle a , b\rangle = \langle b ,  a\rangle  $ for all  $a , b\in  \mathbb{A}$. A symmetric bilinear form $ \langle . , . \rangle$ on $\mathbb{A}$ is said to be {\it nondegenerate} if given any $ a \in \mathbb{A}$, $ \langle a , b\rangle = 0 $ for all $ b\in  \mathbb{A}$ implies $ a = 0$.  There is a delicate relation between quadratic algebras and nondegenerate  symmetric bilinear forms via the so-called notion of composition. A bilinear form $ \langle . , . \rangle$ on an algebra over a field $F$ is said to permit {\it composition} if its quadratic form is multiplicative, namely, $ n (ab) = n(a) n(b)$ for all  $a , b\in  \mathbb{A}$.

A real quadratic algebra is said to be {\it locally complex} if  it is locally a field extension of reals, equivalently, every nonscalar element of the algebra is contained in a copy of the complex numbers within the algebra. Note that the algebra of real numbers is trivially  locally complex. Locally complex algebras  were introduced and studied in \cite{BSS}. For a real quadratic algebra $\mathbb{A}$, {\it the real part of} an element $ a \in \mathbb{A}$, denoted by $ {\rm Re}(a)$, is defined to be the unique real number such that $ a = {\rm Re}(a) 1$ or   $ a^2 = 2{\rm Re}(a) a + s 1$ for a unique  $ s \in   \mathbb{R}$ depending on whether or not  $ a \in \mathbb{R}1$.

An element $ a \in  \mathbb{A}$ is said to be a {\it left (resp. right) zero-divisor} if $ ab = 0$  (resp. $ ba = 0$) for a nonzero element $b \in \mathbb{A}$. The element $ a \in  \mathbb{A}$ is said to be a {\it zero-divisor} if it is a left or a right zero-divisor.  The element $ a \in  \mathbb{A}$ is said to be a {\it joint zero-divisor} if there exists a nonzero $ b \in \mathbb{A}$ such that $ ab = ba = 0$. {\it The trivial (joint) zero-divisor}, by definition, is zero.  We use the symbols $ {\rm j.z.d.} (\mathbb{A}) $ (resp. $ {\rm l.z.d.} (\mathbb{A}) $, $ {\rm r.z.d.} (\mathbb{A}) $) to denote the set of all joint (resp.  left, right) zero-divisors  of the algebra $\mathbb{A}$. It is plain that {\it any quadratic real algebra with no nontrivial joint zero-divisors is locally complex. } Also, it is quite straightforward to check that {\it a quadratic real algebra is locally complex if and only if the unital subalgebra generated by any element of the algebra has  no nontrivial joint zero-divisors.}

For $ a, b , c \in \mathbb{A}$, by definition 
$$ [a, b ] := ab - ba, \ \ \ [a, b , c] :=(ab)c - a(b c),$$
are called {\it the commutator of the elements} $a, b$  and  {\it the associator of the  elements} $ a, b , c$, respectively. It is plain that the commutator and the associator are bilinear and trilinear functions, respectively. An $ a \in \mathbb{A}$ is said to be a {\it central element} if $ [ a, \mathbb{A} ] :=  \{ [ a , b ] : b \in \mathbb{A} \} = \{0\}$. An $ a \in \mathbb{A}$ is said to be a {\it nuclear element} if $ [ a, \mathbb{A}, \mathbb{A} ] = [  \mathbb{A}, a,  \mathbb{A} ]  =  [  \mathbb{A},  \mathbb{A}, a ]  =  \{0\}$. The symbols $Z(\mathbb{A}) $ (resp. $N(\mathbb{A} ) $) are used to denote the {\it centre (resp. nucleus)} of $  \mathbb{A}$, namely, the set of all central (resp. nuclear) elements of the algebra   $  \mathbb{A}$. Note that the subalgebra generated by any nuclear element of an arbitrary algebra is both commutative and associative.

By a {\it topological real (resp. complex) algebra}, we mean a   real (resp. complex) algebra  $\mathbb{A}$ together with a Hausdorff topology with the property that the operations addition and scalar product of the algebra $A$ are continuous, the  multiplication is separately continuous, and that the inversion, defined on $\mathbb{A}^{-1}$ the set of all invertible elements of $A$, is continuous.  For a topological algebra $\mathbb{A}$ over $\mathbb{F}$,  we use the symbols $\mathcal{L}_c(\mathbb{A})$  (resp. $ \mathbb{A}^* := \mathcal{L}_c(\mathbb{A}, \mathbb{F}) $) to denote the set of all continuous $\mathbb{F}$-linear operators (resp. functionals)  from  $\mathbb{A}$  into $\mathbb{A}$ (resp.  $\mathbb{F}$);  $ \mathbb{A}^* $ is called the {\it (topological) dual} of the algebra $\mathbb{A}$. 
The topological algebra  $\mathbb{A}$ is said to be {\it locally convex} if its  topology is generated by a family   $\{s_j\}_{j \in J}$ of   real (resp. complex) vector space seminorms. A locally convex topological algebra is said to be a {\it Fr\'echet algebra} if its topology comes from a complete invariant metric.   Given an $M > 0$, a vector space norm $\| .\|$ on a vector space $X$  over $\mathbb{F}$ is said to be {\it $M$-submultiplicative} (resp. {\it $M$-supermultiplicative}) if  $ \|ab\| \leq M \|a \| \| b\|$ (resp.  $ M \|a \| \| b\|  \leq \| ab \| $) for all $ a,  b \in \mathbb{A}$. {\it Submultiplicative} (resp. {\it supermultiplicative}) {\it norms}, by definition, are {\it $1$-submultiplicative} (resp. {\it $1$-supermultiplicative}). {\it Multiplicative norms}, by definition, are both sub and supermultiplicative. Recall that a vector space (semi)norm is said to be an {\it algebra (semi)norm} if it is submultiplicative. A multiplicative algebra norm is called an {\it absolute value} and an algebra equipped with an absolute value is called an {\it absolute-valued algebra}. A submultiplicative algebra norm $\|.\|$ on an algebra $\mathbb{A}$  is called a {\it nearly absolute value} if   $\|.\|$ is $c$-supermultiplicative for some $ c > 0$, called  {\it a constant} of the nearly absolute value of $\mathbb{A}$, i.e.,  $ c \|a \| \| b\|  \leq \| ab \| $ for all $ a,  b \in \mathbb{A}$; clearly, $ c \leq 1$ if $\mathbb{A}$ is unital. An algebra equipped with a nearly absolute value is called a {\it nearly absolute-valued algebra}.  It is easily checked that the completion of any (nearly) absolute-valued algebra is a (nearly) absolute-valued algebra. Also, the algebra norm of a  normed algebra $( \mathbb{A}, \|.\|)$ is a nearly absolute value if and only if $ \inf_{||a||=1, ||b||=1} ||ab|| > 0$, in which case this positive number is called {\it the constant} of the nearly absolute value of the algebra $\mathbb{A}$; furthermore,   $ \inf_{||a||=1, ||b||=1} ||ab|| = \sup \{ c > 0 : c \|a \| \| b\|  \leq \| ab \|  \ \forall \ a,  b \in \mathbb{A} \}$.

 A locally convex unital topological algebra is said to be a {\it continuous inverse algebra} if   $\mathbb{A}^{-1}$  is open and that the inversion function from $\mathbb{A}^{-1}$ into $\mathbb{A}^{-1}$ is continuous.  By a {\it locally multiplicatively convex (lmc) algebra}, {\it locally sm. convex algebra} or  {\it multinormed algebra}, we mean a locally convex vector space whose seminorms are algebra seminorms; see for instance  \cite[Definition 2.1]{F},  \cite[Definition 4.4.1]{B}, \cite[Definition I.2.4]{H}. Clearly, lmc algebras are locally convex algebras but, by \cite[Exercise I.2.22]{H}, the converse does not hold. By a {\it real (resp. complex) normed algebra}, we mean a real (resp. complex)   algebra $\mathbb{A}$,  together with a real (resp. complex) algebra norm on $\mathbb{A}$.    By   the counterpart of   \cite[Theorem 16.12]{Wa}  for left (resp. right) alternative algebras,  in view of Proposition \ref{1.3}(iv),  {\it left (resp. right) alternative normed  algebras, more generally   left (resp. right)  alternative  lmc  algebras, are topological algebras.} 

A complex Banach algebra together with an involution $*$ satisfying the identity $ \| a^* a \| = \|a\|^2 $ is said to be a  {\it $C^*$-algebra}. A {\it real $C^*$-algebra}, by definition, is a closed real $*$-subalgebra of a  ${\rm C}^*$-algebra. By a theorem of A. Rodr\'iguez, \cite[Theorem 3.2.5]{GP}, {\it every unital complex, and hence real,  $C^*$-algebra is alternative. }

An algebra norm $||.||$  of a unital algebra  $\mathbb{A}$ is said to be {\it unital} if $||1||= 1$,  where the first $1$ denotes the identity  element of the algebra  $\mathbb{A}$.  An element $ a \in  \mathbb{A}$ of a normed algebra is said to be a {\it left (resp. right) topological zero-divisor}  if there exist  $b_n \in \mathbb{A}$ with $ ||b_n|| = 1$ ($n \in \mathbb{N}$) such that   $ \lim_n a b_n  = 0$ (resp.  $ \lim_n  b_n a = 0$).  An element $ a \in  \mathbb{A}$  is called a {\it topological zero-divisor} if it is a left or a right topological zero-divisor.  The element $ a \in  \mathbb{A}$ is said to be a {\it joint topological  zero-divisor} if there are $ b_n \in \mathbb{A}$ with $ ||b_n|| = 1$ ($n \in \mathbb{N}$)  such that $\lim_n ab_n =\lim_n b_na = 0$.   Again,  by definition, zero is {\it the trivial  (joint)  topological zero-divisor}.  We use the symbols $ {\rm j.t.z.d.} (\mathbb{A}) $ (resp. $ {\rm l.t.z.d.} (\mathbb{A}) $, $ {\rm r.t.z.d.} (\mathbb{A}) $) to denote the set of all joint (resp. left, right) topological zero-divisors  of the algebra $\mathbb{A}$. 

 Clearly,  (one-sided) zero-divisors are (one-sided) topological zero-divisors. However, the converse of this does not hold in general unless the algebra is finite-dimensional. It is plain that  an element $ a $ of an algebra  $\mathbb{A}$ is not a left (resp. right) divisor of zero if and only if the mapping $L_a: \mathbb{A} \longrightarrow a \mathbb{A} $  defined by $ L_a x = ax$ (resp.  $R_a: \mathbb{A} \longrightarrow \mathbb{A}a $ defined by $R_a x = xa$)  ($ x \in \mathbb{A}$) is one-to-one, equivalently, left invertible as a linear function. We leave it as an exercise to the interested reader to verify that {\it an element $ a $ of a normed algebra  $\mathbb{A}$ is not a left (resp. right) topological divisor of zero if and only if the mapping $L_a: \overline{\mathbb{A}} \longrightarrow a\overline{\mathbb{A}} $ defined by $ L_a x = ax$ (resp.  $R_a: \overline{\mathbb{A}} \longrightarrow \overline{\mathbb{A}}a$  defined by $R_a x = xa$)  ($ x \in \overline{\mathbb{A}} $)  has a continuous left inverse as a linear operator;} here $\overline{\mathbb{A}}$ denotes the completion of the normed algebra $\mathbb{A}$.

An algebra $\mathbb{A}$ is said to be a {\it left (resp. right) division algebra} if the mappings $  L_a : \mathbb{A} \longrightarrow  \mathbb{A}  $ (resp. $  R_a : \mathbb{A} \longrightarrow  \mathbb{A}  $)  are one-to-one and onto for each nonzero element $ a \in \mathbb{A}$.    The algebra $\mathbb{A}$ is said to be a   {\it division algebra} (resp. {\it one-sided division algebra}) if it is both a left and right  (resp. a left or a right)  division algebra. The algebra  $\mathbb{A}$ is said to be a   {\it quasi-division algebra} if for every nonzero element $ a \in \mathbb{A}$, at least one of the mappings $  L_a $ or  $R_a$ is bijective. Finally, the algebra $\mathbb{A}$ is said to be a   {\it classical division algebra} or a {\it division algebra in the classical sense}  if  it is unital and $\mathbb{A}^{-1} = \mathbb{A} \setminus \{0\}$.

It is not difficult to observe  that {\it alternative one-sided division algebras are division algebras in the classical sense.} To this end, let $\mathbb{A}$ be an alternative left  division algebra. First, pick a nonzero $ a \in  \mathbb{A}$. From the surjectivity of $  L_a$, we see that $ ae = a$ for some nonzero element $ e \in \mathbb{A}$. As $\mathbb{A}$ is alternative, we get that $ a e^2 = ae$, from which we obtain $ L_a ( e^2 - e)= 0$. This, in turn, yields $ e^2 = e$ because $L_a $ is one-to-one. Now, for an arbitrary $ b \in  \mathbb{A}$, we can write $e^2 b = e b$ and $ b e^2 = be$. It thus follows from the alternativity of the algebra that  $L_e( eb - b) = 0 = L_{be - b}(e )$, implying that $ eb = b = be$. In other words, $e=1$ is the unital element of the algebra   $\mathbb{A}$. Next, for an arbitrary nonzero $ a \in  \mathbb{A}$, there exists a unique $ b  \in  \mathbb{A}$ such that $L_a (b) =a b = 1$.  Since $\mathbb{A}$ is alternative, this clearly yields $L_{a}(ba - 1) = 0$, from which we obtain $ ba = 1$. This means $b \in  \mathbb{A}$ is the unique inverse of the element $a   \in  \mathbb{A}$. Moreover, {\it alternative  division algebras in the classical sense are division algebras;} see Proposition \ref{1.3}(vi) below.

\bigskip

Naturally, one can define the notion of a t-division algebra as follows.  A topological algebra $\mathbb{A}$ is said to be a {\it left (resp. right) t-division algebra}  if the continuous mappings $  L_a : \mathbb{A} \longrightarrow  \mathbb{A}  $ (resp. $  R_a : \mathbb{A} \longrightarrow  \mathbb{A}  $) are invertible and have continuous inverses for every nonzero element $ a \in \mathbb{A}$. The algebra $\mathbb{A}$ is said to be a {\it t-division algebra} if it is both a left and right t-division algebra.   The algebra  $\mathbb{A}$ is said to be a   {\it quasi-t-division algebra} if for every nonzero element $ a \in \mathbb{A}$, at least one of the mappings $  L_a $ or  $R_a$ is bijective and has a continuous inverse.   

 \bigskip

In what follows, we need the following well-known proposition. It is worth mentioning that part (v) of it, very likely known by the experts, is not taken from anywhere. For the sake of completeness, we  provide a detailed proof of Skornyakov's Identities, taken from \cite[pages 178-179]{Sk}, with a minor typo that occured there corrected in the process; also see    \cite[Lemma 1]{K}.

 \bigskip 

\begin{prop} \label{1.3} 
 Let  $F$ be a field and  $\mathbb A$  be a left (resp. right) alternative  $F$-algebra. The following statements hold.

{\rm (i)} If $ a, b, c \in \mathbb A$, then 
$$ (ab + ba) c = a(bc ) + b(ac), \ \  \Big({\rm resp.} \  c ( ab + ba) = (ca) b + (cb) a \Big),$$
equivalently, 
$$ [a, b, c] = - [ b, a , c] , \ \ \  \Big ({\rm resp.} \   [a, b, c]= -[ a, c, b] \Big). $$
In particular, if $ a, b, c \in \mathbb A$ and $ ab =  -ba$, then 

$$ a (bc) = - b (ac), \ \  \   \Big({\rm resp.}  \ (c a) b = -  (c b) a \Big).$$ 

{\rm (ii)} {\bf (Skornyakov's Identities)}   If $ a, b , c \in  \mathbb A$ and the characteristic of $F$ is not $2$, then 
$$ (a (b a))c = a \big( b (a c) \big)   ,    $$
$$ \Big ({\rm resp.} \       a((bc) b) = \big( (ab) c \big) b. \Big) $$

{\rm (iii)} {\bf (Moufang's Identities)}  If  $\mathbb A$ is alternative, then the following identities, known as the Moufang's identities, hold on  $\mathbb A$ 
$$ (ab a)c = a \big( b (a c)  \big)     , \ \ \   (ab) (ca) = a (bc) a  , \  \ \       a(bc b) = \big( (ab) c \big) b.   $$


{\rm (iv)} Let $ a \in \mathbb A$ and suppose that the characteristic of $F$ is not $2$. Then,   $b \in \mathbb A$ is an inverse of $a$, i.e., $ ab = ba = 1$, if and only if   $L_a L_b = I = L_b L_a $ (resp. $R_a R_b = I = R_b R_a $), if and only if $a(bc) = c = b(ac)$   (resp.  $(cb)a = c = (ca)b$) for all $ c \in \mathbb{A}$. Therefore, $a$ has a unique inverse. 

{\rm (v)}  Let $\mathbb A$ be alternative with characteristic not $2$ and $ a, b, c \in \mathbb A$ with $ a, b \in \mathbb A^{-1}$. Then, $ a(bc) = 1$  (resp. $ (cb)a = 1$) if and only if  $ (ab) c = 1$   (resp. $ c(ba) = 1$).  
Moreover, if $ a, b \in \mathbb A^{-1}$, then $ ab \in \mathbb A^{-1}$ and $ (ab)^{-1} = b^{-1} a^{-1}$.  

{\rm (vi)}  Let $\mathbb A$ be alternative with characteristic not $2$ and $ a \in \mathbb A$. Then, $ a \in  \mathbb A^{-1}$ if and only if $ L_a \in L( \mathbb A)^{-1}$ (resp.  $ R_a \in L( \mathbb A)^{-1}$), in which case, $L_a^{-1} = L_{a^{-1}}$  (resp. $R_a^{-1} = R_{a^{-1}}$).

\end{prop}

\bigskip

\noindent {\bf Proof.} (i) The first assertion is a straightforward conclusion of linearizing the identities $ a(ac) = (aa)c $ and $ c(aa) = (ca)a$ for all $ a, c \in  \mathbb{A}$.  The second assertion is a quick consequence of the first assertion. 

(ii)  
It suffices to prove the assertion for right alternative algebras. The proof for left alternative algebras is done by passing to the corresponding opposite algebras. Let $ a , b, c, d \in \mathbb{A}$ be arbitrary. It is quite straightforward to verify that
$$a[ b , c, d] + [ a, b, c ] d = [ab, c, d] - [ a, bc, d] + [ a, b, cd].$$
This,  in particular, yields 
 $$  c[ a , b, a] + [ c, a, b ] a  = [ca, b, a] - [ c, ab, a]  + [ c, a, ba], $$
$$ c[ a , a, b] + [ c, a, a ] b  =  [ca, a, b] - [ c, a^2, b] + [ c, a, ab].    \eqno  (1)$$ 
Adding the preceding equlatities and using part (i), we obtain 
 $$  [ c, a, b ] a  =  2 [ c, a, ab] - [ c, a^2, b] + [ c, a, ba].   \eqno  (2)  $$
On the other hand, we can write 
\begin{eqnarray*}
- [c, a, b ]a  & = & [ c , b, a] a + 0, \\
& = & [ c , b, a] a  +c [ b, a, a ], \\
 &   = & [cb, a, a] - [ c, ba, a] + [ c, b, a^2]\\
  & =&[ c, a,  ba] - [ c, a^2, b] , 
   \end{eqnarray*}
 implying that  
 $$ - [ c , a, b] a = [ c, a,  ba] - [ c, a^2, b].  \eqno  (3)  $$
This along with $(2)$ and the hypothesis that the characteristic of $F$ is not $2$ yields 
 $$ [ c, a^2, b]  =[ c, a, ab]  +  [ c, a,  ba].    \eqno  (4) $$ 
 From (3) and (4), we get  that 
 $$ [ c, a, ab] =  [ c, a,  b]a,    $$ 
equivalently, 
$$ [ a, bc, b] = - [ a, b,  c]b.    $$ 
Now, expanding and simplifying the above equality entails 
  $$ \big( (ab) c \big) b = a((bc) b) , $$
  as desired. 

(iii) See  \cite[page 28]{Sch} or \cite[Lemma 2.3.60]{GP}.

(iv) See \cite[Lemma on Inverses (ii)]{N}. 

(v)  First, suppose that $ (ab) c = 1$. Thus, $ ( (ab) c) b =  b ( (ab) c)  = b$. By the right Moufang identity, we can write $ a \big((bc)b \big) = b$, and hence, in light of (iv), $ (bc)b = a^{-1}\Big(a \big((bc)b \big) \Big)  = a^{-1} b$ because $\mathbb A$ is left alternative. That is,  $ (bc)b =a^{-1} b$, from which we obtain $ bc = \big( (bc) b\big)b^{-1} = (a^{-1} b) b^{-1} = a^{-1}$. This clearly implies $  a(bc) = 1$, as desired.

For the second assertion, it suffices to show that $ (b^{-1} a^{-1}) (ab) = 1 = (ab) (b^{-1} a^{-1})$. Plainly, $ \big((ab) b^{-1}\big) a^{-1} = 1$. It thus follows from the preceding paragraph that $(ab) ( b^{-1}a^{-1}) = 1$. Next, setting $  z:= (b^{-1} a^{-1}) (ab) $, multiplying both sides by $a^{-1}$ from the right, and using the right Moufang identity, we get that 
$$
 z  a^{-1} = \Big( (b^{-1} a^{-1}) (ab) \Big) a^{-1}  =b^{-1} \Big( \big(a^{-1} (ab) \big) a^{-1} \Big)  = b^{-1} ( b a^{-1}) = a^{-1},
 $$
yielding 
$$z =  (z  a^{-1}) a = a^{-1} a = 1.$$

Finally, suppose   $  a(bc) = 1$. We see that $ bc = a^{-1} \big( a (bc) \big) = a^{-1}$, which yields $c = b^{-1}a^{-1} $. It thus follows from the second assertion that $c  =(ab)^{-1} $, from which we see that $(ab) c  = 1$, as desired. 

(vi) First, suppose  $ a \in  \mathbb A^{-1}$. By (iv), $L_a^{-1} = L_{a^{-1}}$  (resp. $R_a^{-1} = R_{a^{-1}}$). Conversely, if $ L_a \in \mathcal{L}( \mathbb A)^{-1}$ (resp.  $ R_a \in \mathcal{L}( \mathbb A)^{-1}$), then there exists a $ b \in \mathbb A$  such that $ ab = 1$ (resp. $ba = 1$). Then again, we get that $ L_a( ba - 1) = 0$ (resp. $R_a ( ab - 1) = 0$, from which we see that $b$ is an inverse of $a$, and hence the inverse of $a$  by (iv). That is $ a \in  \mathbb A^{-1}$, which is what we want. This completes the proof. 
\hfill \qed

\bigskip 


\bigskip

\begin{section}
{\bf On quadratic algebras over general fields}
\end{section}

\bigskip

We start off with a useful proposition. Part (iii) of the proposition is a generalization of  \cite[Lemma 2.5.5]{GP} to algebras over arbitrary fields and parts (iv)-(vi) of it  are motivated by \cite[Proposition 2.5.38]{GP}. 

\bigskip

\begin{prop} \label{1.1} 
   {\rm (i) }  Let  $F$    be a field and  $\mathbb A$  be an alternative $F$-algebra with 
   $ {\rm l.z.d.} (\mathbb{A})  = \{0\}$ (resp. $ {\rm r.z.d.} (\mathbb{A}) = \{0\}$). If  $\mathbb A$ has a nonzero algebraic element, then  $\mathbb A$   is unital. Therefore, an alternative  algebra with   $ {\rm l.z.d.} (\mathbb{A})  = \{0\}$ (resp. $ {\rm r.z.d.} (\mathbb{A}) = \{0\}$) is  unital if and only if it has a nonzero idempotent, in which case the unital element of the algebra is the only nonzero idempotent element of the algebra, and for $a, b \in \mathbb A$, $ab = 1$ if and only if $ ba =1$.  Moreover, if the algebra $\mathbb A$ is real and algebraic, then it is locally complex  and division in the classical sense.

{\rm (ii)}  Let  $F$    be a field and  $\mathbb A$  be an alternative algebraic $F$-algebra. Then $ {\rm l.z.d.} (\mathbb{A})  = \{0\}$ (resp. $ {\rm r.z.d.} (\mathbb{A}) = \{0\}$) if and only if $\mathbb A$ is a  division algebra. Therefore, if the algebra $\mathbb A$ is real, then  it has  no  nontrivial left (resp. right) zero-divisor if and only if $\mathbb A$ is a quadratic  division  algebra.

{\rm (iii)}  Let  $F$    be a field whose characteristic is not equal to $2$  and  $\mathbb A$  be a power-associative $F$-algebra with $ {\rm j.z.d.} (\mathbb{A})  = \{0\}$. Then, $\mathbb A$ is unital if and only if $\mathbb A$ has a nonzero $F$-algebraic element.

{\rm (iv)}  Let  $F$    be a field whose characteristic is not equal to $2$  and  $\mathbb A$  be a left (resp. right) alternative algebraic $F$-algebra. Then, $\mathbb A$ is a quasi-division algebra,  if and only if $ {\rm j.z.d.} (\mathbb{A})  = \{0\}$,  if and only if it is a  left (resp. right) division  algebra.

{\rm (v)}  Let  $F$  be a field and  $\mathbb A$  be an alternative $F$-algebra. Then, $\mathbb A$ is  a division  algebra, if and only if it is a one-sided division algebra, if and only if it is a  quasi-division algebra,  if and only if it is a division algebra in the classical sense. 

{\rm (vi)}  Let  $F$  be a field whose characteristic is not $2$ and  $\mathbb A$  be a left (resp. right)  alternative $F$-algebra. Then, $\mathbb A$ is  a division  algebra, if and only if it is a one-sided division algebra, if and only if it is a  quasi-division algebra,  if and only if it is a division algebra in the classical sense. 
\end{prop}

\bigskip

\noindent {\bf Proof.} (i) The second assertion is a quick consequence of the first assertion and the observation that for $a, b \in \mathbb A$, $ ab =1$ if and only if $ a(ba - 1) =0$ if and only if $ ba = 1$. This being noted, let $ a \in \mathbb A$ be nonzero and $F$-algebraic. Since $\mathbb A$   has no nontrivial  left (resp. right) zero-divisor, there exists a minimal $ m \in \mathbb{N}$ with $ m > 1$ such that $ a^m + c_{m-1} a^{m-1} + \cdots + c_1 a= 0$, where $ c_i \in F$ ($1 \leq i \leq m-1$). Also, note that $c_1 \not= 0$ because $\mathbb A$  has no  nontrivial left (resp. right) zero-divisor. Let $ e := -\Big( c_1^{-1} a^{m-1} +  c_1^{-1}c_{m-1} a^{m-2} + \cdots + c_1^{-1}c_{2} a \Big)$.  From the minimality of $m$, we get that $ e \not= 0$. Clearly, $ae = a  $, from which   we obtain $a e^2  =ae $, which in turn implies $ e^2  = e $. Now, for an arbitrary $ x \in  \mathbb A$, we can write  $ e^2 x = ex$  and $ x e = x e^2$, which yields     $  e( ex -  x) =0$ and $ (x - xe) e = 0$  because  $\mathbb A$ is alternative.   Thus, $ ex = x = xe$. That is, $e$ is the identity element of $\mathbb A$, as desired.

(ii) The nontrivial part of the assertion is a straightforward consequence of (i) and  the fundamental theorem of algebra for polynomials with real coefficients.

(iii) The nontrivial part of the assertion is a consequence of the proof of (i) and the observation that the counterparts of \cite[Lemmas 2.5.3-2.5.5]{GP} hold for 
 power-associative $F$-algebras with ${\rm ch} (F) \not= 2$; the proofs, adjusted accordingly, go verbatiem.

(iv)  That division algebras are quasi-division  is trivial. It is easily checked that quasi-division algebras have no nontrivial joint zero-divisors. So to complete the proof, it remains to be shown that the left (resp. right) alternative algebraic $F$-algebra  $\mathbb A$ is a division algebra if it has no nontrivial  joint zero-divisors.  Note that $\mathbb A$ is  power-associative by \cite[Lemma 1]{A2}. It thus follows from (iii) that   $\mathbb A$ is unital. This together with the hypothesis that $\mathbb A$ is algebraic implies that every nonzero element $ a \in \mathbb A$ has an inverse $ a^{-1} \in \mathbb A$, which is indeed a polynomial in $a$ with coefficients from $F$.  Now note that, in view of Proposition \ref{1.3}(vi),  $L_a, R_a: \mathbb A \longrightarrow \mathbb A$ are invertible for all nonzero elements  $ a \in \mathbb A$ with their inverses being given by  $L_{a^{-1}}, R_{a^{-1}}: \mathbb A \longrightarrow \mathbb A$, respectively. This completes the proof.

(v) Clearly, in view of Proposition \ref{1.3}(vi), it suffices to show that if  $\mathbb A$  is an alternative quasi-division $F$-algebra, then it  is a division $F$-algebra in the classical sense. To this end, pick  a nonzero  $ a \in \mathbb A$. It follows that $ L_a $ or $R_a$ is invertible. If $ L_a $ is invertible, there exists an $ e \in  \mathbb A$ such that $ ae = a $.   It follows that  $e$   is an idempotent element of $ \mathbb A$, for  $ L_a(e^2) = L_a(e)  $.   Thus, $e^2 = e$, and hence   $e^2 b = e b$ and $ b e^2 = be$  for all $ b \in \mathbb A$. Then again,  $L_e$ or $ R_e$ is invertible and   $L_e( eb - b) = 0 = R_e (be - b)$  for all $ b \in \mathbb A$. Therefore, $e$ is a left or a right identity element of $\mathbb A$. But if $e$ is a left identity element of $\mathbb A$,   then 
$ (be - b)^2 =0 $, implying that $be =b $ for all   $ b \in \mathbb A$, which means $e$ is an identity element of  $\mathbb A$. Likewise, if $e$ is a right identity element of $\mathbb A$, it is also a left identity element, and hence an identity element of  $\mathbb A$. 
If  $R_a$ is invertible, one can similarly get  that  $\mathbb A$ has an identity element, say, $ e \in \mathbb A$. 
In any case,  we have shown that $ e =1$ is the identity element of $\mathbb A$. Now, just as we saw in the observation we made preceding the proposition, we conclude that every nonzero element of $\mathbb A$ has a multiplicative inverse. This completes the proof. 

(vi) Under any of the hypotheses that $\mathbb A$ is  a division  algebra, or  a one-sided division algebra, or a  quasi-division algebra, or a division algebra in the classical sense, $\mathbb A$ does not have a nonzero nilpotent. To see this, we need to use Proposition \ref{1.3}(iv) along with  the fact  that the algebra $\mathbb A$ is power-associative by \cite[Lemma 1]{A2}. It thus follows from \cite[Corollary 1 (Kleinfeld) on page 345]{ZSSS} that $\mathbb A$ is an alternative algebra. Now, the assertion follows from (v). 
\hfill \qed

\bigskip 

\noindent {\bf Remarks.} 
1. Adjusting the proof of part (i) of the proposition, which is quite elementary, we obtain the following. {\it An alternative $F$-algebra $\mathbb A$  with ${\rm j.z.d. }(\mathbb{A}) =\{ 0\}$ is unital if and only if  $\mathbb A$  has a nonzero central $F$-algebraic element.} Compare this and its self-contained proof with part (iii) and its proof, which relies on several theorems. In the same vein, the following holds.  {\it An $F$-algebra $\mathbb A$  with ${\rm j.z.d. }(\mathbb{A}) =\{ 0\}$ is unital if and only if  $\mathbb A$  has a nonzero nuclear $F$-algebraic element.}

2. As pointed out in the proof of part (iii) of the proposition, the counterparts of \cite[Lemmas 2.5.3-2.5.5]{GP} hold for  power-associative $F$-algebras with ${\rm ch} (F) \not= 2$; their proofs are identical to those of   \cite[Lemmas 2.5.3-2.5.5]{GP} except that they should be adjusted naturally.

3. Clearly, one-sided division algebras have no nontrivial joint (resp. one-sided)  zero-divisors. Proposition \ref{1.3}(iv), along with  parts (iv) and (v) of the proposition and  \cite[Corollary 1 (Kleinfeld) on page 345]{ZSSS}, reveals that left (resp. right) alternative topological, e.g., normed,  algebraic algebras with no nontrivial joint zero-divisors are alternative and division in the classical sense, and  hence have no nontrivial joint  (resp. one-sided) topological zero-divisors.

\bigskip

 Part (i) of the following is the counterpart of \cite[Lemma 2.6.29]{GP}  in the setting of algebras over general fields; also see page 49 of \cite{Sch}.  Part (ii) is a Gelfand-Mazur type theorem for algebraic algebras over algebraically closed fields.

\begin{thm} \label{1.2} 
{\rm (i)} Let $F$ be a field with more than two elements and $\mathbb A$ be an $F$-algebra with the property that $0$ is the only element of it whose square is $0$. If every singly generated subalgebra of $\mathbb A$ is at most one-dimensional, then $\mathbb A$ is one-dimensional, and hence isomorphic to $F$. 

{\rm (ii)} 
 Let $\mathbb A$ be an algebraic  algebra over an algebraically closed field $F$  with no nontrivial left (resp. right) zero-divisors. Then, $\mathbb A$  is isomorphic to $ F$. Therefore, any  algebraic   left (resp. right) division algebra  over an algebraically closed field $F$ is isomorphic to $F$.
 
\end{thm}

\bigskip

\noindent {\bf Proof.} 
(i)  It follows from the hypothesis that for  any nonzero element $a \in \mathbb{A}$, there exists a unique nonzero $ \lambda_a \in F$ such that $ a^2 = \lambda_a  a $.    
Clearly,  $ \lambda_{ za } = z \lambda_a $  for all $ z \in F$ and $  a \in  \mathbb A$. Next, we show that $ \lambda_{ a + b} =  \lambda_a +  \lambda_b$ for all  linearly independent elements $  a, b \in  \mathbb A$. Letting  $ z, w \in F$ with $ zw(z - w) \not= 0$, we can write
 \begin{eqnarray*}
(za +  w b)^2  - zw(a  + b)^2  & =& \lambda_{za + w b} (za + w b) - zw\lambda_{a + b} (a  + b),
 \end{eqnarray*}
from which we obtain
$$(z^2 - zw) \lambda_a a + (w^2 - zw)  \lambda_b b = z\big( \lambda_{za + w b} - w\lambda_{a +  b} \big) a + w \big( \lambda_{za + w b} - z \lambda_{ a + b}\big) b. $$
But $  a, b \in  \mathbb A$ are linearly independent. So we get that
 \begin{eqnarray*}
  z\big( \lambda_{za + w b} - w\lambda_{a +  b} \big) & = & (z^2 - zw)  \lambda_a, \\
 w \big( \lambda_{za + w b}  - z \lambda_{ a + b} \big) & = & (w^2 - zw)  \lambda_b,
  \end{eqnarray*}
  which yields 
   \begin{eqnarray*}
  \lambda_{za + w b} - w\lambda_{a +  b}  & = & (z - w)  \lambda_a, \\
   \lambda_{za + w b}  - z \lambda_{ a + b}  & = & (w - z)  \lambda_b.
  \end{eqnarray*}
Subtracting  the above equations, we obtain $ (z - w) \lambda_{a +  b} =  (z-w) (\lambda_a +  \lambda_b)$, and hence $  \lambda_{a +  b} =  \lambda_a +  \lambda_b$,  as desired. Now proceed by way of contradiction and suppose there are linearly independent elements $  a, b \in  \mathbb A$. Letting $ z = \lambda_b $ and $ w = - \lambda_a $, we see that $ \lambda_{za + w b}= z\lambda_a + w\lambda_b=  0 $, which yields $(za +  w b)^2 =0$. Thus  $za +  w b=0$, a contradiction, completing the proof.

(ii)
Plainly, it suffices to prove the assertion for algebras having no nontrivial left zero-divisors. To this end, pick an arbitrary  $a \in  \mathbb A$. Let $ \mathbb{B}$ denote the $F$-algebra generated by $a$, which is finite-dimensional by the hypothesis. We get that the linear mapping $ L_a : \mathbb{B} \longrightarrow \mathbb{B}$ is one-to-one, and hence onto and invertible. Let $ \lambda_a \in F$ be an eigenvalue of the linear operator $ L_a^{-1} L_{a^2} \in  \mathcal{L}(\mathbb{B})$. It follows that there exists a nonzero element $ b \in \mathbb{B}$ such that $ L_a^{-1} L_{a^2}(b)   = \lambda_a b$, from which we obtain $ a^2 b = \lambda_a ab$. This yields $ ( a^2 - \lambda_a a ) b = 0$, and hence $ a^2 = \lambda_a a$, for zero is the only left zero-divisor of $\mathbb A$.   It thus follows from (i) that $\mathbb A$ is one-dimensional, finishing the proof.  The second assertion follows from the first one. This completes the proof. 
\hfill \qed

\bigskip 

\noindent {\bf Remark.}  With $F$ as in  part (ii) of the theorem, if the $F$-algebraic algebra  $\mathbb A$ is assumed to be power-associative, then it is isomorphic to  $ F$ if and only if j.z.d.$(\mathbb{A}) =\{ 0\}$.

\bigskip

Part (ii) of the following theorem is a standard result on quadratic algebras, whose main part to the best of our knowledge is due to L. Dickson; see \cite[pages 63 and 64]{Di}. Its standard proof is presented here for the sake of completeness; see \cite[pages 49 and 50]{Sch}.

\bigskip

\begin{thm} \label{02.3} 
 {\rm (i)} 
 Let $F$ be a field whose characteristic is not $2$ and   $\mathbb{A}$  a flexible quadratic algebra over $F$ relative to a left (resp. right) identity element $e$ of  $\mathbb{A}$ with the property that $ a^2 - f a = 0$ with $ a \in \mathbb{A}$ and $f \in F$ implies $ a = 0$ or $ a = fe $. Then, $e$ is an identity element of  $\mathbb{A}$. In particular,  $\mathbb{A}$ is a quadratic algebra over $F$ that is locally a field extension of $F$.

 {\rm (ii)} 
 Let $F$ be a field whose characteristic is not $2$ and   $\mathbb{A}$  a quadratic algebra over $F$. Then, there exists a unique  symmetric $F$-bilinear form $\langle . , . \rangle :  \mathbb{A} \times \mathbb{A} \longrightarrow  F$, called the bilinear form of $\mathbb{A}$, such that 
$$ a^2 - 2  \langle a , 1 \rangle a + \langle a , a \rangle 1 = 0, $$
for all $ a \in \mathbb{A}$. In fact, 
\begin{eqnarray*}
  \langle a, b \rangle &  = & 2 \langle a , 1 \rangle  \langle b , 1 \rangle -  \frac{1}{2}  \langle ab + ba , 1 \rangle, \\
   &  = & 2 \langle a , 1 \rangle  \langle b , 1 \rangle -  \frac{1}{4}  \langle (a + b)^2 - ( a  - b)^2, 1 \rangle, 
\end{eqnarray*}
for all $ a, b \in \mathbb{A}$.
Moreover,   first, the bilinear form of   $\mathbb{A}$ is  nondegenerate if $0$ is the only element of  $\mathbb{A}$ whose square is $0$,  and, second, $\langle a , a \rangle = 0$ implies $ a = 0$ for all $ a \in \mathbb{A}$ if the quadratic algebra $\mathbb{A}$ is locally a field extension of $F$.

\end{thm}

\bigskip

\noindent {\bf Proof.} (i)  Let $e$ be a left  identity element of  $\mathbb{A}$. The proof should be adjusted accordingly in case $e$ is a right  identity element of  $\mathbb{A}$. 
 Define the functions $ p_1, n: \mathbb{A} \longrightarrow F$ as follows. On $Fe$, for $ a = t e$ with  $t \in F$, set 
$$ p_1 (a) := t, \  n(a) := t^2 ,$$
and for $ a \in \mathbb{A} \setminus Fe$, let $p_1(a)$ and $n(a)$ be the unique elements of $F$ satisfying the relation $ a^2 - 2 p_1(a) a + n(a) e = 0$. Clearly, $p_1$ is $F$-linear on $ Fe$ and $n$  is $F$-homogeneous of degree two on $ Fe$.  First, $p_1$  and $n$  are, respectively, $F$-homogeneous of degrees $1$  and $2$ on $ Fa$ for any  $ a \in \mathbb{A} \setminus Fe$. To see this, given an $ a \in \mathbb{A} \setminus Fe$, for all $ t \in F \setminus \{0\}$, we can write 
\begin{eqnarray*}
t^2a^2  & = & 2t^2 p_1(a) a - t^2n(a) e , \\
t^2a^2  & = & 2 tp_1(ta) a - n(ta) e,
\end{eqnarray*}
from which, we obtain 
\begin{eqnarray*}
0 =  2t \big( t p_1(a) -  p_1(ta) \big)a + \big( n(ta) - t^2n(a) \big) e.
\end{eqnarray*}
This, in turn, yields  $ p_1(ta)=  t p_1(a)$ and that $ n(ta) = t^2n(a) $, as desired. It thus follows that  $p_1$ and $n$ are indeed $F$-homogeneous of degrees one and two on $ \mathbb{A} $, respectively.  Thus, in particular, $p_1$ is linear on    $Fe + Fa$  if and only if $ p_1$ is additive on $Fe + Fa$. Next,  we show that  $p_1 $ is additive on $Fe + Fa$. 
To see this,  fixing an arbitrary $ a \in \mathbb{A}$, we need to get that $ p_1(a + e) = p_1(a) + p_1(e) = p_1(a) + 1$. To this end, with no loss of generality,  we may assume that $ \{ a, e\}$ is linearly independent, equivalently, $ a \in \mathbb{A} \setminus Fe$. First, note that $ n(a) \not= 0$, for otherwise $ a^2 - 2 p_1(a) a  = 0$, from which we obtain $ a \in    Fe$, which is impossible. Now, we can write 
\begin{eqnarray*}
(a+ e)^2a  - ( a - e)^2a  & = &   2 p_1(a +e) (a^2+a) - n(a+e) a \\
 & & - 2 p_1(a -e) (a^2 - a) + n(a-e) a , \\
(a+ e)^2 a + ( a - e)^2 a & = &  2 p_1(a +e) (a^2+a) - n(a+e) a  \\
& &  + 2 p_1(a -e) (a^2 -a) - n(a-e) a,
\end{eqnarray*}
from which, in view of the flexibility of  $\mathbb{A}$, we see that 
\begin{eqnarray*}
4a^2  & = &   2\big[ p_1(a +e)(2 p_1(a) + 1)   -  p_1(a -e)(2 p_1(a) - 1) \big]a \\
 & &   +\big[ n(a-e) - n(a+e)  \big] a  - 2 n(a)  \big[ p_1(a +e)  -  p_1(a -e)  \big] e , 
 \end{eqnarray*}
 \begin{eqnarray*}
2 a^2 a + 2 a    & = &   2\big[ p_1(a +e)(2 p_1(a) + 1)   +  p_1(a -e)(2 p_1(a) - 1)  \big]a\\
 & &   - \big[ n(a+e) + n(a-e)  \big] a  - 2 n(a)  \big[ p_1(a +e)  + p_1(a -e)  \big] e . 
\end{eqnarray*}
But $ 4a^2 = 8 p_1(a) a - 4 n(a) e$ and 
$$ 2 a^2 a + 2 a  = ( 8 p_1(a)^2 - 2 n(a) + 2 ) a - 4 p_1(a) n(a) e. $$
It thus follows that 
\begin{eqnarray*}
8 p_1(a)  & = &   2\big[ p_1(a +e)(2 p_1(a) + 1)   -  p_1(a -e)(2 p_1(a) - 1)\big] \\
 & &   + \big[ - n(a+e) + n(a-e)  \big], \\
 - 4 n(a) & = & -2 n(a)  \big[ p_1(a +e)  -  p_1(a -e)  \big], 
 \end{eqnarray*}
 \begin{eqnarray*}
 8 p_1(a)^2 - 2 n(a) + 2     & = &  2\big[ p_1(a +e)(2 p_1(a) + 1)   +  p_1(a -e)(2 p_1(a) - 1)  \big]\\
 & &   + \big[- n(a+e) - n(a-e)  \big], \\
 - 4 p_1(a) n(a) & = &  -2 n(a)  \big[ p_1(a +e)  + p_1(a -e)  \big].
\end{eqnarray*}
But  $n(a) \not= 0$. Thus,  
\begin{eqnarray*}
2 & = &  p_1(a +e)  -  p_1(a -e), \\
2 p_1(a)  & = &  p_1(a +e)  + p_1(a -e) , 
\end{eqnarray*}
implying that  $ p_1(a + e) = p_1(a) + p_1(e)$, as desired. This, in turn, yields 
$ n (a) = 2 p_1(a)^2 - p_1( a^2)$, from which we easily obtain $ n ( a + e) + n ( a - e) = 2 (n(a) + n(e)) = 2 ( n(a) + 1)$. Then again, since  $ p_1(a + e) = p_1(a) + 1$, we get that $ n (a + e) - n(a - e) = 4 p_1(a)$. The relations we have shown so far,  together with $( a + e)^2  - ( a - e)^2 = 2 ae + 2a $, reveal that $ ae = a$. That is, $ e$ is also a right identity element of $\mathbb{A}$, for $ a \in \mathbb{A}$ was arbitrary. 

It remains to be shown that $\mathbb{A}$ is locally a field extension of $F$. Since $\mathbb{A}$ is quadratic, it suffices to show that $ ( a - f_1e ) ( a - f_2e ) = 0$ with $ a \in \mathbb{A}$ and $ f_1 , f_2 \in F$ implies $ a= f_1e $ or $ a= f_2e$. We can write 
$$ ( a - f_1 e)^2 - ( f_2 - f_1)  ( a - f_1e ) = ( a - f_1 e) ( a - f_2 e)  = 0, $$
which yields $ a - f_1e = 0$ or  $ a - f_1e  = (f_2 - f_1) e$, which means $ a = f_1e $ or  $ a  = f_2e $, as desired. This completes the proof.  

(ii) 
 First, we prove the existence. To this end, likewise as done in (i), define the functions $ p_1, n: \mathbb{A} \longrightarrow F$ as follows. On $F1$, for $ a = t 1$ with  $t \in F$, set 
$$ p_1 (a) := t, \  n(a) := t^2 ,$$
and for $ a \in \mathbb{A} \setminus F1$, let $p_1(a)$ and $n(a)$ be the unique elements of $F$ satisfying the relation $ a^2 - 2 p_1(a) a + n(a) 1 = 0$. Just as shown in part (i),   $p_1$ and $n$ are indeed $F$-homogeneous of degrees one and two on $ \mathbb{A} $, respectively.  In particular, $p_1 \in  \mathbb{A}'$  if and only if $ p_1$ is additive on $ \mathbb{A} $.
Next,  we show that  $p_1 \in  \mathbb{A}'$. To see this,  fixing $ a, b \in \mathbb{A}$, we need to get that $ p_1(a + b) = p_1(a) + p_1(b)$. To this end, there are two cases to consider. (a) $ \{ a, b , 1\}$ is linearly dependent; and (b) $ \{ a, b , 1\}$ is linearly independent. First, if $ \{ a, b , 1\}$ is linearly dependent, either $ a, b \in F1$, in which case the assertion follows because $p_1$ is $F$-linear on $ F1$, or else, with no loss of generality,  we may assume that $ \{ a, 1\}$ is linearly independent and $ b= r_0 1 + s_0 a$ for some $ r_0 ,  s_0 \in F$. In this subcase, the assertion is reduced to showing that $ p_1(a + 1) = p_1(a) + p_1(1)$ for all $ a \in \mathbb{A} \setminus F1$, equivalently, $ \{ a, 1\}$ is linearly independent. A standard argument similar and in fact simpler than  that presented in part (i) establishes this.  That is,  $ p_1(a + 1) = p_1(a) + p_1(1)$ for all $ a \in \mathbb{A} $, as desired. Next, if $ \{ a, b , 1\}$ is linearly independent, we have 
\begin{eqnarray*}
(a+ b)^2  + ( a - b)^2  & = &   2 p_1(a +b) (a+b) - n(a+b) 1 \\
 & & + 2 p_1(a -b) (a-b) - n(a-b) 1 , 
 \end{eqnarray*}
 and hence 
 \begin{eqnarray*}
4 p_1(a) a   + 4 p_1(b) b-  2\big( n(a) + n(b)  \big) 1  & = &  2 \big( p_1(a +b) + p_1(a -b) \big) a  \\
& & +  2 \big( p_1(a +b)  -  p_1(a -b) \big) b \\
 & &  - \big( n(a+b)  + n(a-b) \big) 1,
\end{eqnarray*}
from which, we get 
 \begin{eqnarray*}
 4 p_1(a) & = & 2 \big( p_1(a +b) + p_1(a -b) \big), \\
4  p_1(b)  & = & 2 \big( p_1(a +b)  -  p_1(a -b) \big).
\end{eqnarray*}
Thus, $ p_1(a + b) = p_1(a) + p_1(b)$, as desired.  Therefore,   $p_1 \in  \mathbb{A}'$.
 This, along with $ a^2 - 2 p_1(a) a + n(a) 1 = 0$, yields  
$$n(a) =2 p_1(a)^2 - p_1(a^2), $$
 for all $ a  \in  \mathbb{A}$, and hence 
$$  n(a+ b) + n( a -b) = 2 n(a) + 2n(b), $$ 
 for all $ a , b \in  \mathbb{A}$. At this point, we are ready to  define 
$$ \langle . , . \rangle: \mathbb{A} \times \mathbb{A} \longrightarrow F,  \  \langle a, b \rangle := \frac{1}{4}\Big( n(a+b) - n(a - b) \Big).$$
With what we have shown so far, it is easily checked that $  \langle a, b \rangle =2 p_1(a) p_1(b) - \frac{1}{2}p_1(ab + ba)$ and that $p_1(a) = \langle a, 1 \rangle$ and 
 $ n(a)  =\langle a, a \rangle $ for  all $ a , b \in  \mathbb{A}$. Then again, as  $p_1 \in  \mathbb{A}'$, it is evident that 
$ \langle . , . \rangle$ is a bilinear and  symmetric form, and that 
$$ a^2 - 2  \langle a , 1 \rangle a + \langle a , a \rangle 1  =  0,$$
\begin{eqnarray*}
  \langle a, b \rangle &  = & 2 \langle a , 1 \rangle  \langle b , 1 \rangle -  \frac{1}{2}  \langle ab + ba , 1 \rangle, \\
   &  = & 2 \langle a , 1 \rangle  \langle b , 1 \rangle -  \frac{1}{4}  \langle (a + b)^2 - ( a  - b)^2, 1 \rangle, 
\end{eqnarray*}
for all $ a, b \in \mathbb{A}$. This proves the existence part of the first assertion.

Next, for the uniqueness, suppose $\langle . , . \rangle_1 :  \mathbb{A} \times \mathbb{A} \longrightarrow F$ is a bilinear form such that 
$$ a^2 - 2  \langle a , 1 \rangle_1 a + \langle a , a \rangle_1 1 = 0, $$
for all $ a \in \mathbb{A}$. Note first that $ \langle a , a \rangle_1 = 2 \langle a , 1 \rangle_1^2 - \langle a^2 , 1 \rangle_1  $  and  hence $ \langle 1 , 1\rangle_1 = 1$ because $ 1 \not = 0$, from which we see that $\langle . , . \rangle_1|_{F \times F} = \langle . , . \rangle|_{F \times F} $, for  $\langle . , . \rangle $ and   $\langle . , . \rangle_1$   are bilinear forms. Once again,  since $\langle . , . \rangle_1$  is bilinear, we have 
\begin{eqnarray*}
 \langle a, b \rangle_1   & = &   \frac{1}{4} \Big(\langle a + b,  a +b \rangle_1 - \langle a - b,  a -b \rangle_1 \Big) , \\
 & = &  2 \langle a , 1 \rangle_1  \langle b , 1 \rangle_1 -  \frac{1}{2}  \langle ab + ba , 1 \rangle_1,
\end{eqnarray*}
for all $ a, b \in \mathbb{A}$. Thus, $\langle . , . \rangle = \langle . , . \rangle_1$  if and only if  $\langle . , 1 \rangle = \langle . , 1 \rangle_1$. To see this, if $ a \in F 1$, the equality trivially follows from  $ \langle 1 , 1\rangle_1 = 1 = \langle 1 , 1\rangle$. If $ a \in   \mathbb{A}\setminus F 1$, we can write 
$$
0 = a^2 - a^2 = 2  \big(  \langle a , 1 \rangle_1  - \langle a , 1 \rangle \big) a + \big( \langle a , a \rangle_1  - \langle a , a \rangle \big) 1 . 
$$
This, in particular, yields $ \langle a , 1 \rangle_1  = \langle a , 1 \rangle$ for  all $ a \in  \mathbb{A}\setminus F 1$. Thus, $ \langle a , 1 \rangle_1  = \langle a , 1 \rangle$ for  all $ a \in  \mathbb{A}$, proving the uniqueness, and hence completing the proof of the first assertion.

Finally,  first, suppose that $0$ is the only element of  $\mathbb{A}$ whose square is $0$ and that $b \in  \mathbb{A}$ is given such that $  \langle b , a\rangle = 0$ for all $ a \in \mathbb{A}$. In particular, we get that $  \langle b , 1\rangle = 0$ and $  \langle b , b\rangle = 0$, from which we obtain $ b^2 = 0$. Thus, $ b =0$, as desired. Next, 
suppose that $\mathbb{A}$ is locally a field extension of $F$. Pick an $ a \in  \mathbb{A}$ with $ \langle a, a \rangle =  n(a) =0$. It follows that $ a ( a - 2 p_1(a) 1) = 0$, from which we get $ a = 0$ or $ a = 2 p_1(a) 1$ because $\mathbb{A}$ is locally a field extension of $F$. If $ a = 2 p_1(a) 1$, then $ 0 = n(a) = 4 p_1(a)^2$, and hence $ p_1(a) = 0$, which, in turn, yields $ a = 2 p_1(a) 1= 0$. In other words, $ \langle a, a \rangle =0$ implies $ a =0$. This completes the proof. 
\hfill \qed

\bigskip

At this point, a few comments are in order. Quadratic algebras that are locally field extensions of their ground fields naturally remind one of inner-product spaces. 
  Let's be more precise and let $\mathbb{A}$  be a  quadratic algebra over a field $F$ of characteristic not $2$ that is locally a  field extension of its ground field $F$ and $\langle . , . \rangle :  \mathbb{A} \times \mathbb{A} \longrightarrow  F$ be its bilinear form, which is also called the {\it inner product} of $\mathbb{A}$,  as  defined in part (ii) of the preceding theorem. For $ a \in \mathbb{A}$,  put
\[
p(a):=\ip{a}{1},\qquad n(a):=\ip{a}{a},\qquad a^{*}:=2p(a)1-a .
\]
An element $ a \in \mathbb A$ is said to be {\it pure} if $ p(a) = 0$. Clearly, every element $ a \in \mathbb A$ can be uniquely written as $ a = \alpha 1 + w $ with $ \alpha \in F$ and $p(w) = 0$; indeed, 
$ \alpha = p(a)$ and $ w= a - p(a) 1$.     Elements $ a , b \in \mathbb{A}$, are said to be {\it orthogonal} if  $\langle a , b \rangle = 0$.  For $S \subseteq \mathbb{A}$, by definition $ S^\perp := \{ a \in \mathbb{A} : \langle a , s \rangle = 0 , \ \forall s \in S\}$. 
  
  First, a standard argument utilizing Zorn's Lemma ensures that  the algebra $\mathbb{A}$ has an orthogonal basis containing  any given orthogonal subset of the algebra $\mathbb{A}$ consisting of nonzero elements. Naturally, by an {\it orthogonal basis} of $\mathbb{A}$, we mean a maximal orthogonal subset of the algebra whose elements are nonzero and pairwise orthogonal. Note that, in general, it does not make sense to talk about orthonormal bases  for quadratic algebras over general fields unless their ground field is the field of real numbers and the algebras are locally complex.  That being noted, locally complex algebras have orthonormal bases containing  any given orthonormal subset of the algebra.  If the algebra $\mathbb{A}$ is finite-dimensional, orthogonal bases of  $\mathbb{A}$  are seen to be spanning. Likewise,  orthonormal bases of  locally complex algebras are seen to be spanning in the sense of Hilbert spaces provided that the algebras equipped with their inner-products are Hilbert spaces.
 
Next, the linear function $ ^*: \mathbb{A} \longrightarrow  \mathbb{A}$, defined by $ a^* = 2 p(a) - a$ ($ a \in \mathbb{A}$),  is called the {\it canonical involution} of $  \mathbb{A}$. A quadratic algebra $\mathbb{A}$ over a field $F$ of characteristic not $2$ is said to be {\it proper} if $    \langle ab ,  b \rangle =    \langle b  ,  b \rangle  p(a)$ for all  $ a, b \in \mathbb{A}$. We will see that if $\mathbb{A}$ is a proper quadratic algebra, then $    \langle ab ,  b \rangle =    \langle b  ,  b \rangle  p(a)= \langle ba ,  b \rangle $  for all  $ a, b \in \mathbb{A}$. Also, it turns out that commutative (resp. left or right alternative) quadratic algebras are proper. Finally, by Theorem \ref{02.8} below, a quadratic algebra is proper if and only if it is flexible.  By the following proposition, whose main part is well-known, the canonical involution of   a proper quadratic algebra $  \mathbb{A}$ is indeed an algebra involution and vice versa. Given a linear transformation $ T \in \mathcal{L}(\mathbb{A})$, the adjoint of $T$, denoted by $T^* \in \mathcal{L}(\mathbb{A})$ is defined by the relations
  $$ \langle Ta , b \rangle = \langle a , T^*b \rangle, \  (a, b \in  \mathbb{A}).$$ 
  One can see that the adjoint of $T$, if exists, is unique, thus justifying the symbol $T^*$. In finite-dimensions, the adjoint always exists.   In terms of adjoints, properness of a quadratic algebra $\mathbb{A}$ translates into the identities $( L_a)^* = L_{a^*}$, or equivalently, $(R_a)^* = R_{a^*}$ for all $ a \in \mathbb{A}$.

 \bigskip

\begin{prop} \label{02.4} 
 Let $F$ be a field whose characteristic is not $2$,  $\mathbb{A}$  a quadratic algebra over $F$,  $\langle . , . \rangle :  \mathbb{A} \times \mathbb{A} \longrightarrow  F$ and $ ^*: \mathbb{A} \longrightarrow  \mathbb{A}$   its bilinear form and  canonical involution, respectively. Then, for all  $ a, b \in \mathbb{A}$ and $ r , s \in F$, the following hold. 
  \begin{enumerate}
\item[(i)] $(ra+sb)^*=ra^*+sb^*$, $(a^*)^*=a$, $p(a^*)=p(a)$, $\ip{a^*}{b^*}=\ip{a}{b}$ and $\ip{a^*}{b}=\ip{a}{b^*}$. 
\item[(ii)] $a+a^*=2p(a)1$ and $aa^*=a^*a=n(a)1$.
\item[(iii)] $ab+ba=2p(b)a+2p(a)b-2\ip{a}{b}1$. Consequently, $\ip{ab+ba}{a}=2n(a)p(b)$, and $ab+ba=-2\ip{a}{b}1$ whenever $p(a)=p(b)=0$.
\item[(iv)] $b^*a^*-(ab)^*=\bigl(p(ba)-p(ab)\bigr)1$. Hence, $^*$ is an anti-automorphism of $\mathbb{A}$ (equivalently, an involution) if and only if $p(ab)=p(ba)$ for all $a,b\in \mathbb A$.
\item[(v)] The unital subalgebra generated by $a$, i.e.,  $\mathbb{B} :=F1+Fa$,  is a commutative associative subalgebra of $\mathbb A$, $*|_{\mathbb{B}}$ is an algebra automorphism of $\mathbb{B}$, and $n(a^k)=n(a)^k$ for all $k \in \mathbb N$.
\item[(vi)] The algebra $\mathbb A$ is proper if and only if $\ip{ab}{b}=\ip{ba}{b}$ for all $a,b\in \mathbb A$; in that case also $\ip{ba}{b}=\ip{b}{b}p(a)$. Therefore, every commutative quadratic algebra is proper.
\end{enumerate}
\end{prop}

\bigskip

\noindent {\bf Proof.}  
Proofs of (i) and (ii) are easy and omitted for brevity. 

(iii) We can write  
 \begin{eqnarray*}
 ab + ba   & =& \frac{1}{2} (a + b)^2 - \frac{1}{2}( a - b)^2, \\
 & = &    \langle a + b , 1 \rangle ( a+ b) - \frac{1}{2} \langle a + b , a + b  \rangle \\
  & & -    \langle a - b , 1 \rangle ( a - b) + \frac{1}{2}\langle a - b , a - b  \rangle, \\
 & = & 2  \langle b , 1 \rangle  a + 2  \langle a , 1 \rangle  b - 2  \langle a , b \rangle 1.   
 \end{eqnarray*}
The rest easily follows from this identity. 

(iv) 
Using (iii) and $(ab)^*=2p(ab)1-ab$,
\begin{align*}
b^*a^*-(ab)^*&=\bigl(4p(a)p(b)-2p(ab)\bigr)1-2p(b)a-2p(a)b+(ab+ba)\\
&=\bigl(4p(a)p(b)-2p(ab)-2\ip ab\bigr)1 .
\end{align*}
By Theorem \ref{02.3}(ii), we have 
$\ip{a}{b}=2p(a)p(b)-\tfrac12\,p(ab+ba)$, 
and hence, 
$4p(a)p(b)-2\ip ab=p(ab+ba)$, from which the assertion follows.

(v) Since $ \mathbb A$ is quadratic, $\mathbb B$ is clearly  a commutative associative subalgebra of $\mathbb A$. By (iv), $*|_{\mathbb{B}}$ is an algebra automorphism of $\mathbb{B}$. Thus, $(a^k)^*=(a^*)^k$ and, by (ii) applied in $\mathbb B$,
\[
n(a^k)1=a^k(a^k)^*=a^k(a^*)^k=(aa^*)^k=n(a)^k1, 
\]
for all $ k \in \mathbb N$.

(vi) By (iii) with the roles of $a,b$ exchanged, we get $\ip{ab}{b}+\ip{ba}{b}=2n(b)p(a)$. Hence, $\ip{ab}{b}=n(b)p(a)$, iff $\ip{ab}{b}=\ip{ba}{b}$, iff $\ip{ba}{b}=n(b)p(a)$. If $\mathbb A$ is commutative, we trivially have $\ip{ab}{b}=\ip{ba}{b}$, completing the proof. 
\hfill \qed

\bigskip 

\noindent {\bf Remark.} In general, canonical involutions of quadratic algebras are not algebra involutions. To see this, let $V$ be a real inner product space with a nonzero skew-symmetric bilinear form $\sigma$ (so $\dim V\ge2$), and let $\mathbb A :=\mathbb R1\oplus V$ with
\[
(\alpha+x)(\beta+y) :=\alpha\beta-\ip xy+\sigma(x,y)+\alpha y+\beta x .
\]
Since $\sigma(x,x)=0$, $(\alpha+x)^2=2\alpha(\alpha+x)-(\alpha^2+\ip xx)1$, the algebra $\mathbb A$ is quadratic (and indeed locally complex) with $p(\alpha+x)=\alpha$ and $n(\alpha+x)=\alpha^2+\ip xx$. For $V :=\mathbb R^2$ and $\sigma(e_1,e_2)=1$, we get $e_1^2=e_2^2=-1$, $e_1e_2=1$ and $e_2e_1=-1$, and hence  $(e_1e_2)^*=1$. But $e_2^*e_1^*=e_2e_1=-1$, in agreement with Proposition~\ref{02.4}(iv), and $p(e_1e_2)=1\ne-1=p(e_2e_1)$. This algebra is not flexible because $(e_1e_2)e_1=e_1$ while $e_1(e_2e_1)=-e_1$ (and by Theorem~\ref{02.8} below it is not proper either).

\bigskip

 Part (i) of the following theorem is the algebraic counterpart of Theorem \ref{1.10}, and is motivated by \cite[Proposition 2.7.33]{GP}, which is due to A. Rodr\'iguez; see \cite[Theorem 2]{Ro}. The first assertion of part (ii) of the theorem for alternative quadratic algebras  is well-known, see \cite[page 58]{Sch}. Indeed, parts (i) and (ii) of the theorem turn out to be  key ingredients in proving  a slight extension of  a well-known theorem of Hurwitz, see \cite[Theorem 3.25]{Sch} and Theorem \ref{3.9} below. 

\bigskip

\begin{thm} \label{02.5} 
 {\rm (i)} 
Let $F$ be a field whose characteristic is not $2$ and   $\mathbb{A}$  a flexible algebra over $F$ with a left (resp. right) identity element $e$ and with no nontrivial zero-divisors. Let  $\langle . , . \rangle $ be a nondegenerate bilinear form on $\mathbb{A}$ that permits composition, i.e.,  $ \langle ab , ab \rangle = \langle a , a \rangle \langle b , b \rangle$ for all  $ a, b \in \mathbb{A}$. Then, $e$ is an identity element of  the algebra $\mathbb{A}$, the algebra $\mathbb{A}$  is quadratic and locally a field extension of $F$ and alternative, and  $\langle . , . \rangle $ is indeed the  bilinear form of  $\mathbb{A}$. Moreover, $ \langle a  b, c \rangle  =\langle  b, a^* c \rangle  $, $ \langle a  b, c \rangle  =\langle  a, c b^* \rangle  $, $ a^*(ab) = \langle a , a \rangle  b$, and $ (ab)b^* = \langle b , b \rangle a$  for all $ a, b, c \in \mathbb{A}$.

 {\rm (ii)}  
 Let $F$ be a field whose characteristic is not $2$ and $\mathbb{A}$ a  left (resp. right) alternative quadratic algebra  over $F$.  Then,  $\mathbb{A}$ is alternative and its bilinear form permits composition.   Moreover,  if a quadratic algebra $\mathbb{A}$  over  a field $F$ of characteristic not $2$  has no nontrivial joint zero-divisors, then it is locally a field extension of $F$. Finally, if a quadratic algebra $\mathbb{A}$ over a field $F$ of characteristic not $2$ is left (resp. right) alternative and locally a field extension of $F$, then the algebra  $\mathbb{A}$ has no nontrivial  zero-divisors.

\end{thm}

\bigskip

\noindent {\bf Proof.}  
 (i) 
  Let $e$ be a left  identity element of  $\mathbb{A}$. If $e$ is a right  identity element of  $\mathbb{A}$, the proof is carried out in a similar fashion.
The proof of the assertion is an imitation of that of its topological counterpart, which is due to A. Rodr\'iguez; see \cite[Lemma 2.7.32]{GP} and   \cite[Proposition 2.7.33]{GP}. Here, we produce an adjusted proof for the sake of completeness. Let   $ a, b , c   \in \mathbb{A}$ be arbitrary. Set $ u :=  \langle a , e \rangle e - a $ and note that $   \langle u , e \rangle = 0$. For any $ y \in \mathbb{A}$, we can write 
\begin{eqnarray*}
\Big( 1 +  \langle u , u\rangle \Big)  \langle y , y \rangle & =&  \langle  e + u , e + u\rangle \langle y , y \rangle , \\
& =& \langle  y + uy , y + uy\rangle, \\
& =& \langle y , y \rangle + \langle uy , uy \rangle + 2 \langle uy , y \rangle\\
& =& \Big( 1 +  \langle u , u\rangle \Big)  \langle y , y \rangle + 2 \langle uy , y \rangle,
\end{eqnarray*}
which yields $   \langle uy , y \rangle = 0$. Letting $ y = b + c$, we obtain 
$$  \big\langle  \big(\langle a , e \rangle e - a \big) b , c \big\rangle = - \big\langle \big(\langle a , e \rangle e - a \big) c , b \rangle.$$
Simplifying this, we get  $  \langle ab , c \rangle =  \langle b , a^*c \rangle $, where $ a^* = 2 \langle a , e \rangle e - a$.  
Next, since $ \langle ax , ax \rangle = \langle a , a \rangle \langle x , x \rangle$ for all  $ a, x \in \mathbb{A}$, letting $ x = b+c$, where $  b , c   \in \mathbb{A}$ are arbitrary,  we get 
 $ \langle ab , ac \rangle = \langle a , a \rangle \langle b , c \rangle$. By what we just saw, $ \langle b , a^*(ac) \rangle = \langle a , a \rangle \langle b , c \rangle$, which, in turn, implies $ \big\langle b , a^*(ac)  - \langle a , a \rangle  c \big\rangle =0$. Then again, since the bilinear form $ \langle . , . \rangle$ is nondegenerate, we conclude that  $ a^*(ac)  = \langle a , a \rangle  c $, from which, letting $ c = e$, we obtain $ 2  \langle a , e \rangle ae  - a (ae) =  \langle a , a \rangle e  $. Multiplying both sides by $a$ from the right, we get $  2  \langle a , e \rangle (ae)a  - \big( a (ae) \big)a  =  \langle a , a \rangle ea   $. Then, the flexibility of $\mathbb{A}$ yields 
 $  2  \langle a , e \rangle a^2  - a a^2 =  \langle a , a \rangle ea   $ and  $ a a^2 = a^2 a$. It thus follows that 
 $  \big( 2  \langle a , e \rangle a  -  a^2 -  \langle a , a \rangle e \big) a  = 0  $, and hence $  a^2 -  2  \langle a , e \rangle a  + \langle a , a \rangle e = 0$  for all $a \in \mathbb A$, for $\mathbb A$ has no nontrivial zero-divisors. That is,  $\mathbb{A}$ is  a flexible quadratic algebra over $F$ relative to the left identity element $e$. 
 Now, since $\mathbb{A}$ has no nontrivial zero-divisors, we see from Theorem \ref{02.3}(i)  that $ e$ is an identity element, and hence the algebra $\mathbb A$ is quadratic. Next, substituting   for $a^*$  in $ a^*(ac)  = \langle a , a \rangle  c $  and making use of the identity $ x^2 = 2  \langle x , e \rangle x - \langle x , x \rangle e  $ on $\mathbb{A}$, we conclude that $\mathbb A$ is left alternative.  The algebra  $\mathbb A$ being left alternative  and quadratic,  together with  \cite[Theorem 1]{A2},  entails that it is alternative. 
 Now that $\mathbb{A}$ is unital, one can see in a similar  fashion that  $ \langle a  b, c \rangle  =\langle  a, c b^* \rangle  $ and  $ (ab)b^* = \langle b , b \rangle a$
for all $a ,b , c \in \mathbb A$. This completes the proof.

(ii)  For the  first assertion, note that the algebra $  \mathbb{A}$ is alternative by \cite[Theorem 1]{A2}. Thus, we can write 
\begin{eqnarray*}
\langle ab , ab \rangle  & = & \langle (ab) (ab)^* , 1 \rangle, \\
 & = & \langle (ab) (b^* a^*)  , 1 \rangle, \\
 & = & \langle  \big(a (bb^*)  \big) a^*  , 1 \rangle, \\
  & = &  \langle  (aa^*) (bb^*)    , 1 \rangle, \\
 & = & \langle a , a \rangle \langle b , b \rangle,
\end{eqnarray*}
which is what we want.

Next, we  prove the second assertion. To this end, first,  suppose that the quadratic algebra $\mathbb{A}$  has no nontrivial joint zero-divisors and that $ a \in  \mathbb{A}$ is a given nonzero element. We have $ a \big(  2 \langle a , 1 \rangle1  - a \big) = \langle a , a \rangle 1$. From this, we get that $  \langle a , a \rangle \not= 0$, for otherwise 
$$ a \big(  2 \langle a , 1 \rangle1  - a \big) =  \big(  2 \langle a , 1 \rangle1  - a \big)  a=  0,$$ 
 and hence $ a =  2 \langle a , 1 \rangle1 $. This yields  $ (\langle a , 1 \rangle)^2 = 0$, from which we obtain $ a^2 = 0$, implying $ a = 0$, a contradiction. Thus,   $  \langle a , a \rangle \not= 0$, which entails $ a^{-1} = \frac{1}{ \langle a , a \rangle} \big(  2 \langle a , 1 \rangle1  - a \big)$. In other words, the unital algebra generated by $a$ is a field extension of $F$, as desired. 
 
Finally, suppose that the left (resp. right) alternative quadratic algebra $\mathbb{A}$ is locally a field extension of $F$ and that $ ab = 0$ for some $ a, b \in \mathbb{A}$. We need to show that $ a= 0$ or $b = 0$. To this end, we have
$$0 = \langle ab , ab \rangle  = \langle a , a \rangle \langle b , b \rangle,$$
because  $\mathbb{A}$ is alternative by \cite[Theorem 1]{A2}, and hence its bilinear form permits composition. 
It thus follows that $ \langle a , a \rangle = 0  $  or  $  \langle b , b \rangle=0$, which, in turn,  implies $ a \big(  2 \langle a , 1 \rangle1  - a \big) =0$  or $ b \big(  2 \langle b , 1 \rangle1  - b \big) =0$, and hence $ a = 0$  or $2 \langle a , 1 \rangle1  - a  =0$  or $ b = 0$ or $  2 \langle b , 1 \rangle1  - b  =0$, for the algebra generated by any nonzero element of $\mathbb{A}$ has no nontrivial zero-divisors. If $ a =  2 \langle a , 1 \rangle1 $ or $ b =  2 \langle b , 1 \rangle1 $, just as we saw in the above, we get that 
$ a^2= 0$ or $b^2 = 0$. And from this, we conclude that $ a= 0$ or $b = 0$ again because the algebra generated by any nonzero element of $\mathbb{A}$ has no nontrivial zero-divisors, as desired. 
\hfill \qed 

\bigskip

\noindent {\bf Remark.} 
 In part (i), if the algebra $\mathbb{A}$ is assumed to have an identity element, then the hypothesis that $\mathbb{A}$ is flexible is redundant and  the hypothesis that $\mathbb{A}$ has no nontrivial zero-divisors can be replaced by the the weaker hypothesis that $\mathbb{A}$  has no nontrivial joint zero-divisors.

\bigskip

The counterparts of parts (i)-(iv) of the following theorem in finite-dimensional inner-product spaces are standard and well-known. In the spirit of the comments we made preceding Proposition \ref{02.4}, the theorem is quoted here.

\bigskip

\begin{thm} \label{02.6} 
{\rm (i)}  {\bf (The Gram-Schmidt Orthogonalization Process)} Let $F$ be a field whose characteristic is not $2$ and   $\mathbb{A}$  a quadratic algebra over $F$ which is locally a field extension of $F$, $ n \in \mathbb{N}$ and $ \{ a_i\}_{i=1}^n$ be a linearly independent subset of $\mathbb{A}$. Then, $\{ b_i\}_{i=1}^n$ defined inductively by 
$$ b_1 := a_1 , \ b_{i+1} = a_{i+1} - \sum_{j=1}^i  \frac{\langle a_{i+1} , b_{j} \rangle}{\langle b_j , b_j\rangle} b_j , \ ( 1 \leq i < n) $$
is an orthogonal set consisting of nonzero elements with the property that 
$$ \langle \{ a_i\}_{i=1}^j \rangle_F = \langle \{ b_i\}_{i=1}^j \rangle_F, $$
for all $ 1 \leq j \leq n$. Moreover, if  $\{ c_i\}_{i=1}^n$ is an orthogonal set consisting of nonzero elements with the property that 
$$ \langle \{ a_i\}_{i=1}^j \rangle_F = \langle \{ c_i\}_{i=1}^j \rangle_F, $$
for all $ 1 \leq j \leq n$, then there exists nonzero elements $ r_i \in F$ ($ 1 \leq i \leq n$) such that $ c_i = r_i b_i $ for all $ 1 \leq i \leq n$. 

\bigskip

{\rm (ii)} Let $F$ be a field whose characteristic is not $2$ and   $\mathbb{A}$  a quadratic algebra over $F$ which is locally a field extension of $F$, $\langle . , . \rangle :  \mathbb{A} \times \mathbb{A} \longrightarrow  F$ be its bilinear form, and $\mathcal{M}$ be a finite-dimensional subspace of $\mathbb{A}$. Then,  $\mathcal{M} + \mathcal{M}^\perp =\mathbb{A}  $, $\mathcal{M} \cap \mathcal{M}^\perp = \{0\}$, and $ (\mathcal{M}^\perp)^\perp = \mathcal{M}$. Moreover,   any orthogonal subset of $\mathbb{A}$ consisting of nonzero elements is linearly independent and  can be extended to an orthogonal basis of  $\mathbb{A}$. 
 
 \bigskip

{\rm (iii)} {\bf (The Riesz Representation Theorem)}  Let $F$ be a field whose characteristic is not $2$ and   $\mathbb{A}$  a finite-dimensional quadratic algebra over $F$ which is locally a field extension of $F$ and $\langle . , . \rangle :  \mathbb{A} \times \mathbb{A} \longrightarrow  F$ be its bilinear form. Then, corresponding to any $ \varphi \in \mathbb{A}'$, there exists a unique $ a_\varphi \in \mathbb{A}$ such that $ \varphi(a) = \langle a , a_\varphi \rangle$ for all $ a \in \mathbb{A}$. 
 
 \bigskip
 
{\rm (iv)} Let $F$ be a field whose characteristic is not $2$ and   $\mathbb{A}$  a finite-dimensional quadratic algebra over $F$ which is locally a field extension of $F$ and $\langle . , . \rangle :  \mathbb{A} \times \mathbb{A} \longrightarrow  F$ be its bilinear form. Then,  any $ T \in \mathcal{L}( \mathbb{A})$ has a unique adjoint, denoted by $ T^* \in \mathcal{L}( \mathbb{A})$, defined by the relations 
  $$ \langle Ta , b \rangle = \langle a , T^*b \rangle, \    (a, b \in  \mathbb{A}).$$
 
\end{thm}

\bigskip

\noindent {\bf Proof.} The proofs, which are omitted for brevity, are straightforward imitations of those of their counterparts in the context of finite-dimensional inner-product spaces. 
\hfill \qed 

\bigskip

\noindent {\bf Remarks.} 1. The counterpart of part (i) holds for the case when $n = \infty$, namely,  when  $ \{ a_i\}_{i=1}^{\infty}$ is a given countable linearly independent subset of $\mathbb{A}$. 

2. Also, the counterpart of part (i) holds for locally complex algebras. Indeed,  for locally complex algebras, given a linearly independent subset $ \{ a_i\}_{i=1}^{n}$,  the orthogonal subset $\{ b_i\}_{i=1}^n$ can be made to be an orthonormal subset with the property 
 $$ \langle \{ a_i\}_{i=1}^j \rangle_\mathbb{R} = \langle \{ b_i\}_{i=1}^j \rangle_\mathbb{R}, $$
for all $ 1 \leq j \leq n$, in which case we will have $ r_i = \pm 1 $ for all $ 1 \leq i \leq n$.

\bigskip

The following theorem includes  the counterparts of Lemmas 1.1 and 1.3-1.4 and Corollary 1.5 of \cite{M1} for quadratic algebras; also see   \cite[Lemmas 2 and 6]{Sch1}.  

\bigskip

\begin{thm} \label{02.7} 
   Let $\mathbb A$ be a proper quadratic algebra over a field $F$ of characteristic not $2$. Then, for all $a,b,c\in \mathbb A$, the following hold. 
\begin{enumerate}
\item[(i)] $\ip{ab}{1}=\ip{ba}{1}=2p(a)p(b)-\ip ab=\ip{a}{b^*}=\ip{a^*}{b}$. In particular, $ab+ba=-2\ip ab1=2p(ab)1$ whenever $p(a)=p(b)=0$.
\item[(ii)] $\ip{ab}{c}=\ip{b}{a^*c}$ and $\ip{ab}{c}=\ip{a}{cb^*}$; that is, $(L_a)^*=L_{a^*}$ and $(R_b)^*=R_{b^*}$, where adjoints are taken with respect to $\ip\cdot\cdot$.
\item[(iii)] $(ab)^*=b^*a^*$; thus $*$ is an involution of $A$.
\item[(iv)] $\ip{(ab)c}{1}=\ip{a(bc)}{1}$.
\item[(v)] $\ip{ab}{ab}=\ip{a^*b}{a^*b}=\ip{ab^*}{ab^*}=\ip{ba}{ba}$.
\end{enumerate}
  
\end{thm}

\bigskip

\noindent {\bf Proof.} 
Note first that properness gives $\ip{ab}{b}=n(b)p(a)$ and, by Proposition~\ref{02.4}(vi), $\ip{ba}{b}=n(b)p(a)$. Since $\langle . , . \rangle$ is bilinear and symmetric, replacing $b$ by $b+c$ and simplifying,  gives
\begin{equation}\tag{P}
\ip{ab}{c}+\ip{ac}{b}=2p(a)\ip bc,
\end{equation}
\begin{equation}\tag{P$'$}
\ip{ba}{c}+\ip{ca}{b}=2p(a)\ip bc .
\end{equation}

(i) Put $c=1$ in (P) to get $\ip{ab}{1}=2p(a)p(b) - \ip ab$. But the right-hand side is symmetric in $a $ and $b$. Thus, $p(ab)=p(ba)=2p(a)p(b)-\ip ab=\ip a{b^*}=\ip{a^*}{b}$ by Proposition~\ref{02.4}(i). If $p(a)=p(b)=0$, Proposition~\ref{02.4}(iii) gives $ab+ba=-2\ip ab1=2p(ab)1$.

(ii) By (P), we have $\ip{ab}{c}=2p(a)\ip bc-\ip{ac}{b}=\ip b{2p(a)c-ac}=\ip b{a^*c}$. Changing the role of $a$ and $b$ in (P$'$), we similarly get $\ip{ab}{c}=\ip a{2p(b)c-cb}=\ip a{cb^*}$.

(iii) By (i), $p(ab)=p(ba)$. So Proposition~\ref{02.4}(iv) gives $(ab)^*=b^*a^*$. 

(iv) Using (i), (ii) and (iii) in turn, we can write
\[
\ip{(ab)c}{1}=\ip{ab}{c^*}=\ip{a}{c^*b^*}=\ip a{(bc)^*}=\ip{a(bc)}{1}.
\]

(v) Since $p(a-a^*)=0$ and $a+a^*=2p(a)1$, properness gives
\[
\ip{(a+a^*)b}{(a-a^*)b}=2p(a)\ip{b}{(a-a^*)b}=2p(a)\,n(b)\,p(a-a^*)=0 .
\]
Expanding the left-hand side, the two cross terms cancel by symmetry, so $\ip{ab}{ab}=\ip{a^*b}{a^*b}$. Likewise, using the identity $\ip{xy}{x}=n(x)p(y)$, we have
\[
\ip{a(b+b^*)}{a(b-b^*)}=2p(b)\ip{a}{a(b-b^*)}=2p(b)\,n(a)\,p(b-b^*)=0,
\]
whence $\ip{ab}{ab}=\ip{ab^*}{ab^*}$. Applying both equalities, we get $\ip{ab}{ab}=\ip{a^*b^*}{a^*b^*}=n\bigl((ba)^*\bigr)=n(ba)$, where we used (iii) and Proposition~\ref{02.4}(i).
\hfill \qed 

\bigskip

The following theorem asserts that the notions of being proper and flexible in quadratic algebras are equivalent. 

\bigskip

\begin{thm} \label{02.8} 
Let $\mathbb A$ be a quadratic algebra over a field $F$ of characteristic not $2$. Then, 
 the algebra $A$ is proper if and only if it is flexible.
\end{thm}

\bigskip

\noindent {\bf Proof.}  \emph{Proper $\Rightarrow$ flexible.} In the flexible  identity $(ab)a=a(ba)$, writing $a=\alpha1+a_0$ and $b=\beta1+b_0$ with $a_0,b_0 \in \{ 1\}^\perp$, all terms involving $\alpha$ or $\beta$ cancel, and we obtain
\[
(ab)a-a(ba)=(a_0b_0)a_0-a_0(b_0a_0).
\]
So to prove the flexible identity, we may assume that $p(a)=p(b)=0$, i.e., \ $a^*=-a$, $b^*=-b$. Set $t:=a(ba)-(ab)a$. By Theorem~\ref{02.7}(iii), we have
\[
\bigl(a(ba)\bigr)^*=(ba)^*a^*=(a^*b^*)a^*=-(ab)a,
\]
 \[\bigl((ab)a\bigr)^*=a^*(ab)^*=a^*(b^*a^*)=-a(ba).
\]
Thus, $t^*=t$, which means $2p(t)1-t=t$, and hence $t=p(t)1$. By Theorem~\ref{02.7}(iv), we get $p(t)=\ip{a(ba)}1-\ip{(ab)a}1=0$, which in turn, yields $t=0$, as desired.

\emph{Flexible $\Rightarrow$ proper.}

In light of Proposition~\ref{02.4}(vi), it suffices to prove $\ip{ab}{b}=\ip{ba}{b}$ for all $a,b \in \mathbb A$. Write $a=\alpha1+a_0$ and $b=\beta1+b_0$ with $a_0,b_0 \in \{ 1\}^\perp$. A straightforward calculation reveals that 
  \begin{eqnarray*}
\ip{ab}{b}-\ip{ba}{b} & = & \ip{ab -ba}{b},\\
& = &  \ip{a_0b_0 -b_0a_0}{ \beta 1 + b_0}, \\
& = &  \beta p(a_0b_0 -b_0a_0) +   \ip{a_0b_0}{b_0}-\ip{b_0a_0}{b_0}.
  \end{eqnarray*}
By showing that $ p(a_0b_0 ) = p (b_0a_0)$, we may assume that the elements $a , b \in \mathbb A$ are pure.  To this end, suppose that $ a , b  \in \mathbb A$ are pure. Write 
$ ab = p(ab) 1 + w$ with $w$ being pure. By Proposition~\ref{02.4}(iii), we have $ab+ba=-2\ip ab1$. So $ba= \big(-2\ip ab - p(ab)\Big) 1- w$, which yields $p(ba)=-2\ip ab-p(ab)$. Then, flexibility gives
 \begin{eqnarray*}
0 & =& (ab)a-a(ba), \\
& = & \big(p(ab) - p(ba)\big) a+(wa+aw), \\
& =& \big(p(ab) - p(ba)\big)a-2\ip wa1,
  \end{eqnarray*}
by Proposition~\ref{02.4}(iii), since $w,a$ are pure. If $a\neq0$, then $a$ and $1$ are linearly independent, so $p(ab) = p(ba)$; if $a=0$, then obviously $p(ab) =p(ba)=0$.

Finally, assuming that $ a, b$ are pure, by Theorem~\ref{02.3}(ii), we have $\ip xb=-\tfrac12p(xb+bx)$ for every $x$.  Hence,
\[
\ip{ab}{b}-\ip{ba}{b}=-\tfrac12\,p\Bigl(\bigl[(ab)b-b(ba)\bigr]+\bigl[b(ab)-(ba)b\bigr]\Bigr).
\]
The second bracket vanishes by flexibility. For the first braket, $ab=-ba-2\ip ab1$ and $ba=-ab-2\ip ab1$ give $(ab)b-b(ba)=-(ba)b+b(ab)$, which is again $0$ by flexibility. Thus, $\ip{ab}{b}=\ip{ba}{b}$, completing the proof. 
\hfill \qed 

\bigskip

The following is the counterpart of \cite[Corollary 1.6]{M1} for quadratic algebras. For descriptions of zero-divisors of the real Cayley-Dickson algebras, we refer the reader to  \cite{M1} and \cite{BDI}.

\bigskip

\begin{prop} \label{02.9} 
 Let $F$ be a field whose characteristic is not $2$ and   $\mathbb{A}$  a flexible, equivalently, proper, quadratic algebra over $F$ which is locally a field extension of $F$ and $ a, b \in \mathbb{A}$. Then, $ ab = 0$ if and only if $ ba = 0$. Moreover, if  $ a, b \in \mathbb{A} \setminus F 1$ and $ ab = 0$, then  $ a, b \in \{1\}^\perp$.  
 
\end{prop}

\bigskip

\noindent {\bf Proof.}   If   $ a \in F1$ or $ b \in F1$, the assertion trivially holds.  Suppose $ a, b \in \mathbb{A} \setminus F 1$. We have $0 =   \langle ab , a \rangle =  \langle a , a \rangle  \langle b , 1 \rangle $, from which we obtain  $  \langle b , 1 \rangle  =0$, for $  \langle a , a \rangle \not= 0$. Likewise,  $  \langle a , 1 \rangle  =0$. That is, $ a^* = -a$ and $  b^* = -b$. But $ ab = 0$. Taking $*$ of both sides,  or making use of Theorem \ref{02.7}(i), gives $ ba = 0$, as desired. 
\hfill \qed

 \bigskip

   An important class of finite-dimensional quadratic algebras is the Cayley-Dickson algebras; see  \cite{Sch1} and  \cite[pages 5 and 55-61]{Sch}. Consequently, a natural question arises: {\it under what conditions is a finite-dimensional quadratic algebra over a field $F$ isomorphic to a Cayley-Dickson algebra defined over $F$?}   An answer to this question would provide a unified extension of the classical theorems of Frobenius, Zorn, and Hurwitz in finite dimensions.

\bigskip

The following theorem gives necessary and sufficient conditions for a function to preserve the bilinear structure of quadratic algebras that are locally field extensions of their ground fields.

\bigskip

\begin{thm} \label{02.10} 
 {\rm (i)}   Let $F$ be a field whose characteristic is not $2$,  $\mathbb{A}$  and $\mathbb{B}$ quadratic algebras over $F$ that are locally  field extensions of $F$, and $\langle . , . \rangle_\mathbb{A} $ and $\langle . , . \rangle_\mathbb{B} $ their bilinear forms.  If   $ f : \mathbb{A} \longrightarrow  \mathbb{B}$ is a surjective $F$-linear function with $ f(a^2) = f(a)^2$ for all  $ a \in \mathbb{A}$, then $f$ is a linear $*$-isomorphism preserving the bilinear form, i.e.,  $ \langle f(a_1) , f(a_2) \rangle_ \mathbb{B} =  \langle a_1 , a_2 \rangle_\mathbb{A}$ for all $ a_1 , a_2 \in \mathbb{A}$; in particular, $ f : \mathbb{A}^{\rm sym} \longrightarrow  \mathbb{B}^{\rm sym}$ is an $F$-algebra  $*$-isomorphism preserving the bilinear form. Moreover, if  $ f : \mathbb{A} \longrightarrow  \mathbb{B}$ is a surjective  function preserving the bilinear form with $ f(1_\mathbb{A}) = 1_\mathbb{B}$, then $ f $  is an $F$-linear $*$-isomorphism with $ f(a^2) = f(a)^2$ 
   for all $ a \in \mathbb{A}$. 
 
 {\rm (ii)}     Let $F$ be a field whose characteristic is not $2$,  $\mathbb{A}$  and $\mathbb{B}$  quadratic algebras over $F$ that are locally  field extensions of $F$, $\langle . , . \rangle_\mathbb{A} $ and $\langle . , . \rangle_\mathbb{B} $ their bilinear forms, and  $ f : \mathbb{A} \longrightarrow  \mathbb{B}$  a surjective function. Then, the following statements are equivalent.
 
 {\rm (a)}  The function $f$ is $F$-linear and  $ f(a^2) = f(a)^2$ for all  $ a \in \mathbb{A}$;

 {\rm (b)}    The function $f$ preserves  the bilinear form and $ f(1_\mathbb{A}) = 1_\mathbb{B}$; 
 
 {\rm (c)} The function $f$ is additive, $ f(1_\mathbb{A}) = 1_\mathbb{B}$, and $ \langle f(a) , f(a) \rangle_ \mathbb{B} =  \langle a , a \rangle_\mathbb{A}$ for all $ a \in \mathbb{A}$.

 {\rm (iii)}  Let $F$ be a field whose characteristic is not $2$,  $\mathbb{A}$  and $\mathbb{B}$  quadratic algebras over $F$ that are locally  field extensions of $F$, $\langle . , . \rangle_\mathbb{A} $ and $\langle . , . \rangle_\mathbb{B} $ their bilinear forms, and  $ f : \mathbb{A} \longrightarrow  \mathbb{B}$  a surjective function. Then, the function $f$ is an $F$-algebra isomorphism if and only if  $f$ is multiplicative,  $ f(1_\mathbb{A}) = 1_\mathbb{B}$, and $ \langle f(a) , 1_\mathbb{B} \rangle_ \mathbb{B} =  \langle a , 1_\mathbb{A} \rangle_\mathbb{A}$ for all $ a \in \mathbb{A}$, in which case  $f$ is an $F$-algebra $*$-isomorphism. 
 
\end{thm}

\bigskip

\noindent {\bf Proof.} (i) 
  First, suppose that $ f : \mathbb{A} \longrightarrow  \mathbb{B}$ is a surjective $F$-linear function with $ f(a^2) = f(a)^2$. It follows that $ f(1_\mathbb{A}) = f(1_\mathbb{A})^2$, from which we get that $ f(1_\mathbb{A}) = 0$ or $f(1_\mathbb{A}) = 1_\mathbb{B}$, because $\mathbb{B}$ is locally a field extension of $F$. But $ f(1_\mathbb{A}) \not= 0$, for otherwise $ f(1_\mathbb{A}) = 0$ along with $ f\big((a + 1_\mathbb{A})^2\big) = \big(f(a+1_\mathbb{A})\big)^2 $ implies that $f(a) = 0$ for all $ a \in   \mathbb{A}$, which is a contradiction because  $ f $ is surjective. Thus, $f(1_\mathbb{A}) = 1_\mathbb{B}$. We have $ a^2 - 2  \langle a  ,  1_\mathbb{A}  \rangle_ \mathbb{A}  a +   \langle a  ,  a  \rangle_\mathbb{A}  1_\mathbb{A} = 0$ for all $ a \in   \mathbb{A}$. Taking $f$ of both sides and using the identity $ f(a^2) = f(a)^2$ on $ \mathbb{A}$ yield 
  $ f(a)^2 - 2  \langle a  ,  1_\mathbb{A}  \rangle_\mathbb{A}  f(a)  +   \langle a  ,  a  \rangle_\mathbb{A}  1_\mathbb{B} = 0$ for all $ a \in   \mathbb{A}$. This implies $  \langle f(a)  ,  1_\mathbb{B}  \rangle_\mathbb{B} =  \langle a  ,  1_\mathbb{A}  \rangle_\mathbb{A}  $ and $  \langle f(a)  ,  f(a)  \rangle_\mathbb{B} =  \langle a  ,  a  \rangle_ \mathbb{A} $  for all $ a \in   \mathbb{A}$. Now, we can write
  \begin{eqnarray*}
  \langle f(a)  ,  f(a')  \rangle_\mathbb{B}  & = &  2 \langle f(a) , 1_\mathbb{B} \rangle_\mathbb{B}  \langle f(a') , 1_\mathbb{B} \rangle_\mathbb{B}\\
    & &  -  \frac{1}{4}  \langle ( f(a) + f(a') )^2 - ( f(a) - f(a') )^2, 1_\mathbb{B} \rangle_\mathbb{B},\\
   & = & 2 \langle a , 1_\mathbb{A} \rangle_\mathbb{A}   \langle a' , 1_\mathbb{A} \rangle_\mathbb{A}  -  \frac{1}{4}  \langle (a + a')^2 - ( a  - a')^2, 1_\mathbb{A} \rangle_\mathbb{A} ,\\
     & = &\langle a  ,  a'  \rangle_\mathbb{A} , 
  \end{eqnarray*}
  for all $ a  ,  a' \in \mathbb{A} $. That is, $f$ preserves the bilinear form. To see that $f$ is injective, suppose that $ f(a) = 0$. We have 
  $0  =  \langle f(a)  ,  f(a)  \rangle_\mathbb{B} = \langle a  ,  a  \rangle_\mathbb{A} $, from which, in view of Theorem \ref{02.3}(ii), we obtain $ a = 0$ since  $\mathbb{A}$ is locally a field extension of $F$.  Next, we can write 
    \begin{eqnarray*}
  f(a^*)  & = &  f \big( 2 \langle a  ,  1_\mathbb{A}  \rangle_\mathbb{A} 1_\mathbb{A} - a \big) , \\
  & = & 2 \langle a  ,  1_\mathbb{A}  \rangle_\mathbb{A} 1_\mathbb{B} - f(a) , \\
  & = & 2 \langle f(a)  ,  1_\mathbb{B}  \rangle_\mathbb{B} 1_\mathbb{B} - f(a) ,\\
  & = &   f(a)^* , 
    \end{eqnarray*}
   for all $ a   \in \mathbb{A} $. Now that $f$  preserves  the bilinear form, the assertion that $ f : \mathbb{A}^{\rm sym} \longrightarrow  \mathbb{B}^{\rm sym}$ is an $F$-algebra  $*$-isomorphism  follows from the identity 
 $$ab + ba = 2 \langle b , 1 \rangle  a + 2  \langle a , 1 \rangle  b - 2  \langle a , b \rangle 1,$$ 
 on   quadratic algebras, which we proved in Proposition \ref{02.4}.

   Next, suppose that $ f : \mathbb{A} \longrightarrow  \mathbb{B}$ is a surjective function  preserving the bilinear form with $ f(1_\mathbb{A} ) = 1_\mathbb{B}$. 
  First, for each $ r \in F$  and $ a , b , c \in \mathbb{A}$,  we can write 
   \begin{eqnarray*}
   \langle f ( ra + b) - r f(a) -  f(b)    ,  f(c)   \rangle_\mathbb{B} & = & \langle f ( r a + b)     ,  f(c)   \rangle_\mathbb{B}  - r\langle    f(a)     ,  f(c)   \rangle_\mathbb{B}\\ 
   & &  -   \langle  f(b)     ,  f(c)   \rangle_\mathbb{B}, \\
  & = & \langle r a + b     ,  c  \rangle_\mathbb{A}  - r \langle  a     ,  c   \rangle_\mathbb{A} -   \langle b     ,  c   \rangle_\mathbb{A},\\
  & = & 0, 
  \end{eqnarray*}
   implying that $f$ is  $F$-linear, for $f$ is surjective and $ \mathbb{B}$  is locally a field extension of $F$.  Second,   we have 
     \begin{eqnarray*}
   a^2 - 2  \langle a  ,  1_\mathbb{A}  \rangle_ \mathbb{A}  a +   \langle a  ,  a  \rangle_\mathbb{A}  1_\mathbb{A}  & = & 0, \\
    f(a)^2 - 2  \langle f(a)  ,  1_\mathbb{B}  \rangle_\mathbb{B}  f(a)  +   \langle f(a)  , f(a)   \rangle_\mathbb{B}  1_\mathbb{B}   & = & 0,
       \end{eqnarray*}
    for all $ a   \in \mathbb{A} $. Taking $f$ of both sides of the first equality above and using the hypotheses that $ f(1_\mathbb{A} ) = 1_\mathbb{B}$ and that $f$  preserves  the bilinear form, we obtain 
    \begin{eqnarray*}
   f(a^2)  - 2  \langle f(a)  ,  1_\mathbb{B}  \rangle_ \mathbb{B} f(a)   +   \langle f(a)  ,  f(a)  \rangle_\mathbb{B}  1_\mathbb{B}  & = & 0, \\
    f(a)^2 - 2  \langle f(a)  ,  1_\mathbb{B}  \rangle_\mathbb{B}  f(a)  +   \langle f(a)  , f(a)   \rangle_\mathbb{B}  1_\mathbb{B}   & = & 0,
       \end{eqnarray*}
   from which we get  
   $$f(a^2) =  f(a)^2= 2  \langle f(a)  ,  1_\mathbb{B}   \rangle_\mathbb{B}  f(a)  -   \langle f(a)  , f(a)   \rangle_\mathbb{B}  1_\mathbb{B}  $$
    for all $ a   \in \mathbb{A} $, as desired. Now, that $f$ is injective and is a $*$-homomorphism  follows from the first assertion.

 (ii) By (i)  the statements  (a) and (b) are equivalent. The following  identity on quadratic algebras 
 \begin{eqnarray*}
 \langle a  ,  b  \rangle & = &  \frac{1}{2}  \big( \langle a + b  ,  a + b  \rangle - \langle a  ,  a  \rangle - \langle b  ,  b  \rangle \big),
  \end{eqnarray*}
 reveals that the statement (c) implies the statement (b). It thus follows that  (a),   (b), and (c) are equivalent.  

(iii) In view of part (i), it suffices to prove that  $f$ is an $F$-algebra isomorphism if  $f$ is multiplicative,  $ f(1_\mathbb{A}) = 1_\mathbb{B}$, and $ \langle f(a) , 1_\mathbb{B} \rangle_ \mathbb{B} =  \langle a , 1_\mathbb{A} \rangle_\mathbb{A}$ for all $ a \in \mathbb{A}$. From the hypothesis and the following identity on quadratic algebras, proved in Theorem \ref{02.7}, 
 \begin{eqnarray*}
 \langle a  ,  b  \rangle & = & 2  \langle a  ,  1  \rangle  \langle b ,  1  \rangle  -  \langle ab ,  1  \rangle  , \\
        & = & 2  \langle a  ,  1  \rangle  \langle b ,  1  \rangle  -  \langle ba  ,  1  \rangle  , 
  \end{eqnarray*}
  we see that  $f$ preserves  the bilinear form and $ f(1_\mathbb{A}) = 1_\mathbb{B}$. The assertion now follows from the multiplicativity  of $f$ and the second assertion of (i). This completes the proof.   
\hfill \qed

\bigskip

\begin{section}
{\bf Spectrum in one-sided alternative topological algebras}
\end{section}

\bigskip

\bigskip

Let   $\mathbb A$  be a real algebra. The {\it algebraic complexification} of $\mathbb A$, denoted by $\mathbb A_{\mathbb{C}}$, is the set of $\mathbb A \times \mathbb A$ together with addition, product by scalars, multiplication, and the conjugation operation  defined as follows
 \begin{eqnarray*}
 (a, b) + (c, d) & :=&  ( a + c, b + d),\\
 (r + is)(a, b) & :=& (ra - sb, sa + rb),\\
  (a, b) . (c, d) & :=& ( ac - bd, ad + bc), \\
\overline{(a, b)} & :=& (a, -b)
 \end{eqnarray*}
where $ a, b, c, d \in \mathbb A$ and $ r, s \in \mathbb{R}$. 
It is quite straightforward to check that  $\mathbb A_{\mathbb{C}}$ together with the aforementioned operations is a complex algebra equipped with the  conjugation, in particular 
$$ \overline{ (a, b) . (c, d)} =  \overline{ (a, b) } . \overline{  (c, d)},   \ \   \ (  a, b, c, d \in \mathbb A),$$
 and that, as a real algebra, it contains a copy of $\mathbb A$ via the monomorphism $ a \rightarrow (a, 0)$ ($ a \in \mathbb A$)   of real algebras.

 Also, it follows from  Proposition \ref{1.3}(i) that  $\mathbb A_{\mathbb{C}}$ is left (resp. right) alternative whenever $\mathbb A$ is.  It is not difficult to check that if $\mathbb A$ is a topological real algebra, then $\mathbb A_{\mathbb{C}}$ equipped with the product topology is a topological complex algebra.

\bigskip

If $p$ is a real vector space  seminorm on a real vector space   $X$, then $p$ gives rise to a complex vector space seminorm $p_{X_{\mathbb{C}}}: X_{\mathbb{C}}  \longrightarrow \mathbb{R}$ on  $ X_{\mathbb{C}}$, the complexification of $X$, defined by
$$ p_{X_{\mathbb{C}}} ( a, b) := \frac{1}{\sqrt{2}} \sup_{t \in \mathbb{R}} p_1 \big( (\cos t + i \sin t) ( a, b) \big), $$
 where $p_1$ is the real vector space   seminorm  on  $X_{\mathbb{C}}$ defined  by $ p_1 ( a, b) := p(a) + p(b)$. (See \cite[Theorem 1.3.1]{Ri}; note that the theorem and its proof should be rewritten for seminorms rather than for norms.) Likewise, one can rewrite  \cite[Theorem 1.3.2]{Ri} to show that any given real algebra seminorm on a real associative algebra $\mathbb A$   gives rise to a complex algebra seminorm on its complexification $ \mathbb A_{\mathbb{C}}$.
More generally, rewriting \cite[Proposition 13.1]{BD}, one can show that any algebra seminorm $p$ on a real algebra $\mathbb A$   gives rise to a complex algebra seminorm on its complexification $ \mathbb A_{\mathbb{C}}$.

 With these observations at our disposal, we see that the complexification  $\mathbb A_{\mathbb{C}}$ of a
 locally convex (resp. an lmc) real algebra  can naturally be viewed as  a locally convex (resp. an lmc) complex algebra.

\bigskip

Let  $\mathbb A$ be a complex topological algebra with identity. The spectrum of an element $ a \in \mathbb A$ is denoted by $ \sigma (a ; {\mathbb A} )$, or simply by $ \sigma(a)$, and is defined as follows
$$ \sigma(a) = \sigma (a ; {\mathbb A} ) := \{ z \in \mathbb{C}: z1 - a \notin \mathbb{A}^{-1}\}.$$

Let $\mathbb A$ be  a real or complex algebra lacking an identity element. The unitization of  $\mathbb A$, denoted by  $\mathbb A_1$, is the set $ \mathbb{F} \times \mathbb A$ together with 
 addition, product by scalars, and multiplication  defined as follows
 \begin{eqnarray*}
 (r, a) + (s, b) & :=&  ( r + s, a + b),\\
 r(s, b) & :=& (rs,  rb),\\
  (r, a) . (s, b) & :=& ( rs, sa + rb + ab), 
 \end{eqnarray*}
where $ a, b \in \mathbb A$ and $ r, s \in \mathbb{F}$. 
It is easily verified that  $\mathbb A_1$  contains a copy of $\mathbb A$ as an ideal of it via the monomorphism $ a \rightarrow (0, a)$ ($ a \in \mathbb A$)   of $\mathbb F$-algebras and that $ 1 = (1, 0)$ is the identity element of $\mathbb A_1$. It is also quite straightforward to check that $\mathbb A_1$  is left (resp. right) alternative (resp. power associative) whenever $\mathbb A$ is, and that $\mathbb A_1$  is  normed (resp. lmc) whenever $\mathbb A$ is. To see this latter fact, just note that any algebra (semi)norm $p$ on $\mathbb A$ gives rise to an algebra (semi)norm  $p_1$ on $\mathbb A_1$, which is defined as follows 
$$p_1( r , a ) := |r| + p(a), \ \ \   ( r \in \mathbb{F}, \ a \in \mathbb{A}).$$
By definition, 
$$ \sigma(a) = \sigma (a ; {\mathbb A} ) :=\sigma (a ; {\mathbb A}_1 ).$$

\bigskip

Motivated by the definition of the approximate point spectrum for bounded linear operators on complex Banach spaces, e.g., \cite[Definition VII.6.3]{C}, for an element $ a$ in a unital complex normed algebra   $( \mathbb A, \| .\|) $,  we use the symbol  $ \sigma_{lap}  (a ; {\mathbb A} )$, or simply $ \sigma_{lap} (a)$, to denote {\it the left   approximate spectrum} of $ a \in \mathbb A$, which we define by
\begin{eqnarray*}
\sigma_{lap} (a)  & := &    \Big\{ z \in \mathbb{C}:  \exists \ a_n \in \mathbb A,   ||a_n|| = 1 \  (n \in \mathbb{N}) , \lim_n  (z1 - a)   a_n = 0 \Big\}.
\end{eqnarray*}
The {\it  right  approximate spectrum} of $ a \in \mathbb A$, denoted by  $ \sigma_{rap}  (a ; {\mathbb A} )$ or simply by $ \sigma_{rap} (a)$, is defined  in a similar fashion. The {\it   approximate spectrum}  of $ a \in \mathbb A$, denoted by  $ \sigma_{ap}  (a ; {\mathbb A} )$ or simply by $ \sigma_{ap} (a)$, is defined as $ \sigma_{lap} (a) \cap  \sigma_{rap} (a)$. 

More generally, for a complex topological algebra with identity, {\it the left approximate spectrum} of $ a \in \mathbb A$ is defined by 
\begin{eqnarray*}
 \sigma_{lap} (a) & :=&    \Big\{ z \in \mathbb{C}:  \exists  \ 0 \in U \in \tau_\mathbb A,    a_n \in \mathbb A \setminus  U  \   (n \in \mathbb{N}),\\
 & & \ \    \lim_n (z1 - a) a_n  = 0  \Big\},
 \end{eqnarray*} 
where $\tau_\mathbb A$ denotes the topology of $\mathbb A$. It is not difficult to see that the two definitions coincide for complex normed algebras. Likewise, the notions of right  approximate spectrum and the approximate spectrum of an element $ a \in \mathbb A$ are naturally defined.

Let $\mathbb A$ be a real topological left (resp. right) alternative algebra with identity,   $ a \in \mathbb A$, and ${\mathbb A}_{\mathbb{C}}$ be its complexification. By definition
\begin{eqnarray*}
  \sigma(a)   =   \sigma (a ; {\mathbb A})  & :=&    \sigma(a; {\mathbb A}_\mathbb{C}) =   \{ z \in \mathbb{C}: z1 - a \notin {\mathbb A}_{\mathbb{C}}^{-1}\}, \\
   \sigma_{lap}(a)  =   \sigma_{lap} (a ; {\mathbb A}) & :=&  \sigma_{lap}(a; {\mathbb A}_\mathbb{C}) \\
 & =&  \big\{ z \in \mathbb{C}:  \exists  \ 0 \in U \in \tau_\mathbb A,    a_n \in {\mathbb A}_\mathbb{C} \setminus U \ (n \in \mathbb{N}),  \\
 & & \  \ \lim_n (z1 - a) a_n  = 0  \big\}.
\end{eqnarray*}
 The symbols $ \sigma_{rap}(a) $ and $\sigma_{ap}(a)$ are similarly defined. 
 Since $ a \in \mathbb A$, we get that   $z1 - a \in {\mathbb A}_{\mathbb{C}}^{-1}$ for a $  z \in \mathbb{C}$ if and only if  $\overline{z1 - a}   = \overline{z}1 - a \in {\mathbb A}_{\mathbb{C}}^{-1}$. This, in turn, yields  $z1 - a \in {\mathbb A}_{\mathbb{C}}^{-1}$ for a $  z \in \mathbb{C}$ if and only if
  $a^2 - 2 {\rm Re} (z) a + |z|^21 \in  {\mathbb A}^{-1}$. Therefore, 
  $$ \sigma(a) =  \sigma(a; {\mathbb A}_\mathbb{C}) = \{  z \in \mathbb{C} : a^2 - 2 {\rm Re} (z) a + |z|^21 \notin  {\mathbb A}^{-1} \}.$$
Likewise, one can see that 
 \begin{eqnarray*}
 \sigma_{lap}(a) & = &  \sigma_{lap}(a; {\mathbb A}_\mathbb{C})  =  \big\{  z \in \mathbb{C} :   \exists  \ 0 \in U \in \tau_\mathbb A,    a_n \in \mathbb A \setminus  U (n \in \mathbb{N}), \\
 & & \  \   \ \ \ \ \ \ \ \ \ \ \ \ \ \ \ \ \ \ \lim_n  (a^2 - 2 {\rm Re} (z) a + |z|^21) a_n = 0 \big\}.
   \end{eqnarray*}
Obviously, similar relations hold for $ \sigma_{rap}(a) $ and $\sigma_{ap}(a)$.  

\bigskip

At this point, we make a  useful, simple and crucial observation.  Proposition \ref{1.3}(iv) clearly implies that for invertible elements  $a, b$ of a left (resp. right) alternative algebra  $\mathbb A$, we have 
 \begin{eqnarray*}
  \  \ \ \ \ \ \ \ \ a^{-1} - b^{-1} & = &  a^{-1} ( b b^{-1}) - (a^{-1} a) b^{-1},\\
  &  = &    a^{-1} ( b b^{-1}) - a^{-1} ( a b^{-1}), \\
  &  = &  a^{-1}\big(  ( b-a) b^{-1} \big), \\
\Big({\rm resp.} \  a^{-1} - b^{-1} & = &   (  a^{-1} b)b^{-1} -  (a a^{-1} )b^{-1},\\
  &  = &    (  a^{-1}b -  a^{-1}a ) b^{-1}, \\
  &  = & \big(a^{-1}  ( b-a)  \big) b^{-1}.  \Big)
   \end{eqnarray*}

With this observation at one's disposal, one can develop the spectral theory in the setting of one-sided alternative algebras. 
For instance, with a bit of care, one can prove the counterparts of most of the standard results of Section 1 and some of the results of Sections 2 and 3 of Chapter 1 of \cite{GP}  as well as some of the results Chapter 1 of \cite{BD}, e.g., holomorphic functional calculus, differentiability of the inversion function, and the Gelfand-Beurling formula for the spectral radius to name a few, in the setting of one-sided alternative (complete) normed algebras. For the counterparts of  several standard results in the setting of continuous inverse algebras, we refer the reader to \cite{G}. We freely make use of the standard theorems in the setting of (complete) one-sided  alternative normed algebras when needed. 
As an example, let's prove the existence of the spectrum in the setting of one-sided  alternative complex topological unital algebras. 

\bigskip

\begin{thm} \label{1.4} 
Let $\mathbb A$ be a left (resp. right) alternative complex topological  algebra with identity whose (topological) dual separates its points, e.g. a  complex  left (resp. right) alternative unital lmc  algebra or a complex  left (resp. right) alternative continuous inverse  algebra. Then, $ \sigma (a)\not= \emptyset$  for all $ a \in \mathbb A$.
\end{thm}

\bigskip

\noindent {\bf Proof.}  Proceed by contradiction and suppose that $ \sigma (a) = \emptyset$ so that $ (x + i y) 1 - a \in \mathbb A^{-1}$ for all $ x, y \in \mathbb{R}$.  Fix an arbitrary $ \phi \in  \mathbb A^*$ and define $ f : \mathbb{R}^2 \longrightarrow \mathbb{C}$ by
$ f (x, y) = \phi \Big(  \big( (x + i y) 1 - a \big)^{-1}  \Big) $.
Let $ f = u + i v$, where $ u = {\rm Re}(f), v = {\rm Im}(f)$. In view of the exhibited identities preceding the theorem, it is pretty straightforward to verify that
$$\frac{\partial f}{\partial x}= \phi \Big(  - \big((x + i y) 1 - a \big)^{-2}  \Big) , \  \frac{\partial^2 f}{\partial x^2}= \phi \Big( 2 \big((x + i y) 1 - a \big)^{-3}  \Big), $$
$$\frac{\partial f}{\partial y}= \phi \Big( -i \big((x + i y) 1 - a \big)^{-2}  \Big) , \  \frac{\partial^2 f}{\partial y^2}= \phi \Big( -2 \big((x + i y) 1 - a \big)^{-3}  \Big), $$
from which we obtain
$ \frac{\partial^2 f}{\partial x^2} + \frac{\partial^2 f}{\partial y^2}= 0$, and hence $ \frac{\partial^2 u}{\partial x^2} + \frac{\partial^2 u}{\partial y^2}= 0$ and $ \frac{\partial^2 v}{\partial x^2} + \frac{\partial^2 v}{\partial y^2}= 0$. On the other hand, we have
$$ \lim_{(x, y ) \rightarrow \infty} |f(x, y)|  =  \lim_{(x, y ) \rightarrow \infty} \frac{1}{|x + i y|} \Big| \phi \Big( \big(   1 - \frac{1}{x + i y}a \big)^{-1} \Big) \Big| =0 \times  |\phi(1)| = 0 ,  $$
implying that
 $$ \lim_{(x, y ) \rightarrow \infty} |u(x, y)|= 0 =  \lim_{(x, y ) \rightarrow \infty} |v(x, y)|,$$
and hence $f, u , v$ are all bounded functions.
Therefore, $ u, v :\mathbb{R}^2 \longrightarrow \mathbb{R}$ are zero functions because they are bounded harmonic functions, and hence constant, and whose limits  at infinity are zero.  This, in turn, yields  $ f (x, y) = \phi \Big(  \big((x + i y) 1 - a \big)^{-1}  \Big)=0 $
for all  $ x , y \in \mathbb{R}$ and $ \phi \in  A^*$, from which we get that $ \big(  (x + i y) 1 - a \big)^{-1} =0 $ for all  $ x , y \in \mathbb{R}$, which is absurd. This completes the proof.  
\hfill \qed

\bigskip

\noindent {\bf Remarks.} 1.  In the theorem, one can drop the hypothesis that the complex topological algebra $\mathbb A$  is one-sided alternative. However to compensate for this, one can assume that for any  $ a \in \mathbb A$, the function  $ (z 1 - a)^{-1}$ is weakly differentiable on its domain, meaning that the complex-valued function $\phi\big( (z 1 - a)^{-1} \big) $  is differentiable on its domain, which is a subset of complex numbers,  for all  $\phi \in \mathbb A^*$.

2. The following question, with an affirmative answer in the setting of associative algebras, is an example of a standard result that holds in the setting of associative algebras. However, we do not know if its counterpart holds in the setting of alternative algebras. {\it Let $\mathbb A$ be as in the theorem and $ a, b \in \mathbb A$. Is it true that $ \sigma (ab) \cup \{ 0\} =  \sigma (ba) \cup \{ 0\}$?} The question has an affirmative answer if the algebra $\mathbb A$ is algebraic.

 \bigskip

 An application of the existence of spectrum is the existence of nontrivial hyperinvariant linear manifolds for continuous linear operators acting on arbitrary Banach spaces.   Part (i) of the corollary below is the counterpart of  \cite[Proposition 8.1]{RR} in the setting of complex Hilbert spaces. It is also worth mentioning that the following corollary can be rephrased in the language of {\it strictly transitive algebras}    in the context of arbitrary Banach spaces; see \cite[Section 8.1]{RR}. For a standard reference on invariant subspaces, see \cite{RR}.

\bigskip

\begin{cor} \label{3.2} 
 {\rm (i)} Let $X$ be a complex Banach space and $ T \in \mathcal{B} (X)$ be nonscalar. Then,  $ T $ has a nontrivial hyperinvariant linear manifold, i.e., there exists a nontrivial linear manifold in  $X$ that is invariant under all continuous linear operators that commute with $T$.

{\rm (ii)}  Let $X$ be a real Banach space and $ T \in \mathcal{B} (X)$ be nonquadratic, i.e., $ T$ does not have  a vanishing  polynomial with real coefficients of degree at most two. Then,  $ T $ has a nontrivial hyperinvariant linear manifold i.e., there exists a nontrivial linear manifold in  $X$ that is invariant under all continuous linear operators that commute with $T$. 
\end{cor}

\bigskip

\noindent {\bf Proof.} (i) By  the preceding theorem that $ \sigma \big(T;  \mathcal{B} (X) \big) \not= \emptyset$. Now, pick a $ \lambda \in  \sigma \big(T;  \mathcal{B}(X) \big)$ and set $ S := T - \lambda I_X$. We get  that $ 0 \not= S \notin   \mathcal{B} (X)^{-1}$. But $X$ is a  Banach space. So, in view of the Open Mapping Theorem, \cite[Theorem 2.11]{Ru}, we see that $S$ is not injective or it is not surjective. Therefore, one of the linear manifolds $ \ker S$ or $ SX$ will be  a nontrivial  hyperinvariant linear manifold  for $T$ depending on whether or not  $S$ is not injective, as desired.   

 (ii) The proof is almost identical to that of (i) except that one should set $ S := T^2 - 2 {\rm Re}( \lambda) T + |\lambda|^2 I_X$, where $ \lambda \in  \sigma \big(T;  \mathcal{B} (X) \big)$. This completes the proof. 
\hfill \qed 

\bigskip

\noindent {\bf Remarks.} 1. The counterpart of part (ii) of the corollary  for nonquadratic continuous left (resp. right) linear operators acting on left (resp. 
right) quaternion  Banach spaces remains valid. We will deal with this and other results in a forthcoming note based on a work in progress.

2. With the corollary and a conjecture of Victor Lomonosov, \cite[page 338]{L}, at disposal, one is tempted to suggest the following {\it hyperinvariant subspace problem conjecture}.  {\it Let $X$ be a complex (resp. real or a left  or right quaternion) Banach space and $ T \in \mathcal{B} (X)$  be nonscalar (resp. nonquadratic). Then, $ T^* \in  \mathcal{B} (X^*)$  has a nontrivial hyperinvariant subspace. }

\bigskip

\begin{section}
{\bf On algebras over the fields of real and complex numbers}
\end{section}

\bigskip

The following proposition, which might be of independent interest, turns out to be useful. 

\bigskip

\begin{prop} \label{1.6} 
   Let $( \mathbb A, |.|)$ be an absolute-valued  algebra and  let  $\mathcal{B} (\mathbb A) $ be equipped with $\|.\|$, the operator norm of it induced by $|.|$. 
   
{\rm (i)} If $L_a , R_a \in \mathcal{B} (\mathbb A)^{-1}$ for some $ a \in \mathbb A$, then $\| L_a^{-1}\| =  |a|^{-1}  = \| R_a^{-1}\|$.  
   
{\rm (ii)}    The mappings  $L, R : \mathbb A  \longrightarrow \mathcal{B} (\mathbb A) $ defined by $L(a) = L_a$ and $R(a) = R_a$ are isometric linear mappings. Moreover, for each $n \in \mathbb{N}$, if $ a_i  \in \mathbb A$ is such that $L_{a_i} , R_{a_i} \in \mathcal{B} (\mathbb A)^{-1}$  ($ 1 \leq i \leq n$), then 
$$  \|S_1 \cdots S_n \| = \|S_1 \cdots S_{n-1} \| \| S_n \| = \|S_1 \| \cdots \|S_n\|  , $$
where $ S_i \in \{ L_{a_i}, L_{a_i}^{-1},   R_{a_i}, R_{a_i}^{-1}\}$  ($ 1 \leq i \leq n$). Also, if $ S_i \in \{ L_{a_i},  R_{a_i}\}$  ($ 1 \leq i \leq n$), the hypothesis that 
$L_{a_i} , R_{a_i} \in \mathcal{B} (\mathbb A)^{-1}$ ($ 1 \leq i \leq n$) is redundant. 
\end{prop}

\bigskip

\noindent {\bf Proof.} (i) Since $L_a  \in \mathcal{B} (\mathbb A)^{-1}$, we can write 
\begin{eqnarray*}
\| L_a^{-1}\|  & =&  \sup_{|x| = 1} | L_a^{-1} (x ) | = \sup_{|L_a y| = 1} | L_a^{-1} (L_a y )|,\\
&  = &  \sup_{|y| = \frac{1}{|a|}} |  y|  = |a|^{-1},  
\end{eqnarray*}
 as desired. Likewise, $\| R_a^{-1}\|=  |a|^{-1}  $. 

(ii) That the mappings $L$ and $ R $  are isometric linear mappings is easy to check.  For the rest, it suffices to prove that
$$  \|S_1 \cdots S_n \| = \|S_1 \cdots S_{n-1} \| \| S_n \|,  $$
 for all $n \in \mathbb{N}$. To this end, first suppose  $ S_n  = L_{a_n}$ for some $ a_n  \in \mathbb A$. Since $L_{a_n}  \in \mathcal{B} (\mathbb A)^{-1}$, we have 
$$ \{ a_n x : x \in \mathbb A,  |x| = 1 \} = \{|a_n | y : y  \in \mathbb A,  |y| = 1 \} .$$
So we can write 
\begin{eqnarray*} 
 \| S_1 \cdots S_n \|  & =&  \sup_{|x| = 1} |S_1 \cdots S_n (x)| = \sup_{|x| = 1} |S_1 \cdots S_{n-1} (a_nx)|, \\
 &  = & \sup_{|y| =1} |S_1 \cdots S_{n-1} ( |a_n| y)| = |a_n| \sup_{|y| =1} |S_1 \cdots S_{n-1} ( y)|, \\
 & =&  \|S_1 \cdots S_{n-1} \| \| S_n \| . 
\end{eqnarray*}
In the three other cases, namely, $ S_n \in \{  L_{a_n}^{-1},   R_{a_n}, R_{a_n}^{-1}\}$, the assertion is proved in a similar fashion. The last assertion is trivial, for $|.|$ is  an absolute value on $ \mathbb A$.  This completes the proof. 
\hfill \qed

\bigskip

Part (i) of the following proposition is known by the experts. Part (ii)-(a) of it is motivated by \cite[Corollary 2.7.29]{GP}.

\bigskip

\begin{prop} \label{1.7} 
{\rm (i)}   Let $\mathbb A$ be a normed  algebra. Then the norm of $\mathbb A$ restricted to any subspace $ \mathcal{M}$ of $\mathbb A$ on which the multiplication is commutative and the norm is multiplicative comes from an inner product. In particular, the absolute value of any commutative absolute-valued algebra comes from an inner product on the algebra.

\bigskip

{\rm (ii)}  The absolute value of any absolute-valued algebra $\mathbb A$ comes from an inner product on the algebra if one of the following conditions holds.

{\rm (a)} There exists an element  $a_0 \in \mathbb A$  such that $ L_{a_0} $ or $ R_{a_0} $ has dense range, for instance, if $ \{L_{a_0},  R_{a_0}\} \cap \mathcal{B}(\mathbb A)^{-1} \not= \emptyset$, or   $\mathbb A$ is left or right unital, or  $\mathbb A$ is finite-dimensional;

 {\rm (b)} the  algebra $\mathbb A$  is power-associative. 

\end{prop}

\bigskip

\noindent {\bf Proof.}(i) Let $\mathbb A$ be an arbitrary normed algebra equipped with a norm algebra $\|.\|$ and $ \mathcal{M}$ be a subspace  of $\mathbb A$ on which the multiplication is commutative and the norm is multiplicative. In view of a well-known result of Schoenberg, \cite[Theorem 1]{Scho}, it suffices to show that $\|.\|$ satisfies Ptolemy's inequality on $ \mathcal{M}$, namely$$ \|a+ b \| \|a + c\| \leq   \|b\| \|c\|  + \|a\| \|a + b + c\|, $$
for all $ a, b, c \in \mathcal{M}$. But this is obvious because 
\begin{eqnarray*}
 \|a+ b \| \|a + c\| & =& \| ( a+ b) (a + c) \|, \\
 & \leq & \| bc \| + \| a^2 + ac + ba\|, \\
 & =& \| bc \| + \| a^2 + ac + ab\|, \\
 & = &   \|b\| \|c\|  + \|a\| \|a + b + c\|, 
\end{eqnarray*}
for all $ a, b, c \in \mathcal{M}$. This completes the proof. 

(ii) If the algebra is complex, indeed the assertion holds under the weaker hypothesis that the given algebra  $\mathbb A$ is nearly absolute-valued, for by Theorem \ref{1.8}(iv) the algebra   $\mathbb A$ is isomorphic to  $\mathbb C$.  It remains to prove the assertion for real algebras. To this end, let $\mathbb A$ be an absolute-valued real algebra. Clearly, the completion of $\mathbb A$, denoted by  $ \overline{\mathbb A}$, is also absolute-valued. We prove the assertion by showing that the absolute value of  $ \overline{\mathbb A}$ comes from an inner product. 
First, suppose that $ L_{a_0} $ has dense range for some $a_0 \in \mathbb A$.
Since $L_{a_0}: \mathbb A \longrightarrow \mathbb A  $ is a bounded linear operator, it can uniquely be extended to $ \overline{\mathbb A}$; which is $L_{a_0}$ but viewed as the left multiplication by $ a_0$ on $ \overline{\mathbb A}$, i.e., $L_{a_0} a = a_0 a $ ($a \in \overline{\mathbb A}$). We see that $L_{a_0}: \overline{\mathbb A} \longrightarrow \overline{\mathbb A}  $  is surjective because $L_{a_0}: \mathbb A \longrightarrow \mathbb A  $ has dense range. Consequently,  $L_{a_0} \in \mathcal{B}(\overline{\mathbb A} )^{-1} $, for  $ \overline{\mathbb A}$  is complete and $L_{a_0}$ is bijective. 
 It thus follows from \cite[Lemma  2.7.18]{GP} that $L_{a} \in \mathcal{B}(\overline{\mathbb A} )^{-1} $ for all $ a \in \mathbb{A} \setminus \{0\}$. Let $ a, b \in \mathbb{A} $ with $\{ a, b\}$ linearly independent be arbitrarily given. Set $ S : =L_{b} L_{a}^{-1} \in \mathcal{B}(\overline{\mathbb A} )^{-1} $. Clearly, the subalgebra  $ \mathbb{S}$ generated by $ \{ S , I_{\overline{\mathbb A}} \}$ in $\mathcal{B}(\overline{\mathbb A} )$  is a commutative algebra. By showing that the operator  norm of $\mathcal{B}(\overline{\mathbb A} )$ is multiplicative on $ \mathcal{M}$, the subspaces spanned by  $ \{ S , I_{\overline{\mathbb A}} \}$, we get that the absolute value of the algebra  $\mathbb A$, say $|.|$,  comes from an inner product. To this end, in view of Proposition \ref{1.6}(ii), for all $ c, d, e , f \in \mathbb{F}$ with $(c, d) \not= 0$ and $(e , f) \not= 0$, we can write 
\begin{eqnarray*}
\|\big( c S + d I_{\overline{\mathbb A}} \big) \big(e S + f I_{\overline{\mathbb A}} \big)\| & = & \| \big( c L_{b}  + d L_{a}\big) L_{a}^{-1}\big( e L_{b}  + fL_{a}\big) L_{a}^{-1}  \| , \\
& =& \|  L_{cb + da}   L_{a}^{-1}  L_{eb + fa}  L_{a}^{-1}  \| , \\
& =& \|  L_{cb + da} \| \| L_{a}^{-1}\| \| L_{eb + fa} \| L_{a}^{-1}  \| , \\
& =& \|  L_{cb + da}   L_{a}^{-1} \| \|L_{eb + fa}  L_{a}^{-1}  \| , \\
& =& \| (c L_{b}  + d L_{a} ) L_{a}^{-1}\| \| (e L_{b}  + f L_{a}) L_{a}^{-1}\|, \\
& =&  \| c S + d I_{\overline{\mathbb A}} \| \|e S + f I_{\overline{\mathbb A}} \|.
\end{eqnarray*}
Consequently, 
$$ \| S + I_{\overline{\mathbb A}} \|^2 +   \| S - I_{\overline{\mathbb A}} \|^2 = 2 \big( \| S \|^2 + \|I_{\overline{\mathbb A}} \|^2 \big) , $$
from which, simplifying just  as done in the above paragraph, we see that 
$$\| ( L_{b}  +  L_{a} ) L_{a}^{-1}\|^2  + \| (L_{b}  -  L_{a}) L_{a}^{-1}\|^2 =  2 \big( \|L_{b} L_{a}^{-1}  \|^2 +  \|L_{a} L_{a}^{-1}\|^2 \big) , $$
 equivalently, 
$$\|  L_{a + b} \|^2 \| L_{a}^{-1}\|^2  + \| L_{b-a} \|^2  \|L_{a}^{-1}\|^2 =  2 \big( \|L_{b}\|^2 \| L_{a}^{-1}  \|^2 +  \|L_{a}\|^2 \| L_{a}^{-1}\|^2 \big) , $$
 and hence, in view of the preceding proposition, we obtain
$$ | b+ a|^2 + |b - a|^2 = 2 \big( |b|^2 + |a|^2\big).$$ 
Therefore, the parallelogram identity holds for $|.|$. The assertion now follows from a well-known theorem of Jordan and von Neumann, \cite[Theorem 1]{JvN}. 
 
  Next, suppose that (b) holds. We prove the assertion by showing that   $\mathbb A$ is unital, which implies that  (a) holds, completing the proof. Let $ a \in \mathbb A$ be arbitrarily given. Set $ \mathbb{B}:= {\rm Alg}(a)$, the real subalgebra  generated by $a$, let $ \overline{ \mathbb{B}}$ denote its completion,  and  $\overline{ \mathbb{B}}_\mathbb{C}$ be a complexification of $ \overline{ \mathbb{B}}$. Consider $ L_a$ as an element of $ \mathcal{B} (  \overline{ \mathbb{B}}_\mathbb{C})$. By \cite[Proposition VII.6.7]{C}, $ \sigma_{ap}\big(  L_a; \mathcal{B} (  \overline{ \mathbb{B}}_\mathbb{C}) \big) \not= \emptyset$. Thus, there exist $ \lambda \in \mathbb{C}$ and  
  $ \bar a_n + i \bar b_n \in \overline{ \mathbb{B}}_\mathbb{C}$ with $\| \bar a_n + i \bar b_n \| = 1$ ($n \in \mathbb{N}$)  such that  
 $$ \lim_n (L_a  - \lambda I_{\overline{ \mathbb{B}}_\mathbb{C}}) (   \bar a_n + i \bar b_n) = 0, $$
 implying that
 $$ \lim_n  a (L_a  - \overline{\lambda} I_{\overline{ \mathbb{B}}_\mathbb{C}})  (L_a  - \lambda I_{\overline{ \mathbb{B}}_\mathbb{C}}) (   \bar a_n + i \bar b_n) = 0, $$
 and hence
\begin{eqnarray*}
\lim_n ( a^3 - 2 {\rm Re}(\lambda)  a^2 + |\lambda|^2 a )    \bar a_n  & = & 0 , \\
\lim_n ( a^3 - 2 {\rm Re} (\lambda)  a^2 + |\lambda|^2 a )   \bar  b_n  & = & 0 .
\end{eqnarray*}
This yields $  a^3 - 2 {\rm Re} (\lambda)  a^2 + |\lambda|^2 a= 0$, for $\overline{ \mathbb{B}}$ is absolute-valued and $\| \bar a_n + i \bar b_n \| = 1$ ($n \in \mathbb{N}$).
Therefore,  $\mathbb A$  is a power-associative algebraic algebra. It thus follows from \cite[Lemma 2.5.5]{GP} or Proposition \ref{1.1}(iii) that $\mathbb A$ is unital, and hence locally complex. From this point on one may argue that the assertion now follows from (a), or that the assertion follows from Theorem \ref{2.13}.
This completes the proof.  
\hfill \qed

\bigskip

Here are some  Gelfand-Mazur type theorems in several settings; see \cite[Proposition 2.5.49(i)]{GP}, \cite[Corollary 2.5.54(i)]{GP}, \cite[Corollary 2.5.56(i)]{GP}, \cite[Lemma 2.6.34]{GP}, \cite[Corollary 2.7.3]{GP},  \cite[Theorem 2.7.16]{GP},  \cite[Theorem 2.7.17]{GP}, and \cite[Characterization 2]{Za}. Note that in parts (iii) and (iv) of the theorem the algebra $\mathbb A$ is not necessarily normed; the algebra $\mathbb A$ is equipped with a vector space norm with respect to which the  multiplication operation of the algebra $\mathbb A$ is separately continuous.

\bigskip 

\begin{thm} \label{1.8}

\bigskip

 {\rm (i)}   Let $\mathbb A$ be a  left (resp. right) alternative complex topological  algebra  whose (topological) dual separates its points. If $\mathbb A$  is a quasi-division algebra, equivalently, a classical division algebra, then it is isomorphic to $\mathbb C$.

\bigskip

{\rm (ii)}  Let $\mathbb A$ be a complex algebra endowed with a vector space norm such that $ {\rm l.t.z.d.} (\mathbb{A})  = \{0\}$ (resp. $ {\rm r.t.z.d.}(\mathbb{A}) = \{0\}$). Assume further that the  multiplication is separately continuous. Then, $\mathbb A$  is isomorphic to $\mathbb C$ provided that  $ L_{a_0} $  (resp. $ R_{a_0} $)  has dense range for some $ a_0  \in \mathbb{A}$. 

\bigskip

{\rm (iii)}  Let $\mathbb A$ be a complex  algebra endowed with a vector space norm such that the algebra generated by any nonzero element of $\mathbb{A}$ has no nontrivial topological zero-divisors. Assume further that the  multiplication is separately continuous. Then, $\mathbb A$  is isomorphic to $\mathbb C$ if $\mathbb A$   is power-associative. In particular, every   alternative quasi-division  complex normed algebra is isomorphic to $\mathbb C$. 

\bigskip

{\rm (iv)} 
Let $\mathbb A$ be  a complex  algebra endowed with a $c$-supermultiplicative vector space norm for some $ c > 0$ relative to which the  multiplication is separately continuous, e.g., a complex nearly absolute value. Then, the complex algebra $\mathbb A$  is isomorphic to $\mathbb C$ if one of the  following conditions holds. 

{\rm (a)} There exists an element  $a_0 \in \mathbb A$  such that $ L_{a_0} $ or $ R_{a_0} $ has dense range, for instance if $ \{L_{a_0},  R_{a_0}\} \cap \mathcal{B}(\mathbb A)^{-1} \not= \emptyset$,   or  $\mathbb A$ is left or right unital, or  $\mathbb A$ is finite-dimensional; 

 {\rm (b)} the complex algebra $\mathbb A$  is power-associative.

 {\rm (c)} the complex algebra $\mathbb A$  is algebraic.

\bigskip

{\rm (v)} 
Let $\mathbb A$ be an absolute-valued complex algebra. If  {\rm (a)},  {\rm (b)}, or  {\rm (c)} in  {\rm (iv)}  holds, then $\mathbb A$  is isometrically isomorphic to $\mathbb C$. 

\bigskip

{\rm (vi)} Let $\mathbb A$ be a power-asociative complex  algebra endowed with an inner-product $ \langle . , . \rangle$, whose norm satisfies the identity $\|a^2\| = \|a\|^2 $ on $\mathbb A$. Then, $\mathbb A$  is isometrically isomorphic to $\mathbb C$. 

\bigskip

{\rm(vii)} Let $\mathbb A$ be a power-asociative  complex  unital algebra endowed with an inner-product $ \langle . , . \rangle$, whose norm is a unital algebra norm on $\mathbb A$. Then, $\mathbb A$  is isometrically isomorphic to $\mathbb C$. 

\end{thm}

\bigskip

\noindent {\bf Proof.} (i)  Since $\mathbb A$ is quasi-division, we see from Proposition \ref{1.1}(vi) that the algebra $\mathbb A$ is a classical division algebra, and, in particular, unital. Now, by Theorem \ref{1.4}, every $ a \in \mathbb A$  has a nonempty spectrum, implying that $a = \lambda 1 $ for some $ \lambda \in \mathbb{C}$. In other words, $\mathbb A$ is one-dimensional. This proves the assertion. 

(ii) We prove the assertion for algebras $\mathbb A$ having no nontrivial  left  topological zero-divisors. The assertion is proved similarly if the algebra $\mathbb A$ has no nontrivial  right topological zero-divisors.   First, suppose that $ L_{a_0} $  has dense range for some $ a_0  \in \mathbb{A}$. Use  $ \overline{\mathbb A}$ to denote the completion of  $\mathbb A$.  The linear operator $L_{a_0}: \overline{\mathbb A}  \longrightarrow \overline{\mathbb A}  $ is indeed the unique extension of the bounded linear operator $L_{a_0}: \mathbb A \longrightarrow \mathbb A  $ to $ \overline{\mathbb A}$; note that this linear operator is bounded  since the multiplication of $\mathbb A$ is separately continuous. Also note  that $L_{a_0}: \overline{\mathbb A}  \longrightarrow \overline{\mathbb A}  $ is continuous, for the multiplication of  $ \overline{\mathbb A}$ is separately continuous since  such is that of $\mathbb A$.  It thus follows that   $L_{a_0} \in \mathcal{B}(\overline{\mathbb A} )^{-1}$ because   $ \overline{\mathbb A}$ is complete and  $L_{a_0}: \overline{\mathbb A}  \longrightarrow \overline{\mathbb A}  $ is bijective, for  $\mathbb A$  has no nontrivial left  topological zero-divisors and that $L_{a_0}: \mathbb A \longrightarrow \mathbb A  $ has dense range. Let $ b \in \mathbb A$ be arbitrary. By  \cite[Proposition VII.6.7]{C}, $ \sigma_{ap}\big(L_{a_0}^{-1} L_{b} ;  \mathcal{B}(\overline{\mathbb A} ) \big) \not= \emptyset $. It follows that there exist a $ \lambda \in \mathbb{C}$ and $ \overline{a}_n \in  \overline{\mathbb A}$ with $\|\overline{a}_n\| =1 $ ($n \in \mathbb{N}$) such that $ \lim_n (L_{a_0}^{-1}  L_{b} - \lambda I_{\overline{\mathbb A}}) \overline{a}_n = 0$, equivalently, $ \lim_n (b - \lambda a_0) \overline{a}_n  = 0$. But $\mathbb A$ is dense in $ \overline{\mathbb A} $. So, we get that there exist $a_n \in  \mathbb A$ with $\|a_n\| =1 $ ($n \in \mathbb{N}$) such that   $ \lim_n (b - \lambda a_0) a_n  = 0$. This implies that $ b = \lambda a_0$ because $\mathbb A$ has no nontrivial left  topological zero-divisors. In other words,  $\mathbb A$  is one-dimensional, from which the assertion follows. 

(iii)  Given a nonzero $a \in \mathbb{A}$, let   $\mathbb{B}$ be the subalgebra generated by $a$, and $ \overline{\mathbb{B}}$ be its completion, which is  commutative and associative. Just as argued  in part (iii),  $L_{a}: \overline{\mathbb B}  \longrightarrow \overline{\mathbb B}  $ is continuous. Then again, by  \cite[Proposition VII.6.7]{C}, $ \sigma_{ap}\big( L_a ;  \mathcal{B}(\overline{\mathbb{B}} ) \big) \not= \emptyset $. So we conclude that there are elements $ \lambda \in \mathbb{C}$ and $ \bar a_n \in  \overline{\mathbb{B}}$ with $\|\bar a_n\| =1 $ ($n \in \mathbb{N}$) such that $ \lim_n ( L_a - \lambda I_{\overline{\mathbb{B}}}) \bar a_n = 0$, equivalently, $ \lim_n (a\bar a_n - \lambda \bar a_n)  = 0$, and hence 
$ \lim_n (a^2 - \lambda a) \bar a_n  = 0$. Since $\mathbb{B}$  is dense in $ \overline{\mathbb{B}}$, for any $ n \in \mathbb{N}$, we may pick $ a'_n \in \mathbb{B}$ such that $ \| a'_n - \bar a_n \| < \frac{1}{n}$.  It thus follows that  $ \lim_n (a^2 - \lambda a)  a'_n  = 0$. If necessary, by passing to a proper subsequence, we may assume that $ \lim_n \|a'_n \| = c > 0$. From this we get that  $\lim_n  a_n (a^2 - \lambda a)  = \lim_n (a^2 - \lambda a)  a_n  = 0$, where $ a_n = \frac{a'_n}{\|a'_n \|}$. This implies that $ a^2 = \lambda a$ because $ \|a_n \| = 1$ for all $ n \in \mathbb{N}$ and  $ \mathbb{B}$ has no nontrivial  topological zero divisors.  The assertion now follows from Theorem \ref{1.2}(i).

(iv)  The first assertion, namely, part (a),  follows from (ii).

As for the second assertion, given an $a \in \mathbb{A}$, let   $\mathbb{B}$ be the subalgebra generated by $a$, and $ \overline{\mathbb{B}}$ be its completion, which is  associative and is endowed with the completion norm. The associative complex algebra  $ \overline{\mathbb{B}}$ equipped with the completion norm is  easily seen to be $c$-supermultiplicative and relative to which the multiplication of $ \overline{\mathbb{B}}$ is separately continuous.  Now, by  \cite[Proposition VII.6.7]{C}, we have $ \sigma_{ap}\big( L_a ;  \mathcal{B}(\overline{\mathbb{B}} ) \big) \not= \emptyset $. So we conclude that there are elements $ \lambda \in \mathbb{C}$ and $ \bar a_n \in  \overline{\mathbb{B}}$ with $\|\bar a_n\| =1 $ ($n \in \mathbb{N}$) such that $ \lim_n ( L_a - \lambda I_{\overline{\mathbb{B}}}) \bar a_n = 0$, equivalently, $ \lim_n (a\bar a_n - \lambda \bar a_n)  = 0$, and hence
$ \lim_n (a^2 - \lambda a) \bar a_n  = 0$.
From this and
$$ c||a^2 - \lambda a ||  = c ||a^2 - \lambda a || ||\bar a_n|| \leq  ||(a^2 - \lambda a) \bar a_n|| , ( n \in \mathbb{N}) ,$$
we get that
$$ c||a^2 - \lambda a ||   \leq \lim_n ||(a^2 - \lambda a) \bar a_n||=0 ,$$
 implying that $ a^2 = \lambda a$.  The assertion now follows from Theorem \ref{1.2}(i).
 
 Finally, if  the algebra  $ \mathbb{A}$  is algebraic, the assertion follows from Theorem \ref{1.2}(ii)  because $ \mathbb{A}$ has no nonzero left zero-divisors.

(v) This is a quick consequence of (iv).

(vi) We prove the assertion by showing  that $ \mathbb{A}$ together with its addition and its Jordan multiplication and its inner product is isometrically isomorphic to $ \mathbb{C}$. Recall that  the symbol $ \mathbb{A}^{{\rm sym}}$ is used to denote the algebra $ \mathbb{A}$ together with its vector space operations and its Jordan product.  It is easily checked that the algebra $ \mathbb{A}^{{\rm sym}}$ is  power-associative and, clearly, commutative.  Now, for all $ a, b \in \mathbb{A}$, we can write 
\begin{eqnarray*}
\|a\|^2 + \|b\|^2 + 2 {\rm Re} \langle a, b \rangle & = & \| a + b \|^2 ,\\
& = & \| a^2  + b^2 + 2 a \circ b  \|, \\
& \leq & \| a\|^2  +\| b \|^2+ 2 \| a \circ b  \|, 
\end{eqnarray*}
from which we obtain  $ {\rm Re}\langle a, b \rangle \leq \| a \circ b  \|$. If necessary, replacing $a$ by $ua$, where $ u \in \mathbb{C}$ is such that  $|u| = 1$ and $ u\langle a, b \rangle = |\langle a, b \rangle |$,  we get that  $ | \langle a, b \rangle | \leq  \| a \circ b  \|$ for all $ a, b \in \mathbb{A}$. On the other hand, 
in view of the parallelogram identity, for all $ a, b \in \mathbb{A}$ with $ \| a \| = \|b \| =1$, we have
\begin{eqnarray*}
 4 \| a \circ b  \| & = &  \| 2 (ab + ba) \|, \\
 & = & \| ( a + b)^2 - (a - b)^2 \|, \\
 & \leq & \|  a + b\|^2  +  \| a - b\|^2 , \\
  & = & 2 ( \|a\|^2 + \| b\|^2 ) = 4.
\end{eqnarray*}
This yields $  \| a \circ b  \|  \leq \| a \| \| b \|$ for all $ a, b \in \mathbb{A}$. Thus, 
$$  | \langle a, b \rangle | \leq \| a \circ b  \|  \leq \| a \| \| b \|$$
for all $ a, b \in \mathbb{A}$. Next, we complete the proof by showing that $ \mathbb{A}^{{\rm sym}}$, and hence $\mathbb{A}$ is one-dimensional. To this end, pick a nonzero $ a \in  \mathbb{A}^{{\rm sym}}$ and form the algebra generated by $ a$ in $ \mathbb{A}^{{\rm sym}}$ and denote it by $\mathbb{B}$. We note that  $\mathbb{B}$ is a commutative and associative complex normed algebra equipped with the inner-product $ \langle . , .  \rangle$ satisfying the relations 
 $  | \langle a, b \rangle | \leq \| a \circ b  \|  \leq \| a \| \| b \|$ for all $ a, b \in \mathbb{B}$.  We claim that the algebra $\mathbb{B}$ has non nontrivial topological zero-divisors. Suppose on the contrary that there exist  norm-one elements $ b,  b_n \in \mathbb{B}$ ($ n \in \mathbb{N}$) such that $ \lim_n b \circ b_n = 0$.  Now, we have 
 \begin{eqnarray*}
 \| b + b_n \|^4 & = & \Big( \| b\|^2 + \| b_n\|^2 +  2 {\rm Re} \langle b, b_n \rangle\Big)^2 \\
 & =&  \big( 2 +  2 {\rm Re} \langle b, b_n \rangle \big)^2,
  \end{eqnarray*}
  and 
   \begin{eqnarray*}
  \| b + b_n \|^4 & = & \| ( b + b_n) \circ  ( b + b_n) \|^2 = \|  b^2 + b_n^2  + 2  b \circ   b_n \|^2, \\
   & = &  \Big( \|  b^2 \|^2 + \|b_n^2\|^2   + 2 \| b \circ   b_n \|^2 +  2 {\rm Re} \langle b^2, b_n^2 \rangle \\
    & & +  2 {\rm Re} \langle b^2, (b \circ   b_n)^2 \rangle +  2 {\rm Re} \langle b_n^2, (b \circ   b_n)^2 \rangle  \Big) ,  
 \end{eqnarray*}
which, respectively,  imply $\lim_n  \| b + b_n \|^4 =4 $ and  $ \lim_n \| b + b_n \|^4 =2$, a contradiction because 
 \begin{eqnarray*}
| {\rm Re} \langle b, b_n \rangle| & \leq &  \| b \circ b_n  \|, \\
| {\rm Re} \langle b^2, b_n^2 \rangle| & \leq & |  \langle b^2, b_n^2 \rangle| \leq  \| b^2 \circ b_n^2  \|=  \| (b \circ b_n)^2  \| , \\
 & \leq &  \| b \circ b_n  \|^2,   \\
    | {\rm Re} \langle b^2, (b \circ   b_n)^2 \rangle | & \leq &   |\langle b^2, (b \circ   b_n)^2 \rangle | ,\\
    & \leq &    \|  b^2 \circ  (b \circ   b_n)^2 \|  \leq   \|  b^2\|  \|  (b \circ   b_n)^2 \|, \\
     & \leq &  \| b \circ   b_n\|^2 ,\\
   | {\rm Re} \langle b_n^2, (b \circ   b_n)^2 \rangle | & \leq & |  \langle b_n^2, (b \circ   b_n)^2 \rangle | , \\
    & \leq &    \|b_n^2, (b \circ   b_n)^2 \| \leq  \|b_n^2\| \|(b \circ   b_n)^2 \| , \\
     & \leq &  \| b \circ   b_n\|^2 
 \end{eqnarray*}
for all $ n \in \mathbb{N}$. Therefore, $\mathbb{B}$ has non nontrivial topological zero-divisors. It thus follows from part (iii) that $\mathbb{B}$  is isomorphic to $ \mathbb{C}$. That is, we have shown that every singly  generated subalgebra of $ \mathbb{A}^{{\rm sym}}$ is at most one-dimensional. This together with Theorem \ref{1.2}(i) implies  that $ \mathbb{A}^{{\rm sym}}$, and hence $\mathbb{A}$, is one-dimensional, and hence both, which are the same,  are isometrically isomorphic to  $ \mathbb{C}$ because $\mathbb{A}$ is, in particular, unital and $\|1 \| = 1$.

(vii) If necessary, by passing to the completion of $\mathbb{A}$, we may assume that $\mathbb{A}$ is a complex Banach algebra equipped with a unital norm and whose norm comes from an inner product. Just as in (vi), we prove the assertion by showing  that  $ \mathbb{A}^{{\rm sym}}$  is one-dimensional and hence isometrically isomorphic to $ \mathbb{C}$ because $\|1 \| = 1$.   To this end, pick a nonzero $ a \in  \mathbb{A}^{{\rm sym}}$ and form the closed algebra generated by $ a$ in $ \mathbb{A}^{{\rm sym}}$ and denote it by $\mathbb{B}$. We note that  $\mathbb{B}$ is a commutative and associative complex Banach algebra whose norm comes from an inner-product $ \langle . , .  \rangle$ that it inherits from $\mathbb{A}$. By \cite[Corollary 2]{I1}, $\mathbb{B}$ is isomorphic to $ \mathbb{C}$, and hence, in particular, one-dimensional.  Also, since $\mathbb{B}$ is isomorphic to $ \mathbb{C}$, we see that $ a^2 = 0$ implies $ a =0$ for all $ a \in \mathbb{A}$.  Therefore, the algebra generated by any  $ a \in  \mathbb{A}$ is at most one-dimensional and $0$ is the only element of $\mathbb{A}$ whose square is zero. It thus follows from Theorem \ref{1.2}(i), that $\mathbb{A}$ is isomorphic, and hence  isometrically isomorphic, to $ \mathbb{C}$, for $\|1 \| = 1$. This completes the proof. 
\hfill \qed

\bigskip

If the algebras in the preceding theorem happen to be real, the algebras turn out to be locally complex; see \cite[Proposition 2.5.49(ii)]{GP}, \cite[Corollary 2.5.54(ii)]{GP}, \cite[Corollary 2.5.56(ii)]{GP}, \cite[Corollary 2.7.23]{GP},   \cite[Corollary 2.7.24]{GP}, and \cite[Characterization 2]{Za}. Just like the preceding theorem, note that in parts (ii) and (iii) of the theorem the algebra $\mathbb A$ need not be normed; the algebra $\mathbb A$ is only equipped with a vector space norm with respect to which the  multiplication operation of the algebra $\mathbb A$ is separately continuous.

\bigskip

\bigskip 

\begin{thm} \label{1.9} 
{\rm (i)}  Let $\mathbb A$ be a   left (resp. right)  alternative  topological  real algebra  whose (topological) dual separates its points. If $\mathbb A$  is a quasi-division algebra, equivalently, a classical division algebra, then it is locally complex.

{\rm (ii)}   Let $\mathbb A$ be a  one-sided t-division  real algebra equipped with a real vector space norm with respect to which the  multiplication operation of the algebra $\mathbb A$ is separately continuous. Then, $\mathbb A$  is locally complex if it is power-associative.

{\rm (iii)}  Let $\mathbb A$ be a  real algebra equipped with a real vector space norm with respect to which the  multiplication operation of the algebra $\mathbb A$ is separately continuous and with ${\rm j.t.z.d.}(A)= \{0\}$. Then, $\mathbb A$  is locally complex if it  is power-associative. 

{\rm (iv)} 
Let $\mathbb A$ be a  real algebra  endowed with a $c$-supermultiplicative vector space norm for some $ c > 0$ relative to which the  multiplication is separately continuous, e.g., a real nearly absolute value. Then, $\mathbb A$  is locally complex if  it is power-associative.

{\rm (v)} Let $\mathbb A$ be a real  algebra endowed with an inner-product $ \langle . , . \rangle$, whose norm satisfies the identity $\|a^2\| = \|a\|^2 $ on $\mathbb A$. Then, $\mathbb A$  is locally complex if  it is power-associative and unital,  and $\mathbb A$  is a locally complex division algebra if it is left (resp. right) alternative.

{\rm (vi)}  Let $\mathbb A$ be a real  unital algebra endowed with an inner-product $ \langle . , . \rangle$, whose norm is a unital algebra norm on $\mathbb A$. Then, $\mathbb A$ is locally complex  if it is power-associative.

\end{thm}

\bigskip

\noindent {\bf Proof.} (i) Let $\mathbb A_{\mathbb{C}} := \mathbb A \times \mathbb A$ be the complexification of $\mathbb A$ equipped with the product topology. By Proposition \ref{1.1}(vi) and its proof, the algebra $\mathbb A$ is unital and alternative and indeed a classical division algebra. We note that 
 the (topological) dual of $\mathbb A_{\mathbb{C}} $  separates its points. This follows from the hypothesis that the (topological) dual of $\mathbb A$ separates its points and that every $ f \in \mathbb A^*$, gives rise to an $  f_{\mathbb{C}} \in (\mathbb A_{\mathbb{C}})^*$, which is defined by $  f_{\mathbb{C}} ( a + ib) = f(a) + i f(b)$. Then again, by Theorem \ref{1.4}, every $ a \in \mathbb A$  has a nonempty spectrum, implying that $a^2 - 2 {\rm Re}(\lambda) a + |\lambda|^2  1 =0 $ for some $ \lambda \in \mathbb{C}$ because  $\mathbb A$  is a divsion algebra in the classical sense. That is, $\mathbb A$ is quadratic, and hence locally complex for it is a classical division algebra, as desired. 

 (ii) We prove the assertion for left t-division algebras. Likewise, one can prove the assertion for right t-division algebras.  To this end, given an arbitrary $a \in \mathbb{A}$, let   $\mathbb{B}$ be the real subalgebra generated by $a$,  $ \overline{\mathbb{B}}$ be its completion, and $ \overline{\mathbb{B}}_\mathbb{C}$ be a complexification of $ \overline{\mathbb{B}}$.  Just as we argued  in part (iii) of the preceding theorem,  $L_{a}: \overline{\mathbb B}  \longrightarrow \overline{\mathbb B}  $ is continuous. By  \cite[Proposition VII.6.7]{C}, $ \sigma_{ap}\big( L_a ;  \mathcal{B}(\overline{\mathbb{B}}_\mathbb{C} ) \big) \not= \emptyset $ so that    $ \lim_n ( L_a - \lambda I_{\overline{\mathbb{B}}_\mathbb{C}}) ( \bar a_n + i \bar b_n )  = 0$  for some $ \lambda \in \mathbb{C}$ and   $ \bar a_n + i \bar b_n \in \overline{ \mathbb{B}}_\mathbb{C}$ with $\| \bar a_n + i \bar b_n \| = 1$   ($n \in \mathbb{N}$). From this, just as we saw in the proof of Theorem \ref{1.7}(ii),  we get that 
\begin{eqnarray*}
\lim_n ( a^3 - 2 {\rm Re}(\lambda)  a^2 + |\lambda|^2 a )    \bar a_n  & = & 0 , \\
\lim_n ( a^3 - 2 {\rm Re} (\lambda)  a^2 + |\lambda|^2 a )   \bar  b_n  & = & 0 ,
\end{eqnarray*}
which, in turn, since $\| \bar a_n + i \bar b_n \| = 1$   ($n \in \mathbb{N}$), implies that 
\begin{eqnarray*}
\lim_n ( a^3 - 2 {\rm Re}(\lambda)  a^2 + |\lambda|^2 a )     c_n  & = & 0 , 
\end{eqnarray*}
for some $ c_n \in \mathbb{B}$ with $\| c_n \| = 1$  ($n \in \mathbb{N}$) because $\mathbb{B}$ is dense in $ \overline{\mathbb{B}}$. Consequently, $ a^3 - 2 {\rm Re}(\lambda)  a^2 + |\lambda|^2 a   =  0$.  Therefore,   $\mathbb A$  is a power-associative algebraic algebra, and hence, from \cite[Lemma 2.5.5]{GP}, we see that $\mathbb A$ is unital, for  $\mathbb A$ has a nonzero idempotent by the proof of  Proposition \ref{1.1}(i). Then again, $\mathbb A$ is a t-division real algebra. Thus, $ a^2 - 2 {\rm Re}(\lambda)  a + |\lambda|^2 1   =  0$,  as desired.

(iii) The proof, which is omitted for the sake of brevity, is similar to that of part (ii) except that it must be adjusted while making use of the hypothesis of (iii).

(iv) This is a consequence of (iii).

(v)  First, suppose that the algebra is unital and power-associative. Pick a nonzero $ a \in  \mathbb{A}$ and form the algebra generated by $ a$ in $ \mathbb{A}$ and denote it by $\mathbb{B}$. The real algebra  $\mathbb{B}$ is commutative and associative. Clearly, the algebra  $\mathbb{B}$ inherits an inner product from  $\mathbb{A}$ whose norm satisfies  the identity $\|b^2\| = \|b\|^2 $ on $\mathbb B$. It follows from the proof of part (vii) of the preceding theorem rewritten  in the setting of real algebras or from the proof of \cite[Lemma 2.6.34]{GP} that the algebra   $\mathbb{B}$  has no nontrivial joint topological zero-divisors and that the norm induced by its inner-product is indeed an algebra norm on $\mathbb{B}$, and hence relative to which the  multiplication operation of the algebra $\mathbb A$ is  continuous. This together with part  (iii) of the theorem implies that $\mathbb{B}$ is locally complex. The proof is complete at this point if the algebra is unital and power-associative. 

Next, suppose that the algebra is left (resp. right) alternative. Since the algebra    $\mathbb{B}$  is locally complex, we see that $ 0$ is the only nilpotent element of the algebra $ \mathbb{A}$ and that $ \mathbb{A}$ is an algebraic algebra. It thus follows from  \cite[Corollary 1 (Kleinfeld) on page 345]{ZSSS} or  the comment following \cite[Theorem, page 944]{K}  that $\mathbb A$ is alternative.
   Once again, given $ a, b \in \mathbb{A}$,  we have  $ \langle a, b \rangle = 0$ whenever $ a b = ba = 0$. To see this, we can write 
\begin{eqnarray*}
\|a\|^2 + \|b\|^2 \pm 2  \langle a, b \rangle & = & \| a \pm b \|^2 ,\\
& = & \| a^2  + b^2 \pm 2 a  b  \|, \\
& \leq & \| a\|^2  +\| b \|^2, 
\end{eqnarray*}
from which we obtain  $ \langle a, b \rangle = 0 $. From this, we conclude that  $\mathbb A$  has no nontrivial joint zero-divisors. To see this, suppose on the contrary that there are norm-one  elements $ b_1, b_2 \in \mathbb{A}$ with $ b_1 b_2 = b_2 b_1 = 0$. Since  $ \langle b_1, b_2 \rangle = 0$, by what we just saw, we get that
$\|b_1^2\| + \|b_2^2\| =  \| b_1^2  + b_2^2 \|$, from which we obtain $ b_1^2 = b_2^2$. But   $ b_1 b_2 = b_2 b_1 = 0$ and  $\mathbb A$ is  alternative. Thus,  
$ b_1^2 b_2^2 = (b_1 b_2)^2 = 0$, from which we obtain $ b_1^4 = 0$, implying that $ b_1 = 0$, a contradiction.  Therefore, $ \mathbb{A}$ is an  algebraic algebra with ${\rm j.z.d.}(\mathbb{A}) = \{ 0\}$. It thus follows from Theorem \ref{1.1}(iv) that $ \mathbb{A}$ is a classical division algebra. 
This completes the proof. 

(vi) Pick a nonzero $ a \in  \mathbb{A}$ and form the unital algebra generated by $ a$ in $ \mathbb{A}$ and denote it by $\mathbb{B}$. The real algebra  $\mathbb{B}$ is unital, commutative, and associative. Also, the algebra  $\mathbb{B}$  inherits an inner product from  $\mathbb{A}$ whose norm is a unital algebra norm on $\mathbb{B}$. It thus follows from \cite[Theorem 2]{I2} that $\mathbb{B}$ is isomorphic  to  $ \mathbb{R}$ or $ \mathbb{C}$, completing the proof.  
\hfill \qed

\bigskip

Part (i) of the following theorem is motivated by \cite[Proposition 2.7.33]{GP}, which is due to A. Rodr\'iguez; see \cite[Theorem 2]{Ro}.

\bigskip 

\begin{thm} \label{1.10} 
{\rm (i)}  Let $ ( \mathbb{A}, | .|)$ be a unital absolute-valued real algebra and $\langle . , . \rangle$ be its inner product. For $ a \in   \mathbb{A}$, set $a^* := 2 \langle a,  1 \rangle   1 - a$. 
Then, $ a^2 - 2 \langle a,  1 \rangle   a + |a|^2 1 = 0 $, and hence  $ a a^* = a^* a = |a|^2 1$. Moreover,  
$ \langle a  b, c \rangle  =\langle  b, a^* c \rangle  $, $ \langle a  b, c \rangle  =\langle  a, c b^* \rangle  $, $ a^*(ab) = |a|^2b$, $ (ab)b^* = |b|^2 a$,  and hence
$ a(ab) = a^2b$ and $ (ab)b = ab^2$
 for all $ a, b, c \in \mathbb{A}$. In particular,  the algebra $\mathbb A$ is  locally complex and alternative.

{\rm (ii)}  Let $ ( \mathbb{A}, | .|)$ be a power-associative  absolute-valued algebra. Then, $\mathbb A$ is unital, and hence its absolute value comes from an inner product. In particular,  $\mathbb A$ is  alternative.  

{\rm (iii)}  Let $ ( \mathbb{A}, | .|)$ be a flexible absolute-valued real algebra with a left (resp. right) identity element $ e \in  \mathbb{A} $. Then, $\mathbb A$ has an identity element, and hence all the assertions of part (i) hold. 

\end{thm}

\bigskip

\noindent {\bf Proof.} (i)  First, clearly, for the first assertion, it suffices to prove it subject to the additional hypothesis that $ |a| =1$. So suppose $|a| = 1$. It is easily checked that  
  \begin{eqnarray*}
| a( a-  2 \langle a,  1 \rangle 1)|^2 &  =&  1, \\
2  \langle a^2,  1 \rangle  & = & 2 - |a^2 - 1|^2, \\
 4  \langle a,  1 \rangle^2 & =&  ( 2 - |a-1|^2)^2 .
 \end{eqnarray*}
Thus, we can write 
  \begin{eqnarray*}
| a^2 -  2 \langle a,  1 \rangle   a + 1|^2 & =& \Big\langle a( a-  2 \langle a,  1 \rangle 1)    + 1,  a(a -  2 \langle a,  1 \rangle 1)   + 1 \Big\rangle, \\
& =& | a( a-  2 \langle a,  1 \rangle 1)|^2 + 2 \big\langle a( a-  2 \langle a,  1 \rangle ), 1  \big\rangle + 1, \\
& =& 1 + 2  \langle a^2,  1 \rangle - 4  \langle a,  1 \rangle^2 + 1, \\
& =& 2 + (2 - |a^2 - 1|^2) - ( 4 - 4 |a-1|^2 + | a-1|^4) ,\\
& =& - |a - 1|^2 \big( |a + 1|^2 + |a-1|^2 - 4 \big) , \\
& =& - |a - 1|^2 \big( 2 (|a|^2 + 1) - 4 \big) , \\
& =& 0,
 \end{eqnarray*}
proving the first assertion.  The rest is a quick consequence of  Theorem \ref{02.5}(i) together with the second remark following the theorem.

(ii) If the algebra   $\mathbb{A}$ is complex, by Theorem \ref{1.8}(v),   $\mathbb{A}$  is isometrically isomorphic to $\mathbb{C}$, in which case we have nothing to prove. If the algebra $\mathbb{A}$ is real,  by  Theorem \ref{1.9}(iv),  $\mathbb{A}$ is locally complex, and hence unital. Then again, this along with  Theorem \ref{02.3}(ii) (or  Theorem \ref{1.7}(ii-a)) and (i) completes the proof.

(iii)  By Theorem \ref{1.7}(ii-a), the absolute value of $\mathbb A$, namely, $ |.|$,  comes from an inner product
  $\langle . , . \rangle$ on $\mathbb A$. Since the inner-product $\langle . , . \rangle$  permits composition, we see from Theorem \ref{02.5}(i) that $ e $ is indeed an identity element of $\mathbb A$. That is, $ e = 1$. Thus, we are done by Theorem \ref{02.5}(i) or part (i) applies, completing the proof.  
\hfill \qed

\bigskip

We conclude this section with two useful theorems. 

\bigskip 

\begin{thm} \label{1.11} 
 {\rm (i)}   Let $\mathbb{A}$ be a locally complex  algebra and $ k >1$ be a  given natural number. Then, there exists at most one  nonzero vector space norm $\| .\|$ on $\mathbb{A}$  satisfying the identity $ \|x^ k \| = \|x\|^k$ on  $  \mathbb{A}$. Therefore,  the only nonzero real vector space norm $\|.\|$ on  $\mathbb{R}_n$, the 
real Cayley-Dickson algebra of dimension $2^n$ with $n$ being a nonnegative integer, satisfying the identity $ \|x^ k \| = \|x\|^k$ on $\mathbb{R}_n$ is the Euclidean norm of it, denoted by $|.|$, which is given by $ | x| = \sqrt{\sum_{i=1}^{2^n} x_i^2}$, where $ x= (x_i)_{i=1}^{2^n} \in \mathbb{R}_n$.

{\rm (ii)}  There exists at most one absolute value on any locally complex  algebra. In particular, the Euclidean norm is the only absolute value on the numerical division rings, namely, on $\mathbb{R}$, $\mathbb{C}$,  $\mathbb{H}$, and $\mathbb{O}$.

{\rm (iii)} Let $\mathbb{A}$ be a locally complex  algebra. Then, up to multiplication by a positive scalar, there exists at most one vector space norm $\| .\|$ on $\mathbb{A}$ satisfying the relation 
$$ \inf \{ \| xy\| :  \| x\| = \| y\| = 1 \}   = \sup\{  \|xy\| :    \| x\| = \| y\| = 1 \} .$$

\end{thm}

\bigskip

\noindent {\bf Proof.} (i) The assertion trivially holds if the algebra $\mathbb{A}$ consists of scalars. Let $\mathbb{A}$ be a locally complex quadratic real algebra with a nonscalar element and $\|.\|_1$ and $\|.\|_2$ be vector space norms on $\mathbb{A}$  satisfying the desired identities
$$ \|x^ k \|_1 = \|x\|_1^k ,   \  \|x^ k \|_2 = \|x\|_2^k . $$
on $ \mathbb{A}$. 
 It suffices to show that $\|a\|_1 = \|a\|_2$ for all 
nonzero elements $ a \in  \mathbb{A}$. It is plain that $\| 1 \|_1 =\| 1 \|_2 = 1$.  Clearly,  given an arbitrary $ a \in  \mathbb{A}$, there exists a $ j \in \mathbb{A}$ with $ j^2 = -1$ such that  $ a \in \mathbb{C}_j$, where $  \mathbb{C}_j := \{ r1 + s j : r  , s \in \mathbb{R}\} $. Let $|.|$ denote the absolute value on the copy of complex numbers generated by $ j$ in $ \mathbb{A}$ and $T$ denote the unit circle in $  \mathbb{C}_j$   relative to $|.|$. Since $  \mathbb{C}_j$ is two-dimensional, the three norms, namely, $\|.\|_1$, $\|.\|_2$, and  $|.|$, on $  \mathbb{C}_j$ are equivalent. But the set  $U_k$ consisting of all the $k^m$-th roots of unity in $  \mathbb{C}_j$ with $m$ ranging over $\mathbb{N}$ is dense in $T$ with respect to the topology induced by $|.|$, and hence by those induced by $\|.\|_1$ and $\|.\|_2$, on $T$. Clearly,  $ \|u\|_1 = 1 = \|u\|_2 $ for all $ u \in U_k$ and $ b := \frac{a}{|a|} \in T$.  Thus, there exists a sequence $ (u_n)_{n=1}^\infty $ in $U_k$  such that $ b = \lim_n u_n $ with respect to  $|.|$, and hence relative to $\|.\|_1$ and $\|.\|_2$. Then again,  we can write 
$$ \|b\|_1 =  \| \lim_n u_n\|_1 = \lim_n  \| u_n\|_1 = 1 , $$
implying that $ \|a\|_1 = |a|$. Likewise, $ \|a\|_2 = |a|$, from which we conclude that $ \|a\|_1 = \|a\|_2$, as desired. 
The second assertion follows from the fact that the Euclidean norm of $\mathbb{R}_n$, denoted by $|.|$, indeed satisfies the identity $ |x^ k | = |x|^k$ on $\mathbb{R}_n$  for all $ k \in \mathbb{N}$.  
Just as we saw in the proof of Proposition  \ref{02.4},  this identity follows from the following well-known facts: that $\mathbb{R}_n$ is quadratic, and hence power-associative, and that 
$$   |x|^2 1= x^* x = x x^*, \ 2 {\rm Re}(x) = x + x^* , \ x^2 = 2 {\rm Re}(x) x - |x|^21 ,$$
for all $ x \in \mathbb{R}_n$.

(ii) This is a quick consequence of (i).

(iii) Let $\|.\|_1 $ and $\|.\|_2$ be vector space norms on $\mathbb{A}$ with 
$$ c_1 := \inf \{ \| xy\|_1 :  \| x\|_1 = \| y\|_1 = 1 \}   = \sup\{  \|xy\|_1 :    \| x\|_1 = \| y\|_1 = 1 \} ,$$
$$ c_2 := \inf \{ \| xy\|_2 :  \| x\|_2 = \| y\|_2 = 1 \}   = \sup\{  \|xy\|_2 :    \| x\|_2 = \| y\|_2 = 1 \} .$$
It is quite straightforward to check that $c_1\|.\|_1 $ and $c_2\|.\|_2$ are absolute values on $\mathbb{A}$. It thus follows from (ii) that $c_1\|.\|_1 =c_2\|.\|_2$,  which is what we want. 
This completes the proof. 
\hfill \qed

\bigskip 

It turns out every locally complex algebra is equipped with a unique  inner-product whose  norm commutes with the $k$th-power function for all $ k \in \mathbb{N}$. We call this inner product and its norm {\it the natural inner product} and {\it the natural norm of a locally complex algebra}. 

\bigskip

\begin{thm} \label{2.13} 
 Let $\mathbb{A}$ be a locally complex algebra. Then, $\mathbb{A}$ is  equipped with a unique vector space norm $\|.\|$ that comes from an inner-product satisfying the identity $ \|x^k \| = \|x\|^k$ on  $  \mathbb{A}$ for all $ k \in \mathbb{N}$.

\end{thm}

\bigskip

\noindent {\bf Proof.}  Uniqueness follows from the preceding theorem, indeed, under the weaker hypothesis that the vector space norm $\|.\|$ satisfies the identity $ \|x^k \| = \|x\|^k$ on  $  \mathbb{A}$ for some $ k \in \mathbb{N}$ with $ k > 1$. So we only need to prove the existence. Once again, the assertion trivially holds if the algebra $\mathbb{A}$ consists of scalars. We prove the assertion for the case when $\mathbb{A}$  contains a nonscalar element. By Theorem \ref{02.3}(ii), there exists a symmetric bilinear form 
$ \langle . , . \rangle: \mathbb{A} \times \mathbb{A} \longrightarrow \mathbb{R}$ such that 
$$ a^2 - 2  \langle a , 1 \rangle a + \langle a , a \rangle 1 = 0, $$
and that $ \langle t1 , 1 \rangle = t $ and $ \langle t1 , t1 \rangle = t^2$ for all $t \in \mathbb{R}$. We show that $  \langle . , . \rangle$ is a real inner-product on $  \mathbb{A}$ by showing that it is positive definite.  Define the functions $ p_1, n: \mathbb{A} \longrightarrow \mathbb{R}$ by
$$ p_1 (a) := \langle a , 1 \rangle , \  n(a) :=  \langle a , a \rangle ,    \  ( a \in \mathbb{A}).$$
Indeed,  $p_1(a)$ and $n(a)$ are the unique real numbers satisfying the relation $ a^2 - 2 p_1(a) a + n(a) 1 = 0$. 
Now, note that $ n(a) \geq 0$ for all $ a \in  \mathbb{A}$. To see this, by way of contradiction, pick an $ a \in  \mathbb{A}$ with $ n(a) < 0$.  We then have 
 $ \big( a - p_1(a)1 \big)^2 = \big(p_1(a)^2 - n(a)\big)1$, from which we obtain $ a = \big(p_1(a) \pm \sqrt{p_1(a)^2 - n(a)}\big)  1 \in \mathbb{R}1 $, and hence $n(a) \geq 0$, which is impossible. Thus, $ n \geq 0$ on  $\mathbb{A}$. 
 That $ \langle a, a \rangle =0$ implies $ a =0$ follows from the last assertion of Theorem  \ref{02.3}(ii). Therefore, $ \langle . , . \rangle$ is a real inner-product on $\mathbb{A}$. Finally, $ \|a^k \| = \|a\|^k$, equivalently, $ n(a^k) = n(a)^k$, for  all $ a \in \mathbb{A}$ and $ k \in \mathbb{N}$ also follows from Proposition  \ref{02.4}. This completes the proof. 
\hfill \qed

\bigskip

\bigskip

\begin{section}
{\bf Revisiting Theorems of Frobenius, Hurwitz, and Zorn and their topological counterparts}
\end{section}

\bigskip

We start off this section with a useful lemma.

\bigskip

\begin{lem} \label{2.1} 
  {\rm (i)} Let $F$ be a field whose characteristic is not $2$,   $\mathbb{A}$  a left (resp. right) alternative  quadratic algebra over $F$ that is locally a field extension of $F$, $\langle . , . \rangle$ the symmetric bilinear form of  $\mathbb{A}$, and $ p, q \in \mathbb{A} \setminus F1$.  Then, either $pq = qp$, in which case 
$\langle \{ 1, p \}\rangle_F = \langle \{1, q \}\rangle_F$, or else $ pq + qp \in \langle \{1, p, q \} \rangle_F$ and $ \{ 1, p , q , pq \}$ is linearly independent over $F$. Moreover, if $ p^2 ,  q^2 \in F 1$, then either $pq = qp$, in which case $ p = \pm q$ and $ pq + qp= \pm 2 n_1 1   \in F1$ where $ n_1 = \langle p , p \rangle =  \langle 	q , q \rangle  $, or else $ pq + qp= 2 f 1   \in F1$, where $ f =  \langle pq , 1 \rangle = \langle qp , 1 \rangle =  -\langle p , q \rangle$.

{\rm (ii)}   Let $\mathbb A$  be an algebraic left (resp. right) alternative real algebra with no nontrivial joint zero-divisors. Then,   $ \mathbb{A}$ is a locally complex algebra.  If $ p, q \in \mathbb{A} \setminus \mathbb{R}1$, then, either $pq = qp$, in which case 
$\langle \{ 1, p \} \rangle_{\mathbb{R}} = \langle \{ 1, q \} \rangle_{\mathbb{R}}$, or else $ pq + qp \in \langle \{ 1, p, q \} \rangle_{\mathbb{R}}$ and $ \{ 1, p , q , pq \}$ is linearly independent over $\mathbb{R}$. Moreover, if $ p^2 = q^2 \in \mathbb{R} 1$, then $pq = qp$, in which case $ p = \pm q$, or else $ pq + qp= 2 f 1   \in \mathbb{R}1$, where $f=  {\rm Re}(pq) = {\rm Re}(qp)=   -\langle p , q \rangle$.

\end{lem}

\bigskip

\noindent {\bf Proof.}  (i)  By \cite[Theorem 1]{A2}, the algebra  $\mathbb A$ is alternative.  By  Theorem \ref{02.5}(ii), the algebra $\mathbb A$ has no nontrivial zero-divisors.  First, suppose $ p, q \in \mathbb{A} \setminus F1$ and  $ pq = qp$. Let 
 $$ p_1 := \langle p , 1 \rangle  ,   q_1 :=  \langle q , 1 \rangle ,  $$
 $$n_1 :=  \langle p - p_1 1 , p - p_1 1 \rangle   ,    n_2 :=\langle q - q_1 1, q - q_1 1\rangle $$
 $$ n := \big\langle  p - q - (p_1 - q_1) 1,  p - q - (p_1 - q_1) 1 \big\rangle .$$ 
 It follows that 
 $$(p - p_1 1)^2 + n_1 1 = 0 = (q - q_1 1)^2 + n_2 1 , $$
 which, in view of  $ pq = qp$, yields 
 $$ \big(  p - q - (p_1 - q_1) 1 \big) \big(  p + q - (p_1 + q_1) 1 \big) = ( n_2 - n_1) 1.$$
 If  $n_2 - n_1= 0$, then  $ p - q - (p_1 - q_1)1 = 0$ or  $  p + q - (p_1 + q_1) 1 = 0$, implying that $\langle \{ 1, p \}\rangle_F = \langle \{1, q \}\rangle_F$, as desired.  If $n_2 - n_1\not= 0$, then $r:=  p - q - (p_1 - q_1)1 \not= 0$, and hence $ n \not= 0$. But 
 $$ r^{-1} = \big(  p - q - (p_1 - q_1) 1 \big)^{-1} = \frac{1}{n} \big( (p_1 - q_1) 1  - (p - q)  \big)$$
  and   $ r \big(  p + q - (p_1 + q_1) 1 \big) = ( n_2 - n_1)r r^{-1} $ and $\mathbb A$ has no nontrivial zero-divisors. 
So we get that 
   $$    p + q - (p_1 + q_1) 1  = \frac{n_2 - n_1}{n} \big( (p_1 - q_1) 1  - (p - q)  \big)  .$$
 Note that  $  \frac{n_2 - n_1}{n} \not= \pm 1$, for otherwise $ p \in F 1 $ or $ q \in F 1$, which is impossible. Thus, $  \frac{n_2 - n_1}{n} \not= \pm 1$, which along with the above equality yields $\langle \{ 1, p \}\rangle_F = \langle \{1, q \}\rangle_F$, as desired. Next, suppose  $ pq \not= qp$. By  Proposition \ref{02.4}, we can write 
 $$ pq + qp =  2\langle q , 1 \rangle  p + 2  \langle p , 1 \rangle  q - 2  \langle p , q \rangle 1 \in  \langle \{1, p, q \} \rangle_F, $$
   as desired. 
  Next, to see that $ \{ 1, p , q , pq\}$ is linearly independent over $F$, suppose $ f_01+ f_1 p + f_2q + f_3 pq =0$ for some $ f_0, f_1 , f_2 , f_3 \in F$. We must show that $  f_i= 0$ ($0 \leq i \leq 3$).  To this end, it suffices to show that $ f_2=f_3 = 0$ because  $ p \notin F1$. Suppose by way of  contradiction that $ ( f_2, f_3)  \not = 0$ so that $   f_01+ f_1 p + ( f_21 + f_3 p) q   =0 $. It thus follows from the alternativity of $\mathbb A$ that $ q = -(f1 +  f'p)( f_01+ f_1 p) \in  \langle \{ 1, p \} \rangle_{F}$, where $f1 +  f'p = ( f_21 + f_3 p)^{-1}$,  a contradiction because $ pq \not= qp$. 
This shows that $  f_i= 0$ ($0 \leq i \leq 3$), completing the proof of the first assertion.

 For the second assertion,  suppose $ p^2 ,  q^2 \in F 1$.  Since $ p, q \in \mathbb{A} \setminus F1$, we see that  $ p_1 = q_1 = 0$ so that $ p^2 = -n_11 $,  $ q^2 = - n_2 1 $, and $ (p - q)^2 = - n 1$.   First,  if  $pq = qp$, then $ ( p- q ) ( p + q) = p^2 - q^2= (n_2 - n_1)1 $. It follows that $ n_1 = n_2$. To see this, suppose on the contrary that  $ n_2 - n_1 \not= 0$. This implies $ p - q  \not= 0$, from which, we obtain $ p + q =\frac{n_2 - n_1}{n} ( q - p)  $. Then again, $ \frac{n_2 - n_1}{n}  \not= \pm 1$ because $ p, q \in \mathbb{A} \setminus F1$, from which we get that $ pq = qp$, a contradiction. Thus $ n_1 = n_2$, and hence $ ( p- q ) ( p + q) =0$, entailing that  $ p = \pm q$. This, in turn,  yields $ pq + qp= \pm 2 n_1 1   \in F1$, as desired. Second,  if $pq \not= qp$, since $ p_1 = q_1 = 0$, from what  we saw in the above $ pq + qp= -2 \langle p , q \rangle  1   \in F1$. But by   Theorem \ref{02.7}, $f :=  \langle pq , 1 \rangle = \langle qp , 1 \rangle =  -\langle p , q \rangle$ and $ pq + qp= 2f 1   \in F1$. This completes the proof.

 (ii)  By \cite[Lemma 1]{A2}, the algebra  $\mathbb A$ is power-associative. It thus follows from  Proposition \ref{1.1}(iii) that  the real algebra  $\mathbb A$ is unital and hence quadratic by the fundamental theorem of algebra for polynomials with real coefficients.  Consequently, $\mathbb A$ is alternative by \cite[Theorem 1]{A2}.  The real algebra $\mathbb A$ is locally complex because  it  has no nontrivial joint zero-divisors. Note that the symmetric bilinear form of  $\mathbb{A}$, denoted by  $\langle . , . \rangle$, is indeed a real inner-product on   $\mathbb{A}$ and $ {\rm Re}(a) = \langle a , 1 \rangle $ for all $ a \in \mathbb A$. Thus, part (i) applies, completing the proof. 
  \hfill \qed 

\bigskip

The following theorem can be thought of as a Gelfand-Mazur type theorem for commutative quadratic algebras over fields whose characteristics are not $2$.  

\bigskip

\begin{thm} \label{5.2} 
 Let $F$ be a field whose characteristic is not $2$ and   $\mathbb{A}$  a  quadratic algebra over $F$. Then,  $\mathbb{A}$ is a field extension of $F$ of degree at most two if and only if  $ \mathbb{A}$ is commutative and has  no nontrivial zero-divisors. Therefore, the quadratic algebra $\mathbb{A}$ is a field extension of $F$ of degree at most two if and only if $ ab + ba = 0$ implies $a = 0$  or $b = 0$ for all $a, b \in  \mathbb{A}$.
\end{thm}

\bigskip

\noindent {\bf Proof.}  First, suppose that $ \mathbb{A}$ is commutative and has no nontrivial zero-divisors.  Let  $\langle . , . \rangle$ denote the symmetric bilinear form of  $\mathbb{A}$.  Assuming that there is an element  $ a \in   \mathbb{A} \setminus F1  $,  let $ b \in \mathbb{A} \setminus F1 $ be arbitrarily given. Since   $ b a = ab$ and $\mathbb{A}$  has no nontrivial zero-divisors, we see from the proof of part (i) of the preceding lemma that  $ \langle \{ 1, b \} \rangle_F  =  \langle \{ 1, a \} \rangle_F$. This, in particular,  yields $ \mathbb{A} = \langle \{ 1, a \} \rangle_F = \{ r1 + s a : r, s \in F \}$, as desired.

 Next, the second assertion follows from the first assertion because  $ \mathbb{A}^{{\rm sym}}$  is a commutative quadratic algebra over $F$ with no nontrivial zero-divisors.
\hfill \qed 

\bigskip

The following can be thought of as the counterparts of the Gelfand-Mazur Theorem for commutative real algebras; see \cite[Proposition 2.5.52]{GP},  \cite[Corollary 2.5.58]{GP}. We would like to stress that in subparts (a) and (b) of part (ii) of the theorem the algebra $ \mathbb{A}$ is endowed with a vector space norm as opposed to an algebra norm.

\bigskip

\begin{thm} \label{2.2} 
 {\rm (i)}  Let $ \mathbb{A}$  be an algebraic power-associative real algebra. If $ \mathbb{A}$ is commutative and has  no nontrivial  zero-divisors, then the algebra  $ \mathbb{A}$ is locally complex and, equipped with its natural norm, is isometrically isomorphic to $\mathbb{R}$ or $\mathbb{C}$. Therefore, if $ ab + ba = 0$ implies $a = 0$  or $b = 0$ for all $a, b \in  \mathbb{A}$, then  the algebraic power-associative real algebra $ \mathbb{A}$  is locally complex and, equipped with its natural norm, is isometrically isomorphic to $\mathbb{R}$ or $\mathbb{C}$.

\bigskip

{\rm (ii)}  Let $ \mathbb{A}$  be a power-associative real algebra. If   $ \mathbb{A}$ is commutative and  one of the following conditions holds, then the algebra $ \mathbb{A}$ is locally complex and, equipped with its natural norm, is isometrically isomorphic to $\mathbb{R}$ or $\mathbb{C}$.

{\rm (a)} $\mathbb A$ is a  one-sided t-division real algebra equipped with a vector space norm relative to which the multiplication of $\mathbb A$ is separately continuous.

 {\rm (b)} $\mathbb A$ is a real algebra  equipped with a vector space norm relative to which the multiplication of $\mathbb A$ is separately continuous and with  ${\rm j.t.z.d.}(\mathbb{A}) = \{ 0\}$.

 {\rm (c)} $\mathbb A$ is a  real algebra  equipped with a $c$-supermultiplicative vector space norm for some $ c > 0$ relative to which the  multiplication is separately continuous, e.g., a real nearly absolute value. 

\bigskip

{\rm (iii)}  Let $\mathbb A$ be a   left (resp. right) alternative  topological  real algebra  whose (topological) dual separates its points. If $\mathbb A$  is a quasi-division algebra, equivalently, a classical division algebra, and commutative, then it is locally complex and, equipped with its natural norm, is isometrically isomorphic to $\mathbb{R}$ or $\mathbb{C}$.

\bigskip

{\rm (iv)}  Let $\mathbb A$ be a power-associative  real  algebra endowed with an inner-product $ \langle . , . \rangle$, whose norm satisfies the identity $\|a^2\| = \|a\|^2 $ on $\mathbb A$. Then, $\mathbb A$ is isometrically isomorphic to $\mathbb{R}$ or $\mathbb{C}$  if  it is commutative and left (resp. right) alternative, in which case the given norm of $ \mathbb{A}$ coincides with its natural norm that comes  from its local complex structure.

\bigskip

{\rm (v)}   Let $ \mathbb{A}$  be a power-associative absolute-valued real algebra. Then, $ \mathbb{A}$ is isometrically isomorphic to $\mathbb{R}$ or $\mathbb{C}$  if and only if  $ \mathbb{A}$ is commutative, in which case the given absolute value of $ \mathbb{A}$ coincides with its natural norm  that comes from its local complex structure.

\end{thm}

\bigskip

\noindent {\bf Proof.} (i) Since  $ \mathbb{A}$ is algebraic,  Proposition \ref{1.1}(iii) or the proof of Proposition \ref{1.1}(i) along with \cite[Lemma 2.5.5]{GP}  reveals that the algebra $ \mathbb{A}$ is unital and hence  locally complex, for it has no nontrivial joint zero-divisors.  The assertion thus follows from the preceding theorem and Theorem \ref{2.13}, for $ \mathbb{A}$ is commutative and locally complex. Likewise, the second assertion follows from that of the preceding theorem and Theorem \ref{2.13} because  $ \mathbb{A}^{{\rm sym}}$  is a commutative quadratic algebra over reals with no nontrivial zero-divisors.

(ii) and (iii) follow from (i) and Theorems \ref{1.9} and \ref{2.13}.

(iv)  follows from (i) and Theorems \ref{1.9}(v) and \ref{2.13}.

(v) follows from part (ii)-(c) and Theorem \ref{2.13}. 
\hfill \qed

\bigskip

The next two lemmas are crucial for our proofs of the theorems of Frobenius, Hurwitz, and Zorn.

\bigskip 

\begin{lem} \label{2.3}  {\rm (i)} Let $F$ be a field whose characteristic is not $2$,   $\mathbb{A}$  a left (resp. right) alternative  quadratic algebra over $F$ that is locally a field extension of $F$,  and $ p, q_0 \in \mathbb{A} \setminus F1  $ with $ p^2 ,  q_0^2 \in  F1$ and  $ q_0 \notin \langle \{ 1, p \} \rangle_{F} $, equivalently, $ pq_0 \not= q_0 p$.  Then,  there exists an element $q = r p + s q_0 \in  \mathbb{A}$ with $ r, s \in F$ such that $ q^2 \in F 1$, $ pq = -q p$, and $(pq)^2 = - p^2 q^2 \in F 1$. Therefore, the algebra $ \mathbb{A}$ contains a copy of the quaternion algebra generated by $\{ 1, p, q, pq\}$ over the field $F$.

{\rm (ii)}   Let $ \mathbb{A}$  be an algebraic left (resp. right) alternative  real algebra with no nontrivial joint zero-divisors and $ p, q_0 \in  \mathbb{A}$ with $ p^2 = q_0^2 = -1$ with $ q_0 \notin \langle \{ 1, p \} \rangle_{\mathbb{R}} $, equivalently, $ pq_0 \not= q_0 p$. Then, there exists an element $q = r p + s q_0 \in  \mathbb{A}$ with $ r, s \in \mathbb{R}$ such that $ q^2 = -1$, $ pq = -q p$, and $(pq)^2 = - 1$. Therefore, the algebra $ \mathbb{A}$ contains a copy of the real quaternions  generated by $\{ 1, p, q, pq\}$. 
\end{lem}

\bigskip

\noindent {\bf Proof.}  (i) By \cite[Theorem 1]{A2}, the algebra  $\mathbb A$ is alternative. Let $\langle . , . \rangle$ be the symmetric bilinear form of  $\mathbb{A}$. Let $ q:= q_0 - \frac{\langle q_0 , p \rangle}{\langle p , p \rangle} p$. We have $ q \in \mathbb{A} \setminus F1  $, $ \langle q ,  1 \rangle = 0 =  \langle q , p \rangle$. This, along with Theorem \ref{02.7}, implies $ q^2 =  - \langle q ,  q \rangle 1 \in F 1$ and $ pq = -q p$.  Now,  since $\mathbb A$ is alternative, we have $(pq)^2 = - (qp) (pq) = - p^2 q^2 \in F 1$, completing the proof. 

(ii) As noted in the preceding lemma,  $\mathbb A$ is alternative and locally complex. Let  $\langle . , . \rangle$ denote the inner-product of $\mathbb A$ and $ |.| $ the absolute-value induced by the inner-product of $\mathbb A$. Recall that, by  Theorem \ref{02.5}(ii),  $ |.| $ is indeed an absolute value because $\mathbb A$ is alternative and locally complex. With this in mind, the proof is almost identical to that of part (i) except that $ q \in \mathbb{A}$ should be defined by $ q := \frac{1}{|\hat q|} \hat q$, where $ \hat q := q_0 - \langle q_0 , p \rangle p$; note that $ \langle p , p \rangle = 1$, for $ p^2 =-1$ and $ q^2 = -1$ because $ \langle q , 1 \rangle = 0$ and $ \langle q , q \rangle =  |q|^2 =1$.  This completes the proof. 
\hfill \qed 

\bigskip

\bigskip 

\begin{lem} \label{2.4} 
 {\rm (i)} Let $F$ be a field whose characteristic is not $2$,   $\mathbb{A}$  a left (resp. right) alternative  quadratic algebra over $F$  that is locally a field extension of $F$,  and $ p, q \in \mathbb{A} \setminus F1 $ with $ p^2 ,  q^2 \in  F1$ and  $ pq = -q p$, equivalently,  $ \mathbb{A}$ contains a copy of the quaternion algebra generated by $\{ 1, p, q, pq \}$.   If $    \langle  \{ 1, p, q, pq \} \rangle_F \not=  \mathbb{A} $, then there exists an $r \in  \mathbb{A}$ with $ r^2 \in  F1 $, $ p r = - r p$,  $ q r = - r q$, and $ (pq)r  = - r (pq)$ such that $\{ 1, p, q, pq, r, pr, qr,  (pq)r\}$ is linearly independent over $ F$.    Moreover, 
$$  \mathbb{A} =  \langle  \{1, p, q, pq, r, pr, qr,  (pq)r \} \rangle_F,$$
which is  the  octonion algebra generated by $\{ 1, p, q, pq, r, pr, qr,  (pq)r\}$ over the field $F$.

{\rm (ii)}  
Let $ \mathbb{A}$  be an algebraic left (resp. right) alternative real algebra with no nontrivial joint zero-divisors  and $ p, q \in  \mathbb{A} \setminus \mathbb{R}1 $ with $ p^2 = q^2 = -1$ and $ p q= - q p$, which means  $ \mathbb{A}$ contains a copy of the real quaternions generated by $\{ 1, p, q, pq \}$. If $    \langle  \{ 1, p, q, pq \} \rangle_\mathbb{R} \not=  \mathbb{A} $, then there exists an $r \in  \mathbb{A}$ with $ r^2 = -1$, $ p r = - r p$,  $ q r = - r q$, and $ (pq)r  = - r (pq)$ such that $\{ 1, p, q, pq, r, pr, qr,  (pq)r\}$ is linearly independent over $ \mathbb{R}$.    Moreover, 
$$  \mathbb{A} =  \langle \{ 1, p, q, pq, r, pr, qr,  (pq)r \} \rangle_\mathbb{R},$$
which is the  real octonions generated by $\{ 1, p, q, pq, r, pr, qr,  (pq)r\}$.
\end{lem}

\bigskip

\noindent {\bf Proof.}  (i)  Let $\langle . , . \rangle$ be the symmetric bilinear form of  $\mathbb{A}$. For $ a \in \mathbb{A}$, use $n_a $ to denote $\langle a , a \rangle $. Note that if $ a = t_01 + t_1 p +  t_2q + t_3 pq \in \mathbb{A}$ is nonzero, then $ a^{-1} = n_a^{-1} ( 2 \langle a , 1 \rangle 1 - a)  = n_a^{-1}( t_0 1 -  t_1 p -  t_2q - t_3 pq) $. 
 Pick an $ r_0 \in  \mathbb{A} \setminus \langle  \{ 1, p, q, pq \} \rangle_F$ and let 
  $ r := r_0 - \frac{\langle r_0 , 1 \rangle}{\langle 1 , 1 \rangle} 1  -  \frac{\langle r_0 , p \rangle}{\langle p , p \rangle} p -\frac{\langle r_0 , q \rangle}{\langle q , q \rangle} q - \frac{\langle r_0 , pq \rangle}{\langle pq , pq \rangle} pq $. Clearly, $ r \in  \mathbb{A} \setminus \langle \{ 1, p, q, pq \} \rangle_F$. This, just as we saw in the proof of part (i) of the preceding lemma,   reveals that 
\begin{equation*}
r^2  =  -n_r  1  , \  (pr)^2   =  - n_p  n_r  1 , \ (qr)^2   =   - n_q  n_r 1  , \  ((pq)r)^2   =  -  n_p  n_q n_r  1 , 
\end{equation*}
\begin{equation*}
   pr  =  - rp , \ qr  =   -rq  , \ (pq)r   =  - r (pq).
\end{equation*}
 That $\{ 1, p, q, pq, r, pr, qr,  (pq)r\}$ is linearly independent over $F$ is easy. If $  t_01 + t_1 p +  t_2q + t_3 pq + t_0' r + t_1' pr + t_2' qr + t_3' (pq)r = 0$, where $t_i, t'_i \in F$ ($ 0 \leq i \leq 3$), then $( t_0', t_1' , t_2' , t_3' ) =  0$, for otherwise,   with $ a' := t_0' 1 + t_1' p + t_2' q + t_3' pq$, we have 
  \begin{eqnarray*}
  r & = & n_{a'}^{-1} ( t'_0 1 -  t'_1 p -  t'_2q - t'_3 pq) ( t_01 + t_1 p +  t_2q + t_3 pq) \\
  & \in & \langle \{ 1, p, q, pq \} \rangle_F,
  \end{eqnarray*}
 which is impossible. Thus, $( t_0', t_1' , t_2' , t_3' ) =  0$, and hence so is $( t_0, t_1 , t_2 , t_3 ) $, proving that $\{ 1, p, q, pq, r, pr, qr,  (pq)r\}$ is linearly independent. 
 
 Next, by obtaining the multiplication table of the generating set 
 \[\{ 1, p, q, pq, r, pr, qr,  (pq)r\},\]
 we see that $ \mathbb{A}$ contains a copy of the  octonion algebra generated by $ \{ 1, p, q, pq, r, pr, qr,  (pq)r\}$.  First, the following equalities are easily checked. 
 \begin{eqnarray*}
r (pr) & =&  n_r  p  = -  (pr)r;\\
  r (qr) & =&    n_rq= - (qr)r   ; \\
r ((pq)r) & =&   n_rpq = -   ((pq)r)r ;\\
 (pr) p&  =&  n_p r = -  p (pr); \\
   (qr)q&  =& n_q  r  = -q (qr) ; \\
  ((pq)r) (pq)&  =& n_p n_q r = - (pq) ((pq)r).
  \end{eqnarray*}
 Next, making use of the second set of identities in Proposition \ref{1.3}(i), we can write  
   \begin{eqnarray*}
   p (qr) & =& - p (rq) = r(pq)  = -(pq)r        ,\\
(qr) p  &  =& - (qp)r  = (pq)r ; \\
 q (pr) & =&      -  q (rp)= r(qp)  =(pq)r  , \\  
     (pr) q & = &  -  (pq)r  ; \\
  q ( (pq)r) & = & -q (r(pq)) = r (q(pq))=  -n_q  pr        , \\   
 ((pq)r)q  & =&   -  ((pq)q)r=  n_q pr ;\\
   p((pq)r) & =& - p(r(pq))  =  r(p(pq))= n_p qr    , \\
    ((pq)r)p & =&   -  ((pq)p)r  =     ((qp)p)r =  - n_p qr  ; \\
  (pr) (qr) & =&  -  (q  (pr)r)   =-n_r pq   ,\\    
 (qr) (pr)&  =&   -  (r  (qp)r)  = n_r  pq ;  \\
 (pq) (pr) & =&   - (pq) (rp) =  - (r ((pq) p))  =- n_p qr      , \\
(pr) (pq) & =& -    (p( pq)r) = n_p qr  ; \\
  (pq) (qr) & =&    -(pq) (rq)=  r ((p q) q)  =n_q pr  , \\
 (qr) (pq) & =&  -  ( q (pq)) r ) = -n_q pr ;\\                   
 (pr) ((pq)r) & =&  - ( (pq) ((pr)r)) = (pq)p = n_p n_r q , \\ 
   ((pq)r)   (pr)  & = &  - ( (pq)(pr) )r=  (qr) r = -n_p n_r q   ; \\ 
 (qr) ((pq)r) & =&  - ( (pq) ((qr)r)) =(pq)q= -n_q n_r p    , \\ 
 ((pq)r)   (qr)  & = &  - ( (pq))   (qr) )r= - (pr) r= n_r n_q  p      . 
\end{eqnarray*}

So far we have shown that  $ \mathbb{A}$ contains a copy of the octonion  algebra  generated by the set $\{ 1, p, q, pq, r, pr, qr,  (pq)r\}$, whose  multiplication table is the same as that of the generator set of an octonion algebra.  In the following table, we have 
\begin{equation*}
n_{e_1} = n_p,  \ n_{e_2} = n_q, \ n_{e_3} = n_p n_q , \ \ n_{e_4} = n_r , 
\end{equation*}
\begin{equation*}
 n_{e_5} = n_p n_r, \  n_{e_6} = n_q n_r, \ n_{e_7} =n_p n_q n_r.
\end{equation*}
\vskip 2em   
{\vbox{   
\begin{center}   
{\small   
\begin{tabular}{|c|c|c|c|c|c|c|c|c|c|c|}                    \hline   
     &  $e_0 =1$ & $e_1= p$ & $e_2 =q$ & $e_3 = pq$ & $e_4= r$  & $e_5= pr$ & $e_6= qr$ & $e_7 = (pq)r$ \\ \hline   
  $e_0$ &  $1$ & $e_1$  & $e_2$ & $e_3$ & $e_4$ & $e_5$ & $e_6$ & $e_7$ \\ \hline   
 $e_1$ &  $e_1$ & $-n_{e_1} $ & $ e_3$ & $-n_{e_1} e_2$ & $e_5$ & $-n_{e_1}e_4$ & $-e_7$ & $n_{e_1}e_6$     \\ \hline   
 $e_2$ & $e_2$ & $-e_3$ & $-n_{e_2}$ & $n_{e_2} e_1$ & $e_6$ & $e_7$ & $-n_{e_2} e_4$ & $-n_{e_2}e_5$   \\ \hline   
 $e_3$ & $e_3$ & $n_{e_1}e_2$ & $-n_{e_2}e_1$ & $-n_{e_3}$ & $e_7$ & $-n_{e_1}e_6$ & $n_{e_2}e_5$  &   $-n_{e_3}e_4$ \\ \hline   
  $e_4$ & $e_4$ & $-e_5$ & $-e_6$ & $-e_7$ & $-n_{e_4}$ & $n_{e_4}e_1$ & $n_{e_4}e_2$  &  $n_{e_4}e_3$  \\ \hline   
  $e_5$ & $e_5$ & $n_{e_1}e_4$ & $-e_7$ & $n_{e_1}e_6$ & $-n_{e_4}e_1$ & $-n_{e_5}$ & $-n_{e_4}e_3$   & $n_{e_5}e_2$  \\ \hline   
  $e_6$ & $e_6$ & $e_7$ & $n_{e_2}e_4$ & $-n_{e_2}e_5$ & $-n_{e_4}e_2$ & $n_{e_4}e_3$ & $-n_{e_6}$   &  $-n_{e_6}e_1$ \\  \hline   
 $e_7$ & $e_7$ & $-n_{e_1}e_6$ & $n_{e_2}e_5$ & $n_{e_3}e_4$ & $-n_{e_4}e_3$ & $-n_{e_5}e_2$ & $n_{e_6}e_1$   & $-n_{e_7}$ \\  \hline   
\end{tabular}} 
\end{center}   
\vskip 1em 
}}   
\vskip 1em

Finally, renaming the generators $1 , p, q, \ldots$'s as in the above table, we complete the proof by showing  that   $    \langle  \{e_i\}_{i =0}^7   \rangle_F =  \mathbb{A}$.  Suppose on the contrary that 
 $   \langle  \{e_i\}_{i =0}^7  \rangle_F     \not=  \mathbb{A} $ and pick an $s_0 \notin      \langle  \{e_i\}_{i =0}^7 \rangle_F $. Let 
 \begin{equation*}
  e_8 := s := s_0 -  \sum_{i=0}^7 \frac{\langle s_0 , e_i \rangle}{\langle e_i, e_i \rangle} e_i 
  \end{equation*}
 Clearly, $e_8  \notin  \langle  \{e_i\}_{i =0}^7 \rangle_F $.  Once again, this reveals that 
\begin{equation*}
e_8^2   =  -n_{e_8} 1 , \  e_8 e_i = - e_i e_8  , (1 \leq i \leq 7).
\end{equation*}
 Now, we obtain a contradiction by showing that  $(pq) (rs)  = (rs)(pq)  =  0$. To this end, first, in view of Theorems \ref{02.3}(ii) and \ref{02.5}(i), and the second remark following Theorem \ref{02.5}, as well as  the second set of identities in Proposition \ref{1.3}(i), we can write 
 \begin{eqnarray*}
 \langle rs ,  1  \rangle & = & -  \langle s ,  r^*  \rangle =  \langle s ,  r  \rangle = 0,
\end{eqnarray*}
\begin{eqnarray*}
\langle (qr)s ,  p  \rangle  & = & \langle  s ,  (qr)^*p  \rangle = \langle   s,  (r^* q^*) p  \rangle  =  \langle  s,  (r q) p  \rangle , \\
&  =&  - \langle   s,  (q r) p \rangle =  \langle s,  (q p) r  \rangle = -  \langle s,  ( pq) r  \rangle , \\
&  =&  0, 
\end{eqnarray*}
and
\begin{eqnarray*}
\langle rs ,  pq  \rangle  & = & \langle  s ,  r^* (pq)  \rangle = - \langle   s,  r( pq)  \rangle  =  \langle  s,  (p q) r  \rangle , \\
&  =&  0, 
\end{eqnarray*}
from which, we obtain $  [  (qr) s] p = - p [  (qr) s] $ and $  (rs) (p q) =  -  (pq) (rs)$, respectively.  
 Now, once again in view of the second set of identities in Proposition \ref{1.3}(i), we have 
\begin{eqnarray*}
   (pq) (rs)    & = & -  [  (pq) r] s   = [  (qp) r] s  = - [  (qr) p] s  , \\  
  & = &  [  (qr) s] p =  - p [  (qr) s] =   p [  (rq) s]  , \\
  & = &  -  p [  (rs) q]  =   (rs) (p q) , \\
  & = &   -  (pq) (rs) ,
\end{eqnarray*}
implying that $ (pq) (rs)  = (rs)(pq)  =  0$, a contradiction; compare this proof with Albert's \cite[page 766]{A3}. Therefore, $     \mathbb{A} = \langle  \{e_i\}_{i =0}^7   \rangle_F$, as claimed. 
This completes the proof.

(ii) As noted in part (ii) of the preceding lemma,  $\mathbb A$ is alternative and locally complex equipped by an inner-product  $\langle . , . \rangle$, whose induced norm, denoted by $ |.| $,  is indeed an absolute value on $\mathbb A$. Just as in part (ii) of the preceding lemma, the proof is almost identical to that of part (i) except that
for a picked  element 
$ r_0 \in  \mathbb{A} \setminus \langle \{ 1, p, q, pq \} \rangle_{\mathbb{R}}$, the norm-one element $ r \in  \mathbb{A} \setminus \langle \{ 1, p, q, pq \} \rangle_\mathbb{R}$ is defined by  $ r := \frac{1}{|\hat r|} \hat r$, where 
 $ \hat r := r_0 - \langle r_0 , 1 \rangle 1  -  \langle r_0 , p \rangle p - \langle r_0 , q \rangle  q -  \langle r_0 , pq \rangle pq $.
 It would then follow that 
\begin{equation*}
r^2  =    (pr)^2   =   (qr)^2   =    ((pq)r)^2   = - 1 , 
\end{equation*}
\begin{equation*}
   pr  =  - rp , \ qr  =   -rq  , \ (pq)r   =  - r (pq).
\end{equation*}
Just as we saw in part (i),  the multiplication table of the generating set 
 \[G : =\{ 1, p, q, pq, r, pr, qr,  (pq)r\},\]
is seen to be according to the following table, and hence $ \mathbb{A}$ contains a copy of the algebra of real octonions generated by the set $G$.

\vskip 2em   
{\vbox{   
\begin{center}   
{\small   
\begin{tabular}{|c|c|c|c|c|c|c|c|c|c|c|}                    \hline   
     &  $e_0 =1$ & $e_1= p$ & $e_2 =q$ & $e_3 = pq$ & $e_4= r$  & $e_5= pr$ & $e_6= qr$ & $e_7 = (pq)r$ \\ \hline   
  $e_0$ &  $1$ & $e_1$  & $e_2$ & $e_3$ & $e_4$ & $e_5$ & $e_6$ & $e_7$ \\ \hline   
 $e_1$ &  $e_1$ & $-1 $ & $ e_3$ & $-e_2$ & $e_5$ & $-e_4$ & $-e_7$ & $e_6$     \\ \hline   
 $e_2$ & $e_2$ & $-e_3$ & $-1$ & $e_1$ & $e_6$ & $e_7$ & $-e_4$ & $-e_5$   \\ \hline   
 $e_3$ & $e_3$ & $e_2$ & $-e_1$ & $-1$ & $e_7$ & $-e_6$ & $e_5$  &   $-e_4$ \\ \hline   
  $e_4$ & $e_4$ & $-e_5$ & $-e_6$ & $-e_7$ & $-1$ & $e_1$ & $e_2$  &  $e_3$  \\ \hline   
  $e_5$ & $e_5$ & $e_4$ & $-e_7$ & $e_6$ & $-e_1$ & $-1$ & $-e_3$   & $e_2$  \\ \hline   
  $e_6$ & $e_6$ & $e_7$ & $e_4$ & $-e_5$ & $-e_2$ & $e_3$ & $-1$   &  $-e_1$ \\  \hline   
 $e_7$ & $e_7$ & $-e_6$ & $e_5$ & $e_4$ & $-e_3$ & $-e_2$ & $e_1$   & $-1$ \\  \hline   
\end{tabular}} 
\end{center}   
\vskip 1em 
}}   
\vskip 1em   
 
For the rest of the proof, the reader can adjust the remaining part of the proof of part (i) or refer to the rest of the proof given in detail in Version 2 of this paper on arxiv.org. This completes the proof. 
\hfill \qed 

\bigskip 

\noindent {\bf Remark.} 
It is worth mentioning  that  if, in parts (ii) of Lemmas \ref{2.1},  \ref{2.3},  and \ref{2.4}, we replace the hypotheses that the algebra $ \mathbb{A}$  is left (resp. right) alternative and that it  has no nontrivial joint zero-divisors by the stronger hypotheses that $ \mathbb{A}$ is alternative and  that it has no nontrivial zero-divisors, then the proofs become more self-contained, and hence more accessible. That is  because by Proposition \ref{1.1}(i), whose proof is quite elementary, the real algebra  $\mathbb A$ would be quadratic, and hence locally complex. 

\bigskip

The  following is a consequence of  part (ii) of Lemmas  \ref{2.3}  and \ref{2.4} and the Skolem-Noether Theorem. 

\begin{cor} \label{3.5} 
 {\rm (i)} Let $ p \in \mathbb{H}$ be such that $ p^2 = -1$. Then, there exists a $ q  \in \mathbb{H}$ with $ q^2 = -1$ such that $ pq = -qp$ and that $ \{ 1, p , q , pq\}$ is a generator set of quaternions. Moreover, corresponding to  any generator set $ \{ 1, I, J , IJ\}$ of quaternions with $ I^2 = J^2 = (IJ)^2 = -1$, there is a $ u \in \mathbb{H}$ with $ |u| = 1$ such that $ \bar u p u = I$, $ \bar u q u = J$, and $ \bar u (pq) u = IJ$. 

{\rm (ii)}  Let $ p_1, p_2  \in \mathbb{O}$ be such that $ p_1^2 = p_2^2 = p_3^2 = -1$, where $ p_3 = p_1 p_2$. Then, there exists a $ p_4 \in   \mathbb{O}$ with $ p_4^2 = -1$ such that $\{  p_i\}_{i=0}^7$ is a generator set for $\mathbb{O}$, where $p_0 =1$, $ p_5 =p_1 p_4$,  $ p_6 =p_2 p_4$, and $ p_7 =p_3 p_4$. 
\end{cor}

\bigskip

\noindent {\bf Proof.} The assertions, whose proofs are omitted for the sake of brevity, are straightforward consequences of the proofs of part (ii) of Lemmas \ref{2.3} and \ref{2.4} and the Skolem-Noether Theorem, \cite[page 39]{D}. 
\hfill \qed 

\bigskip 

In the following theorem,  we revisit and slightly strengthen  Frobenius' Theorem and  present slight extensions of its topological counterparts in several settings; part (iii)-(a) of the theorem must be known by the experts. See  \cite[Theorem 3.1]{K2},  \cite[Theorem 2.5.40]{GP}, \cite[Theorem 4]{A},  \cite[Theorem 2]{A2}, \cite[page 812]{K1}, and \cite[Theorem 3.1]{K2} . It is worth saying that in the next two theorems, as strictly stated  in their statements,  when the  algebras are assumed to be  equipped with norms, we mean vector space norms as opposed  to algebra norms. 

\bigskip 

\begin{thm} \label{2.5} 
 {\rm (i)}   {\bf (Frobenius)}  Let $ \mathbb{A}$  be an algebraic associative real algebra with  no nontrivial joint zero-divisors. Then, the algebra $ \mathbb{A}$ is locally complex and, equipped with its natural norm, is isometrically isomorphic to  $\mathbb{R}$, $\mathbb{C}$, or  $\mathbb{H}$.

\bigskip

{\rm (ii)}  Let $ \mathbb{A}$  be an associative real algebra. Then, $ \mathbb{A}$ is is locally complex and, equipped with its natural norm, is isometrically isomorphic to $\mathbb{R}$, $\mathbb{C}$, or $\mathbb{H}$ if   one of the following conditions holds.

{\rm (a)} {\bf (Kaplansky)}   $\mathbb A$ is a   real algebra  endowed with  a vector space norm relative to which the multiplication of it is separately continuous and   with  ${\rm j.t.z.d.}(\mathbb{A}) = \{ 0\}$.

{\rm (b)} {\bf (Gelfand-Mazur)}  $\mathbb A$ is a  quasi-division  real algebra endowed with  a vector space norm relative to which the multiplication of it is separately continuous.

 {\rm (c)}   {\bf (Arens)} $\mathbb A$ is a  real algebra  endowed with a $c$-supermultiplicative vector space norm for some $ c > 0$ relative to which the  multiplication is separately continuous, e.g., a real nearly absolute value.

  {\rm (d)} $\mathbb A$ is a  topological quasi-division real algebra, equivalently, a classical division algebra,  whose (topological) dual separates its points.

\bigskip

{\rm (iii)}  Let $ \mathbb{A}$  be an associative  real algebra. Then, $ \mathbb{A}$ is isometrically isomorphic to $\mathbb{R}$, $\mathbb{C}$, or $\mathbb{H}$ if  one of the following conditions holds, in which case the given norm of $ \mathbb{A}$ coincides with its natural norm that comes  from its local complex structure.

{\rm (a)} $\mathbb A$ is a real ${\rm C}^*$-algebra with no nontrivial joint topological zero-divisors or a quasi-division real ${\rm C}^*$-algebra.

{\rm (b)} $\mathbb A$ is a quasi-division  real algebra equipped with a vector space norm $\|.\|$ relative to which the multiplication of it is separately continuous and satisfying the identity $\|x^k \| = \|x\|^k $ on $ \mathbb{A}$ for a fixed positive integer $ k > 1$.

 {\rm (c)} $\mathbb A$ is a   real algebra equipped with a vector space norm $\|.\|$ relative to which the multiplication of it is separately continuous and satisfying the identity $\|x^k \| = \|x\|^k $ on $ \mathbb{A}$ for a fixed positive integer $ k > 1$ and  with no nontrivial joint topological zero-divisors.

 {\rm (d)} {\bf (Albert)} $\mathbb A$ is an absolute-valued real algebra.
\end{thm}

\bigskip

\noindent {\bf Proof.} (i)   By Proposition \ref{1.1}(i), the algebra $ \mathbb{A}$ is unital and hence contains a copy of real numbers and is quadratic, and more specifically locally complex, for it is algebraic and it has no nontrivial joint zero-divisors. If the algebra $ \mathbb{A}$ were not equal to its copy of real numbers, it would then contain a copy of complex numbers. Once again, if the algebra  $ \mathbb{A}$ were not equal to its copy of complex numbers, it would then contain a copy of real quaternions by Lemma \ref{2.3}(ii). We now see that $ \mathbb{A}$ is equal to its copy of the real quaternions, and hence is isomorphic to $\mathbb{H}$. That is because if $ \mathbb{A}$ were not equal to its copy of the real quaternions, by  Lemma \ref{2.4}(ii), it would then contain a copy of real octonions, contradicting the associativity of $ \mathbb{A}$.  That $ \mathbb{A}$, equipped with its natural norm, is isometrically isomorphic to $\mathbb{R}$, $\mathbb{C}$, or $\mathbb{H}$ follows from Theorem \ref{2.13}. This completes the proof.

(ii) By Theorem \ref{1.9}, the algebra $ \mathbb{A}$ is locally complex. The assertion now follows from (i).

(iii) By (ii) the algebra $ \mathbb{A}$ is isomorphic to $\mathbb{R}$, $\mathbb{C}$, or $\mathbb{H}$. That the algebra $ \mathbb{A}$ is indeed isometric to $\mathbb{R}$, $\mathbb{C}$, or $\mathbb{H}$ follows from the uniqueness of the norm of real ${\rm C}^*$-algebras if (a) holds. If (b), (c), or (d) holds, the assertion follows from Theorem \ref{2.13},  completing the proof. 
\hfill \qed 

\bigskip

The following theorem revisits and slightly strengthens Zorn's extension of Frobenius' Theorem and gives slight extensions of its topological counterparts in several settings; part (iii)-(a) of the theorem is probably known by the experts. See \cite[Theorem 2.5.29]{GP}, 
\cite[Theorem 2.5.50]{GP}, \cite[Corollary 2.5.51]{GP}, \cite[Corollary 2.5.57]{GP}, and  \cite[Theorem 2]{A2}.

\bigskip 

\begin{thm} \label{2.6} 
 {\rm (i)} {\bf (Zorn)}   Let $ \mathbb{A}$  be an algebraic left (resp. right) alternative real algebra with  no nontrivial joint zero-divisors. Then, the algebra $ \mathbb{A}$ is locally complex and, equipped with its natural norm, is isometrically  isomorphic to $\mathbb{R}$, $\mathbb{C}$,  $\mathbb{H}$, or $\mathbb{O}$.

\bigskip

{\rm (ii)}  Let $ \mathbb{A}$  be a left (resp. right) alternative real algebra.  If one of the following conditions holds, then the algebra $ \mathbb{A}$  is locally complex and, equipped with its natural norm, is isometrically isomorphic to $\mathbb{R}$, $\mathbb{C}$,  $\mathbb{H}$,  or $\mathbb{O}$.

{\rm (a)} {\bf (El-Mallah and Micali)} $\mathbb A$ is a  real algebra  endowed with  a vector space norm relative to which the multiplication of it is separately continuous  and 
 with  ${\rm j.t.z.d.}(\mathbb{A}) = \{ 0\}$.

{\rm (b)} {\bf (Nieto)}  $\mathbb A$ is a  quasi-division real algebra  endowed with  a vector space norm relative to which the multiplication of it is separately continuous.

 {\rm (c)} $\mathbb A$ is a   real algebra  endowed with a $c$-supermultiplicative vector space norm for some $ c > 0$ relative to which the  multiplication is separately continuous, e.g., a real nearly absolute-value.

 {\rm (d)}     $\mathbb A$ is a  topological quasi-division real algebra, equivalently, a classical division algebra,  whose (topological) dual separates its points.

 \bigskip

{\rm (iii)}  Let $ \mathbb{A}$  be a left (resp. right) alternative  real algebra. Then, $ \mathbb{A}$ is isometrically isomorphic to $\mathbb{R}$, $\mathbb{C}$, $\mathbb{H}$,  or $\mathbb{O}$  if  one of the following conditions holds, in which case the given norm of $ \mathbb{A}$ coincides with its natural norm coming from its local complex structure.

{\rm (a)}  $\mathbb A$ is a    real ${\rm C}^*$-algebra with no nontrivial joint topological zero-divisors or a quasi-division real ${\rm C}^*$-algebra.

{\rm (b)}   $\mathbb A$ is a  real algebra  endowed with  a vector space norm $ \| .\|$ relative to which the multiplication of it is separately continuous satisfying the identity $\|x^k \| = \|x\|^k $ on $ \mathbb{A}$ for a fixed positive integer $ k > 1$ and with  no nontrivial joint topological zero-divisors.

 {\rm (c)} $\mathbb A$ is a  quasi-division  real algebra  endowed with  a vector space norm $ \| .\|$ relative to which the multiplication of it is separately continuous and  it satisfies the identity $\|x^k \| = \|x\|^k $ on $ \mathbb{A}$ for a fixed positive integer $ k > 1$.

 {\rm (d)} {\bf (Albert)} $\mathbb A$ is an absolute-valued real algebra.
\end{thm}

\bigskip

\noindent {\bf Proof.} (i) By Proposition \ref{1.1}(iii), the algebra $ \mathbb{A}$ is unital and hence contains a copy of real numbers and is quadratic, for it is algebraic and  it has no nontrivial joint zero-divisors. So, by \cite[Theorem 1]{A2}, we get that $ \mathbb{A}$  alternative. It thus follows that $ \mathbb{A}$ is locally complex. Now, in light of Lemmas \ref{2.3}(ii) and \ref{2.4}(ii) and Theorem \ref{2.13}, the proof is almost identical to that of part (i) of the preceding theorem, which we omit for brevity.

(ii)-(iii) With (i) at our disposal, (ii) and (iii) follow from Theorem \ref{1.9} and from  the uniqueness of the norms of real ${\rm C}^*$-algebras along with Theorem \ref{1.11}. This completes the proof.
\hfill \qed

\bigskip

Part (ii) of the following theorem is the well-known non-commutative Urbanik-Wright theorem; see \cite[Theorem 2.6.21]{GP}, \cite[page 216]{GP}, \cite[Characterization 2]{Za}, and \cite[Theorem 3.1]{I3}.

\bigskip

\begin{thm} \label{2.7}
{\rm (i)} Let $\mathbb A$ be a  unital  real ${\rm C}^*$-algebra with no nontrivial joint topological zero-divisors or a unital quasi-division real ${\rm C}^*$-algebra. Then, $ \mathbb{A}$ is isometrically isomorphic to $\mathbb{R}$, $\mathbb{C}$, $\mathbb{H}$,  or $\mathbb{O}$.

{\rm (ii)} {\bf (Urbanik and Wright)}  Let $ ( \mathbb{A}, | .|)$ be a unital absolute-valued real algebra. Then, $ \mathbb{A}$ is isometrically isomorphic to $\mathbb{R}$, $\mathbb{C}$, $\mathbb{H}$,  or $\mathbb{O}$.

{\rm (iii)}  Let $ ( \mathbb{A}, | .|)$ be a flexible  absolute-valued real algebra with a left (resp. right) identity element. Then, $ \mathbb{A}$ is isometrically isomorphic to $\mathbb{R}$, $\mathbb{C}$, $\mathbb{H}$,  or $\mathbb{O}$. 

{\rm (iv)}   {\bf (Zalar)}  Let $\mathbb A$ be a  left (resp. right) alternative real  algebra endowed with an inner-product $ \langle . , . \rangle$, whose norm satisfies the identity $\|a^2\| = \|a\|^2 $ on $\mathbb A$. Then, $\mathbb A$  is isometrically isomorphic to $\mathbb{R}$, $\mathbb{C}$, $\mathbb{H}$,  or $\mathbb{O}$.

{\rm (v)} {\bf (Ingelstam)} Let $\mathbb A$ be a  left (resp. right) alternative real  unital algebra endowed with an inner-product $ \langle . , . \rangle$, whose norm is a unital algebra norm on $\mathbb A$. Then, $\mathbb A$ isometrically isomorphic to $\mathbb{R}$, $\mathbb{C}$, $\mathbb{H}$,  or $\mathbb{O}$.

\end{thm}

\bigskip

\noindent {\bf Proof.}  (i)  By \cite[Theorem 3.2.5]{GP}, the algebra  $\mathbb A$ is alternative. So the assertion follows from part (iii)-(a) of the preceding theorem.

(ii) The assertion is a quick consequence of Albert's Theorem, namely, part (iii)-(d) of the preceding theorem, and Theorem \ref{1.10}(i). 

(iii) By Theorem \ref{1.10}(iii), the algebra  $\mathbb A$  is unital and alternative. So Albert's Theorem, namely,  part (iii)-(d) of the  preceding theorem, applies, proving the assertion.

(iv) By Theorem \ref{1.9}(v) and its proof, the algebra $ \mathbb{A}$ is a locally complex and alternative division algebra, and hence an absolute-valued algebra by Theorem \ref{02.5}(ii).  Now, Albert's Theorem, namely,  part (iii)-(d) of  the preceding theorem, applies, proving the assertion.  

(v) By Theorem \ref{1.9}(vi), the algebra $ \mathbb{A}$ is a locally complex and alternative algebra, and hence an absolute-valued algebra by Theorem \ref{02.5}(ii).  Again, Albert's Theorem, namely,  part (iii)-(d) of  the preceding theorem, applies, proving the assertion.  
\hfill \qed

\bigskip

We conclude this note with a Hurwitz-type theorem; see \cite[Theorem 3.25]{Sch} and \cite{K3}. 

\bigskip

\begin{thm} \label{3.9}
  {\bf (Hurwitz)}  {\rm (i)}  Let $F$ be a field whose characteristic is not $2$,   $\mathbb{A}$  a left (resp. right) alternative  quadratic algebra over $F$ with the property that 
  the algebra generated by any nonzero element of $\mathbb{A}$ has no nontrivial zero-divisors. Then, $\mathbb{A}$ is isomorphic to one of the following algebras over $F$: the field $F$;  a field extension  of degree $2$ over $F$;  a quaternion algebra over $F$;  an octonion algebra over $F$.

{\rm (ii)}  Let $F$ be a field whose characteristic is not $2$ and   $\mathbb{A}$  a flexible   algebra over $F$ with a left (resp. right) identity element $e$ and with no nontrivial zero-divisors. Let  $\langle . , . \rangle $ be a nondegenerate bilinear form on $\mathbb{A}$ that permits composition. Then, $\mathbb{A}$ is isomorphic to one of the following algebras over $F$: the field $F$;  a field extension  of degree $2$ over $F$;  a quaternion algebra over $F$;  an octonion algebra over $F$. 

\end{thm}

\bigskip

\noindent {\bf Proof.} (i) In light of Lemmas \ref{2.3}(i)  and \ref{2.4}(i), the proof is almost identical to that of Theorem \ref{2.5}(i). 
 If the algebra $ \mathbb{A}$ were not equal to $F1$, it would then contain a field extension  of degree $2$ over $F$, say, $K$. Once again, if the algebra  $ \mathbb{A}$ were not equal to $K$, by Lemma \ref{2.3}(i), it would then contain a quaternion algebra over $F$, say, $Q$. Then again, if the algebra  $ \mathbb{A}$ were not equal to $Q$, by  Lemma \ref{2.4}(i), it would contain an octonion algebra  over $F$, say, $O$. We now see that $ \mathbb{A}$ is equal to $O$. That is because if $ \mathbb{A}$ were not equal to  $O$, by  Lemma \ref{2.4}(i), it would contain a nontrivial joint zero-divisor. This is a contradiction because, by Theorem \ref{02.5}(ii),  $ \mathbb{A}$ has no nontrivial zero-divisors. This completes the proof.

(ii) By Theorem \ref{02.5}(i), the algebra $\mathbb{A}$ is quadratic and alternative. Thus, (i) applies, completing the proof. 
\hfill \qed

\bigskip

\noindent {\bf Remark.}   In part (ii), if the algebra $\mathbb{A}$ is assumed to have an identity element, then the hypothesis that $\mathbb{A}$ is flexible is redundant and  the hypothesis that $\mathbb{A}$  has  no nontrivial zero-divisors can be replaced by the weaker hypothesis that the algebra generated by any nonzero element of $\mathbb{A}$ has no nontrivial zero-divisors.

\bigskip

\noindent {\bf Acknowledgement.}  I would like to thank Professor Heydar Radjavi for reading the manuscript and making helpful comments. I also thank Claude, an AI assistant developed by Anthropic, for a careful reading of the sixth version of this manuscript, which led to the correction of errors in Section 2, the proof of Theorem \ref{02.8}, and in the topological part,  including a  comment made in the first paragraph on page 7 of the sixth version that affected the statement of Corollary \ref{3.2}, as well as to improvements in the presentation; see Proposition \ref{02.4} and Theorem \ref{02.7}. Any remaining errors are mine.

\end{document}